\documentclass[pdflatex,sn-mathphys]{sn-jnl}
\jyear{2021}
\usepackage{rotating}
\usepackage{tabularx}
\usepackage{xcolor}
\usepackage{multicol}
\usepackage{lipsum}
\usepackage{booktabs}
\usepackage{stackengine} 
\newcommand\oast{\stackMath\mathbin{\stackinset{c}{0ex}{c}{0ex}{\ast}{\bigcirc}}}
\usepackage{subfiles}
\usepackage{amsmath,amssymb,amsfonts,stackrel}
\usepackage{graphicx}
\usepackage{textcomp}
\usepackage{url,hyperref,lineno,microtype,subcaption}
\usepackage[onehalfspacing]{setspace}
\usepackage{relsize}
\usepackage{amssymb}
\usepackage{pifont}

\usepackage{lineno,hyperref}
\usepackage{amsfonts,epsfig}
\usepackage{multirow}

\usepackage[mathscr]{euscript} 
\usepackage{enumerate}
\usepackage{amsmath, amscd,amsfonts, amssymb}

\theoremstyle{thmstyleone}%
\theoremstyle{thmstyletwo}%

\usepackage{times}
\usepackage{extarrows}
\usepackage{stmaryrd}
\usepackage{textcomp}
\usepackage[mathscr]{euscript}
\usepackage{multirow}
\usepackage{soul}
\usepackage{amsthm}

\theoremstyle{definition}

\theoremstyle{definition}%
\newtheorem{defn}{Definition}%

\theoremstyle{thmstylethree}%

\begin{document}

\title[Cross Tensor Approximation for Image and Video
Completion]{Cross Tensor Approximation for Image and Video
Completion}

\author*[1]{\fnm{Salman} \sur{Ahmadi-Asl}}\email{s.asl@skoltech.ru}

\author[1]{\fnm{Maame} \sur{Gyamfua Asante-Mensa}}

\author[1,2]{\fnm{Andrzej} \sur{Cichocki}}

\author[1]{\fnm{Anh-Huy} \sur{Phan}}

\author[1]{\fnm{Ivan} \sur{Oseledets}}

\author[1]{\fnm{Jun} \sur{Wang}}

\affil*[1]{\orgdiv{Skolkovo Institute of Science and Technology (SKOLTECH)}, \orgname{Center for
Artificial Intelligence Technology,} 
\orgaddress{\state{Moscow}, \country{Russia}}}

\affil[2]{\orgdiv{Systems Research Institute of Polish Academy of Science}, \orgname{} 
\orgaddress{\state{Warsaw}, \country{Poland}}}


\abstract{This paper proposes a general framework to use the cross tensor approximation or tensor ColUmn-Row (CUR) approximation for reconstructing incomplete images and videos. {The key importance of the new algorithms is their simplicity and ease of implementation with low computational complexity. For the case of data tensors with 1) structural missing components or 2) a high missing rate, we propose an efficient smooth tensor CUR algorithms which first make the sampled fibers smooth and then apply the proposed CUR algorithms. The numerical experiments show the significant benefit of this smoothing procedure. The main contribution of this paper is to develop/investigate improved multistage CUR algorithms with filtering (smoothing ) preprocessing for tensor completion. The second contribution is a detailed comparison of the performance of image recovery for four different CUR strategies via extensive computer simulations. Our simulations clearly indicated that the proposed algorithms are much faster than most of the existing state-of-the-art algorithms developed for tensor completion, while performance is comparable and often even better.} Furthermore, we will provide in GitHub the MATLAB codes which can be used for various applications. Moreover, to our best knowledge, the CUR (cross approximation) algorithms have not been investigated nor compared till now for image and video completion.}




\maketitle

\section{Introduction}\label{sec1}

Tensors are efficient tools for multi-way data processing (analysis) as they can preserve the intrinsic structures of multidimensional data tensors, e.g., images and videos \cite{cichocki2016tensor,cichocki2017tensor,cichocki2015tensor,cichocki2009nonnegative}. In many real applications, the underlying data tensors contain missing components or are corrupted by noise/outliers due to inaccurate data acquisition processes or corruption by artifacts.
The process of recovering (reconstructing) an incomplete data tensor from its partially observed data tensor is called tensor completion. In the past few decades, many algorithms have been developed to solve this problem \cite{long2019low,song2019tensor,asante2021matrix} due to its  importance in several applications such as in recommender systems \cite{frolov2017tensor}, computer vision \cite{liu2012tensor}, chemometrics \cite{tomasi2005parafac}, link prediction \cite{acar2009link}, etc. Tensor completion algorithms are generally categorized into 1) rank minimization and 2) tensor decomposition techniques in which the concepts of tensor rank and tensor decomposition play key roles in the analysis. It is worth mentioning that in contrast to the matrices, the notion of the rank for tensors is not unique, and different tensor ranks can be defined. There are several types of tensor decompositions such as CANDECOMP/PARAFAC decomposition \cite{hitchcock1928multiple,hitchcock1927expression}, Tucker decomposition \cite{tucker1963implications,tucker1964extension,tucker1966some} and its special case, i.e., Higher Order SVD (HOSVD) \cite{de2000multilinear}, Block Term decomposition \cite{de2008decompositionsI,de2008decompositionsII,de2008decompositionsIII}, Tensor Train/Tensor Ring (Chain)  (TT-TR(TC)) decomposition \cite{oseledets2011tensor,zhao2016tensor,espig2012note}, tubal SVD (t-SVD) \cite{kilmer2011factorization,kilmer2013third,braman2010third}, each of which generalizes the notion of the matrix rank to tensors in efficient ways. Each iteration of the completion algorithms usually needs a low-rank tensor approximation of the underlying data tensors. The deterministic algorithms are prohibitive for these computations, especially when many iterations are required for convergence or the underlying incompleted data tensor is very large. These drawbacks motivate developing fast algorithms for the low-rank tensor approximation \cite{ahmadi2021randomized,ahmadiTRrandomized2020}.

The ColUmn-Row (CUR) or cross-skeleton algorithms have been introduced in the numerical multilinear/linear algebra community mainly for fast low-rank matrix/tensor approximation of large-scale data matrices/tensors. 

The main motivations for utilizing the CUR algorithms are 1- fast low tensor rank approximation, 2- higher compression ratio, and 3- data interpretation issues. For example, in the case of matrices (second-order tensors), the SVD of sparse matrices does not provide sparse factor matrices while the cross-skeleton approximation achieves this goal, and this leads to more compact data representation. This paper intends to benefit from the first property of the CUR approaches and develop fast tensor completion algorithms.

Although the CUR algorithms have been extensively utilized for the low-rank matrix/tensor approximation and compression purposes, here we use them for the data completion task. Similar algorithms to matrix completion and tensor completion using the CUR approximation techniques have been studied in \cite{xu2015cur} and \cite{wang2017missing}. Our work differs from them in several aspects. First, our proposed algorithms are simple and easy to be implemented. Moreover, our work proposes a general framework that covers a variety of tensor/matrix CUR approximation techniques e.g., Tucker decomposition, tubal decomposition as special cases. 


We should highlight that the main advantage of our proposed algorithms is their simplicity, where they can be implemented in a few lines of code. They also have lower computational complexity than other completion algorithms because of utilizing CUR algorithms.



Our main contributions in this paper are:
\begin{itemize}
    \item Proposing a general framework of tensor CUR approximation, e.g., Tucker decomposition, tubal decomposition, and tensor CUR approximation, for the tensor completion task. Besides, we extensively compare the performance of the proposed tensor CUR models in the simulation part.
    
    \item Proposing smooth tensor CUR completion algorithms for reconstructing incomplete data tensors with structural missing patterns or a high missing ratio. The quality of different smoothing techniques are also compared in our experiments.
    
    \item {Extensive simulations on different data tensors such as images and videos. The simulation results indicate that the proposed algorithms are very fast and outperform most the state-of-the-art completion algorithms.} 
    
\end{itemize}

This paper is organized as follows. In Section \ref{Sec:Pre}, preliminary concepts and definitions are given. Section \ref{Sec:TCUR} is devoted to introducing different types of matrix and tensor CUR algorithms. In Section \ref{Sec:TCProb}, the tensor completion problem is described, and in Section \ref{Sec:Related} we review the literature on tensor completion algorithms. The proposed algorithm is discussed in Section \ref{Sec:ProApp} and its computational complexity is studied in Section \ref{Sec:CCom}. In Section \ref{Sec:Sim}, extensive simulations are conducted to verify the efficiency and applicability of the proposed algorithm. Finally, a conclusion is given in Section \ref{Sec:Conclu}. 
\section{Preliminary concepts and definitions}\label{Sec:Pre}
In this section, we present basic notations and concepts used throughout the paper. {We follow the notations and concepts used in \cite{cichocki2016tensor,cichocki2009nonnegative} in our presentations.} Tensors, matrices, and vectors are denoted by underlined bold upper case letters, e.g., $\underline{\bf X}$, bold upper case letters, e.g., ${\bf X}$ and bold lower case letters, e.g., ${\bf x}$, respectively. Fibers are first-order tensors produced by fixing all modes except one, while slices are second-order tensors generated by fixing all modes except two of them. For a third-order tensor $\underline{\mathbf X}$, the slices $\underline{\mathbf X}(:,:,k),\,\underline{\mathbf X}(:,j,:),\,\underline{\mathbf X}(i,:,:)$ are called frontal, lateral and horizontal slices, respectively. Fibers ${\mathbf X}(i,:,j),\,{\mathbf X}(:,j,k)$ and ${\mathbf X}(i,j,:)$ are called rows, columns and tubes, respectively \cite{cichocki2016tensor}. The generalizations of these concepts to higher order tensors are straightforward. For example, for an $N$th-order tensor, $N$ types of fibers can be defined and they are referred to $n$-mode fibers for $n=1,2,\ldots,N$. Note that the notations $\mathbf{X}^+$ and $\mathbf{X}^T$ denote the Moore-Penrose pseudoinverse (MP)\footnote{The notion of MP and transpose is also defined for tensors based on the t-product, see Appendix A.} and the transpose of matrix $\mathbf{X}$, respectively. The Frobenius norm and the Hadamard (element-wise) product of tensors/matrices are denoted by ${\left\| . \right\|_F}$ and $\oast$, respectively. We use the notation ${{\bf P}_{\underline{\bf\Omega}} }({\underline{\bf X}})$ to denote the projection operator which keeps fixed the components of the tensor $\underline{\bf X}$ for the indices $\underline{\bf \Omega}$ and zeros-out the others. ${{\underline{\bf \Omega} ^ \bot }}$ denotes the complement of the set $\underline{\bf \Omega}$. 

Given an $N$th-order tensor $\underline{\mathbf X}\in\mathbb{R}^{{I_1}\times I_2 \times \cdots\times {I_N}}$, then the $n$-unfolding of the tensor $\underline{\mathbf X}$ is denoted by ${\mathbf X}_{(n)}\in\mathbb{R}^{{I_n}\times {I_1\cdots I_{n-1}I_{n+1}\cdots I_N}}$, and is constructed by stacking all its $n$-mode fibers.

The {\em mode-$n$ product} is a generalization of matrix-matrix multiplication, defined for $\underline{\mathbf X}\in\mathbb{R}^{I_1\times I_2\times \cdots\times I_N}$ and ${\mathbf B}\in\mathbb{R}^{J\times I_n}$, as
\[
{\left( {{\underline{\mathbf X}}{ \times _n}{\mathbf B}} \right)_{{i_1},\ldots {,i_{n - 1}},j,{i_{n + 1}}, \ldots, {i_N}}} = \sum\limits_{{i_n=1}}^{{I_N}} {{x_{{i_1},{i_2}, \ldots ,{i_N}}}{b_{j,{i_n}}}},
\]
for $j=1,2,\ldots,J$. Here, we have

\[
\underline{\mathbf X}\times_{n}{\mathbf B}\in\mathbb{R}^{I_1\times \cdots \times I_{n-1}\times J\times I_{n+1}\times \cdots \times I_N}.
\]
Let $\underline{\mathbf X}\in{\mathbb{R}^{{I_1} \times I_2\times  \cdots  \times {I_N}}}$ be a given tensor, then the Tucker decomposition model is defined as follows: \cite{tucker1963implications,tucker1964extension,tucker1966some}
\begin{equation}\label{HOSVD}
\underline{\mathbf X} \cong \underline{\mathbf S}{ \times _1}{{\mathbf A}^{\left( 1 \right)}}{ \times _2}{{\mathbf A}^{\left( 2 \right)}} \cdots { \times _N}{{\mathbf A}^{\left( N \right)}},
\end{equation}
where $\underline{\mathbf S} \in {\mathbb{R}^{{R_1} \times  R_2\times \cdots  \times {R_N}}}$ is the core tensor and ${\mathbf A}_n \in {\mathbb{R}^{{I_n} \times {R_n}}},$ $\,{R_n} \le {I_n},\,\,n=1,2,\ldots,N$ are the factor matrices. 
The $N$-tuple $(R_1,R_2,\ldots,R_N)$ is called multilinear or Tucker rank where $R_n$ is the rank of the mode-$n$ unfolding matrix ${\mathbf X}_{(n)},\,n=1,2,\ldots,N$. 

In the next section, we introduce the Cross Matrix Approximation (CMA) and Cross Tensor Approximation (CTA) algorithms. 

\section{Cross Matrix Approximation (CMA) and Cross Tensor Approximation (CTA) methods}\label{Sec:TCUR}
Computing a low-rank approximation of a given matrix based on a part of its individual columns and rows was first developed in \cite{goreinov1997theory} and is known as {\em skeleton or cross approximation}. Since the matrix intersecting the columns and rows is used within the approximation procedure, this framework is called cross approximation. The column or row selection can be performed in either a randomized or deterministic manner. The Maxvol-based algorithm \cite{goreinov2010find,goreinov2001maximal}, Cross2D algorithm \cite{savostyanov2006polilinear,tyrtyshnikov2000incomplete} and Discrete Empirical Interpolatory Method (DEIM) \cite{chaturantabut2009discrete,sorensen2016deim} are known algorithms for such selections. If the columns or rows are selected randomly, this framework is often known as randomized CUR approximation. Compared with conventional algorithms, e.g., SVD, this approach requires less memory usage and floating-point operations. It can also preserve the structure of the original data matrix, such as nonnegativity, sparsity, and smoothness.

\begin{figure}
\begin{center}
\includegraphics[width=8.5 cm,height=3.8 cm]{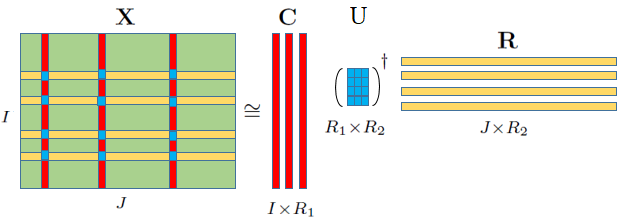}\\
 \caption{\small{Illustration of cross  decomposition for low-rank matrix approximation ${\mathbf X} \cong \mathbf{CUR}$ where ${\mathbf U} = {{\mathbf W}^+ }$ \cite{cichocki2016tensor}.}}\label{CUR}
\end{center}
\end{figure}
The problem is formally formulated as follows. Let ${\mathbf X}\in\mathbb{R}^{I\times J}$ be a given matrix and ${\mathbf C}\in\mathbb{R}^{I\times R_1},\,{\mathbf R}\in\mathbb{R}^{J\times R_2}$ are the selected columns and rows, respectively, and the intersection matrix is ${\mathbf W}\in\mathbb{R}^{R_1\times R_2}$, see Figure \ref{CUR}. The low-rank cross approximation is computed as follows:
\begin{equation}\label{middle}
{\mathbf X} \cong {\mathbf C}{\mathbf U}{\mathbf R}={\mathbf U}{\times _1}{\mathbf C}{\times _2}{\mathbf R}, 
\end{equation}
where ${\mathbf U}\in\mathbb{R}^{R_1\times R_2}$
should be computed to yield the smallest error. The best middle matrix ${\mathbf U}$ in the least-squares sense is ${\mathbf U} = {{\mathbf C}^+ }{\mathbf X}{{\mathbf R}^+ }$ because
\[{{\bf C}^+ }{\bf X}{{\bf R}^+ } = \mathop {\arg \min }\limits_{{\bf U} \in {\mathbb{R}^{ R1\times R2}}} \,\,{\left\| {{\bf X} - {\bf U}{ \times _1}{\bf C}{ \times _2}{\bf R}} \right\|_F}.\]
The approximation computed in the above formulation is exact if ${\rm{rank}}\left( {\mathbf{X}} \right) \le \min \left\{ {{R_1},{R_2}} \right\}$ \cite{goreinov1997theory}.

A faster way for computing the matrix ${\bf U}$ in \eqref{middle} is using the Moore-Penrose of the intersection matrix ${\bf W}$ and computing the approximation ${\bf X} \cong {\bf C}{{\bf W}^+ }{\bf R}$. It is known that the quality of the intersection matrix quite depends on the module of the determinant of the intersection matrix, which is called {\em matrix volume} and an intersection with as much volume as possible should be selected \cite{goreinov2010find}.

A special case of the CMA in which only columns are sampled is called {matrix column selection}, {\it interpolative matrix decompositions}, and in some contexts, it is also referred to as {\it CY decomposition}. This problem is known as column (feature) selection in the field of machine learning and data analysis \cite{frieze2004fast,boutsidis2009improved,boutsidis2014near,deshpande2006adaptive,deshpande2006matrix,deshpande2010efficient,guruswami2012optimal}. Let ${\mathbf X}\in\mathbb{R}^{I\times J}$  be a given data matrix and ${\mathbf C}\in\mathbb{R}^{I\times R}$ be the number of selected columns. Then the matrix column selection is formulated as follows
\begin{equation}\label{ITM}
{\mathbf X} \cong {{\mathbf C}{\mathbf Y}},
\end{equation}
where ${\mathbf C}\in\mathbb{R}^{I\times R}$ is a matrix containing the selected columns and ${\mathbf Y}\in\mathbb{R}^{R\times J}$ should be computed in such a way that the approximate error should be as small as possible. The best solution to the problem \eqref{ITM} in the least-squares (LS) sense is
${\mathbf Y} = {{\mathbf C}^+ }{\mathbf X},$
and if ${\rm rank}({\mathbf X})=R$, then the approximation is exact, i.e., ${\mathbf X}={\mathbf C}{\mathbf C}^+ {\mathbf X}$.

Cross Tensor Approximation (CTA) is a generalization of cross/skeleton matrix or CUR
matrix approximation and is an efficient approach for fast low-rank tensor approximation. The {\em Skeleton} or {\em Cross approximation} or equivalently matrix CUR approximation (CMA) computes a low-rank approximation based on a part of individual columns and rows. It has found applications in deep learning \cite{mai2020vgg}, signal processing \cite{hendryx2018finding,cai2020rapid}, scientific computing \cite{drineas2006fastI,drineas2006fastII,drineas2006fastIII} and machine learning \cite{mahoney2009cur,li2018joint,aldroubi2019cur}. 
In the next subsequent sections, we explain how the CMA can be generalized to tensors. 
In general, there are three main categories under which the CMA techniques are generalized to tensor as follows \cite{ahmadi2021cross}:
\begin{itemize}
    \item Fiber selection \cite{drinea2001randomized,oseledets2008tucker,caiafa2010generalizing},
    \item Slice-tube selection \cite{mahoney2008tensor},
    \item Slice selection \cite{tarzanagh2018fast}.
\end{itemize}
The first category is concerned with the Tucker decomposition in which the factor matrices ${\bf A}^{(n)},\,n=1,2,\ldots,N$ are computed by sampling the $n$-mode fibers of a given data tensor $\underline{\mathbf X}$ and the core tensor is computed in different ways, see Appendices B-D for the details. The essential differences between the algorithms proposed in \cite{drinea2001randomized,oseledets2008tucker,caiafa2010generalizing} are 1) the number of fibers required for approximation and 2) the way that the core tensor is computed. For example, in \cite{drinea2001randomized}, for a given Tucker rank $(R_1,R_2,\ldots,R_N)$, $R_n$ numbers of mode-$n$ fibers are selected to approximate the factor matrix ${\bf A}^{(n)}\,n=1,2,\ldots,N$, and the core tensor is computed as follows 
\begin{eqnarray}
\underline{\mathbf S} = \underline{\mathbf X}{ \times _1}{\mathbf A}_1^{+} { \times _2}{\mathbf A}_2^{+}\cdots { \times _N} {\mathbf A}_N^{+}\in\mathbb{R}^{R_1\times R_2\times \cdots\times R_N},
\end{eqnarray}
while in the Fast Sampling Tucker Decomposition (FSTD) \cite{caiafa2010generalizing}, the core tensor is constructed by sampling a part of the components of the tensor $\underline{\mathbf X}$. 

In the second category, some slices are selected along with the fibers. Here, the following model is considered 
\begin{equation}\label{fs}
\underline{\mathbf X} \cong \underline{\mathbf C}{ \times _3}\left( {\mathbf U{\mathbf R}} \right)^T,
\end{equation} 
where the tensor $\underline{\mathbf C}\in\mathbb{R}^{I_1\times I_2\times L_1}$ and the matrix ${\mathbf R}\in\mathbb{R}^{L_2\times I_3}$ contain the sampled frontal slices and tubes, respectively. The matrix ${\bf U}\in\mathbb{R}^{L_1\times L_2}$ is computed so that the approximation error obtained from model \eqref{fs} be as small as possible. 

The last category is associated with the tubal SVD. More precisely, the tubal cross approximation based on t-product is formulated as follows
\begin{equation}\label{CURTf}
\underline{\mathbf X} \cong \underline{\mathbf C} * \underline{\mathbf U} * \underline{\mathbf R},
\end{equation}
where $*$ stands for the t-product \cite{kilmer2011factorization} (see Appendix A), $\underline{\mathbf C}\in\mathbb{R}^{I_1\times L_1\times I_3}$ and $\underline{\mathbf R}\in\mathbb{R}^{L_2\times I_2\times I_3}$ are some sampled lateral and horizontal slices of the original tensor $\underline{\mathbf X}$ respectively and the middle tensor $\underline{\mathbf U}\in\mathbb{R}^{L_1\times L_2\times I_3}$ is computed in such a way that the approximation \eqref{CURTf} should be as small as possible. For the brevity of the presentation, we have provided the details of different tensor CUR approaches in Appendices B-D.

\section{Tensor Completion}\label{Sec:TCProb}
Tensor completion is the problem of completing (reconstructing) an incompleted data tensor from only a part of observed data components. It includes the matrix completion problem as a special case because a matrix is a second-order tensor. The completion algorithms are mainly categorized as 1) rank minimization and 2) tensor factorization techniques. The former is formally stated as follows:
\begin{equation}\label{MinRankCompl1}
\begin{array}{cc}
\displaystyle \min_{\underline{\bf X}} & {\rm Rank} \left(\underline{\bf X}\right)\\
\textrm{s.t.} & {{\bf P}_{\underline{\bf\Omega}} }({\underline{\bf X}})={{\bf P}_{\underline{\bf\Omega}} }({\underline{\bf M}}),\\
\end{array}
\end{equation}
where $\underline{\bf M}$ is the incompleted data tensor with only observed components and their corresponding indices $\underline{\bf \Omega}$ and $\underline{\bf X}$ is the unknown data tensor that needs to be determined. The rank mentioned in formulations \eqref{MinRankCompl1} can be tensor rank, TT (TC) rank, Tucker rank, etc, and the corresponding minimization problem is considered. For example in the case of Tucker rank which is the $N$-tuple $(R_1,R_2,\ldots,R_N)$,  we have the following problem
\begin{equation}\label{MinRankCom}
\begin{array}{cc}
\displaystyle \min_{\underline{\bf X}} & \sum_{n=1}^{N}R_n\\
\textrm{s.t.} & {{\bf P}_{\underline{\bf\Omega}} }({\underline{\bf X}})={{\bf P}_{\underline{\bf\Omega}} }({\underline{\bf M}}),\\
\end{array}
\end{equation}
while for the TT model, the TT-rank \cite{oseledets2011tensor} which is an $(N-1)$-tuple should be minimized. The rank minimization problem, in general, is an NP-hard problem because this problem for the matrices as a subset of tensors is NP hard \cite{fazel2002matrix}. So convex and tractable surrogates of the rank usually are replaced \cite{song2019tensor,long2019low}. The nuclear norm is the convex envelop of the rank and was extensively used in the literature but it is prohibitive especially for large-scale data due to the utilization of SVD \cite{fazel2002matrix}. The tensor decomposition formulation is a more efficient alternative to the nuclear norm minimization approach. The tensor decomposition variant of formulation \eqref{MinRankCompl1} is presented as follows

\begin{equation}\label{MinRankCompl2}
\begin{array}{cc}
\displaystyle \min_{\underline{\bf X}} & {\|{{\bf P}_{\underline{\bf\Omega}} }({\underline{\bf X}})-{{\bf P}_{\underline{\bf\Omega}} }({\underline{\bf M}})\|^2_F},\\
\textrm{s.t.} & {\rm Rank}(\underline{\bf X})=R,\\
\end{array}
\end{equation}
where the unknown tensor $\underline{\bf X}$ has low tensor rank representation. Here again, different kinds of tensor ranks and associated tensor decompositions can be considered. 
Using an auxiliary variable ${\underline{\bf C}}$, the optimization problem \eqref{MinRankCompl2} can be solved more conveniently by the following reformulation
\begin{equation}\label{MinRankCompl3}
\begin{array}{cc}
\displaystyle \min_{\underline{\bf X},\underline{\bf C}} & {\|{\underline{\bf X}}-{\underline{\bf C}}\|^2_F},\\
\textrm{s.t.} & {\rm Rank}(\underline{\bf X})=R,\\
& {{\bf P}_{\underline{\bf\Omega}} }({\underline{\bf C}})={{\bf P}_{\underline{\bf\Omega}} }({\underline{\bf M}})\\
\end{array}
\end{equation}
and we can alternatively solve optimization problem \eqref{MinRankCompl2} over variables $\underline{\bf X}$ and $\underline{\bf C}$. Thus, the solution to the minimization problem \eqref{MinRankCompl2} can be approximated by the following iterative procedure
\begin{equation}\label{Step1}
\underline{\mathbf X}^{(n)}\leftarrow \mathcal{L}(\underline{\mathbf C}^{(n)}),
\end{equation}
\begin{equation}\label{Step2}
\underline{\mathbf C}^{(n+1)}\leftarrow\underline{\mathbf \Omega}\oast\underline{\mathbf M}+(\underline{\mathbf 1}-\underline{\mathbf \Omega})\oast\underline{\mathbf X}^{(n)},
\end{equation}
where $\mathcal{L}$ is an operator to compute a low-rank tensor approximation of the data tensor $\underline{\mathbf X}^{(n)}$ and $\underline{\mathbf 1}$ is a tensor whose all components are equal to one. Note that Equation \eqref{Step1} solves the minimization problem \eqref{MinRankCompl3} over  
$\underline{\bf X}$ for fixed variable $\underline{\bf C}$. Also Equation \eqref{Step2} solves the minimization problem \eqref{MinRankCompl3} over  
$\underline{\bf C}$ for fixed variable $\underline{\bf X}$. The algorithm consists of two main steps, {\it low rank tensor approximation} \eqref{Step1} and {\it Masking computation} \eqref{Step2}. It starts from the initial incomplete data tensor $\underline{\mathbf X}^{(0)}$ with the corresponding index set $\underline{\mathbf \Omega}$ and  sequentially improves the approximate solution till some stopping criterion is satisfied or the maximum number of iterations is reached. Note that the term $\underline{\mathbf \Omega}\oast\underline{\mathbf M}$ is not required to be computed at each iteration as it is equal to the initial data tensor $\underline{\mathbf X}^{(0)}$. In general, any type of tensor decomposition can be utilized for low-rank approximation in \eqref{Step1}. For example, in \cite{fang2018sequentially}, the Tucker decomposition is exploited at each iteration while in \cite{acar2011scalable} and \cite{wang2017efficient} the CPD and the TR/TC decomposition are used. In the existing papers, usually the deterministic algorithms are utilized for computation of the tensor decomposition which may be expensive especially when a large number of iterations is required for convergence or the data tensor is quite large. In this paper, we propose to use the tensor CUR algorithms instead of the deterministic counterparts in which the actual fibers or slices are used for the low-rank approximation. A main feature of the paper is proposing a general framework to utilize the tensor CUR algorithm for the completion task. More precisely, all types of the tensor CUR algorithms including fiber sampling, slice-tube sampling or slice sampling can be exploited. It is experimentally shown that this can accelerate the basic completion algorithms significantly.

\section{Related works}\label{Sec:Related}
{The tensor completion problem \cite{song2019tensor} as a generalization of the matrix completion problem \cite{candes2009exact} has been extensively studied in the past decade. Similar to the matrix completion, the low-rank property is the main assumption used in formulation of the tensor completion problem. Liu et.al were the first who used the tensor decompsoitions (Tucker and CPD models) for solving the tensor completion problem  \cite{liu2012tensor2,liu2012tensor}. The Bayesian CPD was proposed in \cite{zhao2015bayesian} to automatically select an appropriate tensor rank for reconstructing incomplete data tensors. Besides, to better reconstruct the images the CPD with smoothness constraint was proposed in \cite{yokota2016smooth} where the tensor rank needs to be known in advance. The tensor completion using the tubal rank was studied in \cite{zhang2016exact} and the results showed that it can provide promising results. The TT and TC models were later utilized in \cite{bengua2017efficient} and \cite{wang2017efficient,yuan2018higher}, respectively to recover images/videos with missing pixels and have been shown to efficient for large-scale data tensors. Recently the idea of tensor Hankelization was proposed in \cite{yokota2018missing} to tackle the reconstruction of data images/videos with structural missing pixels, e.g. removing several frames of a video. The Hankelization approaches are known for their slow convergence and a fast version of them were developed in \cite{yamamoto2022fast}. The tensor completion problems are categorised into 1) rank minimization and 2) tensor factorization classes. In the first category, the nuclear norm minimization is used to estimate the tensor rank while in the second category the tensor rank needs to be known in advance. The Alternating Direction Method of Multipliers (ADMM) and Alternating Least Squares (ALS) frameworks are used to solve the first and second categories respectively \cite{boyd2011distributed}. In both categories, at each iterations, we need to compute a SVD or a low-rank approximation of some large-scale unfolding matrices, see \cite{song2019tensor,long2019low} and the references therein for details. The main challenge of completion approaches is reducing the computational complexity of these computations and in this paper we utilize the CUR framework which is a class of randomized algorithms called sampling methods \cite{halko2011finding}. It is interesting to note that recently Deep Learning (DL) algorithms have been also adopted for the completion task \cite{ulyanov2018deep,rombach2022high,qin2021image}. It is known that the DL algorithms are data hungry and need a lot of samples to train the model. On the other hand, the training procedure needs too much computational resources and time\footnote{One should also carefully tune the parameters and regularize the model to hinders the deep learning from over-fitting.}. As a result, since DL models are trained on a lot of images, they learn and capture richer information and can better encode the visual patterns. In contrast, the tensor completion models such as the proposed methods rely only on the low-rank assumption and use a single image to estimate the missing pixels. So they can provide reasonable result in less time and using less computational resources.}

\section{Proposed Approach}\label{Sec:ProApp}
In this section, we describe our proposed algorithm for the tensor completion task. A key computation in our algorithm is the CUR approximation of the underlying data tensors. More precisely, at each iteration of our algorithm, a CUR approximation of the underlying data tensor is computed, after which the mask operator is applied. Indeed, the operator $\mathcal{L}$ in \eqref{Step1} is replaced by the tensor CUR algorithms in Section \ref{Sec:TCUR}. The algorithm is initialized with a random data tensor with zero mean and unite variance and is  sequentially updated to reconstruct the incompleted data tensor. We experimentally confirmed that the same results are achieved if the algorithm starts with the initial incomplete data tensor. Our proposed algorithm is summarized in Algorithm 1. This formulation is quite general and can incorporate different tensor/matrix factorization cases. For example, for the case of matrices, we have the matrix CUR algorithm which we only select some columns and rows (${\bf C}\in\mathbb{R}^{I_1\times R}$ and ${\bf R}\in\mathbb{R}^{R\times I_2}$) and compute the middle matrix as ${\bf U}={\bf C}^+ {\bf X}{\bf R}^+\in\mathbb{R}^{R\times R}$ and a low CUR approximation ${\bf X} \cong {\bf CUR}$ is computed.

Also, for the tensor case, different tensor decompositions can be utilized, such as Tucker decomposition \cite{tucker1964extension,tucker1966some}, tubal decomposition \cite{kilmer2011factorization,kilmer2013third}, tensor CUR approximation \cite{mahoney2008tensor}. To be more precise, let us consider the Tucker decomposition case. Here,  in the first stage, the factor matrices ${\bf C}_n\in\mathbb{R}^{I_n\times R_n}$ are computed by sampling the columns of $n$-unfolding matrices (or $n$-mode fibers) after which the core tensor is computed as follows
\begin{eqnarray}\label{CorInter}
\underline{\mathbf S} = \underline{\mathbf X}{ \times _1}{\mathbf C}_1^{+} { \times _2}{\mathbf C}_2^{+}\cdots { \times _N} {\mathbf C}_N^+\in\mathbb{R}^{R_1\times R_2\times \cdots\times R_N},
\end{eqnarray}
Then, a low Tucker rank approximation is computed at each iteration as follows
\begin{eqnarray}\label{CorInter}
\underline{\mathbf X} = \underline{\mathbf S}{ \times _1}{\mathbf C}_1 { \times _2}{\mathbf C}_2\cdots { \times _N} {\mathbf C}_N.
\end{eqnarray}
Of course, other kinds of CUR Tucker approximations such as FSTD/adaptive fiber sampling \cite{caiafa2010generalizing}, Cross3D \cite{oseledets2008tucker} or TT-Cross approximation \cite{oseledets2010tt} can be applied. In all our experiments, the sampling procedure is performed randomly, and as a result, the proposed algorithm can be considered a randomized tensor completion algorithm. It is also possible to exploit heuristic ones as it is a special case of our algorithm. However, they are slower than the randomized ones. In the case of images, we basically use the Tucker-2 model and only select columns/rows as the third mode is small and it is always considered as 3. In the case of videos,  fibers in all modes are selected because the third mode which is the time (number of frames) is large.

\begin{figure}
\begin{center}
\includegraphics[width=1\columnwidth]{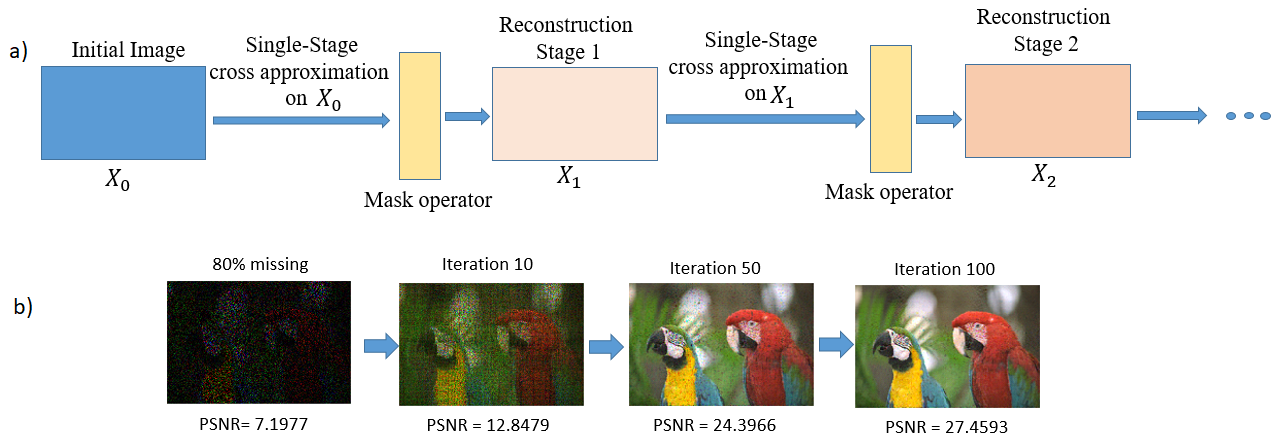}\\
\caption{\small{(a) The procedure of the proposed algorithm as a  multi-stage CUR approximation followed by mask operator. (b) The proposed approach for reconstructing an image with $80\%$ missing pixels for different number of iterations.}}\label{graph}
\end{center}
\end{figure}

The CY approximation can be incorporated in Algorithm 1. We have not considered it in our simulations as it is a special case of the CUR approximation, and there will not be a significant difference in the results. Moreover, our simulations show that for the Tucker CUR and the tubal CUR decompositions, the middle matrix should be computed very carefully and accurately otherwise, the approximation scheme will be unstable, and the results will be very poor. In all our experiments for the Tucker CUR and tubal CUR, we have utilized formulas \eqref{CorInter} and \eqref{TCUR}, respectively. 

{Our proposed algorithm starts with an incomplete image or a random tensor and computes a tensor CUR approximation of this initial data tensor followed by applying the mask operator to keep fixed its known pixels. Let the approximation  computed in the first stage be denoted by $\underline{\bf X}^{(1)}$ which is a single-stage CUR approximation. The single-stage CUR approximation does not work well in practice and we need multi-stage CUR approximation 
by concatenating several single-stage CUR approximation. More precisely, we continue the above-mentioned procedure by computing a CUR approximation of the tensor $\underline{\bf X}^{(1)}$ and applying again the mask operator. This procedure is continued till a maximum number of iterations is reached or a stopping criterion is satisfied. This method also can be interpreted as a combination of a sequence of CUR approximation (linear operator) and a mask function (nonlinear operator). For a graphical illustration on the proposed approach see Figure \ref{graph} (a). In \ref{graph} (b), we have also shown that the multistage CUR approximation provides better results than the single-stage CUR approximation. We can see that from an incomplete image with PSNR $7.1977$ after 10 iterations/stages, we reconstruct an image with $12.8479$ while for 50 and 100 iterations, we achieve PSNR $24.3966$ and $27.4593$, respectively. We experimentally found that  smoothing the selected fibers in the tensor CUR approximation before the tensor reconstruction, can totally improve the performance of the algorithm. In particular, this was more visible for the difficult scenarios, i.e. high missing ratio or structural missing patterns. 
In the next section, we will explain some smoothing approaches which we have used in our experiments.}

\subsection{Smoothing Techniques}
 The smoothing techniques, also known as {\it low-pass filtering} or {\it curve fitting} tools, are useful pre-processing approaches used in the machine learning community.
In many applications, the elements of the underlying data tensors change smoothly (continuously), e.g., images, videos, time series. It is natural to consider the smoothness constraint when we are working on the optimization problems associated with the mentioned data tensors. {The smoothing techniques have been also used with deep learning models for haze visibility
enhancement in some internet of thing applications \cite{liu2022deep}.} There are various techniques  and approaches to make the elements of a data tensor smooth such as  moving average, locally estimated scatterplot smoothing (LOESS) and its weighed variant (LWOESS), the robust LOESS (RLOESS) and robust LWOESS (RLWOESS) etc., see for \cite{savitzky1964smoothing,cleveland1981lowess,cleveland1979robust,garimella2017simple}, for more details. The moving average is the simplest type of smoothing procedure used widely in the signal community as finite impulse response (low-pass) filter. It smooths out the elements of given data, by replacing the elements of the data with a sequence of averages of different subsets of the given data points. To be more precise, assume that $({ x}_1,{y}_1),\,({ x}_2,{y}_2),\ldots,({ x}_N,{ y}_N)$, are given data points, where $x_i$ and $y_i$ stand for independent and dependent variables, respectively. For the above data points, the so-called cumulative moving average gives a new  sequence ${ z}_1,\,{z}_2,\,\ldots,{z}_N$ for the dependent variables $y_i$ as follows
\begin{eqnarray}
\nonumber
{ z}_1&=&{ x}_1,\\
\nonumber
{ z}_1&=&\frac{{ x}_1+{ x}_2}{2},\\
\nonumber
\vdots
\\
{ z}_N&=&\frac{{ x}_1+{ x}_2+\ldots+{ x}_N}{N},
\end{eqnarray}
while the simple moving average is defined as 
\[
z_n=\frac{1}{n}\sum_{i=N-n+1}^{N} y_i,
\]
for $n=1,2,\ldots,N$. Other types of moving average include central, exponential and weighted moving average. The data points can also be made smooth by the regression strategy. In the regression approach, the goal is fitting the given data points by a function $f$, i.e., $y_i=f(x_i,\beta)+e_i,\,i=1,2,\ldots,N$, where $e_i$ are the prediction error and $\beta$ is the coefficients vector used in the data fitting. For example, the linear function $y_i=\beta_0+\beta_1x_i+\epsilon_i,\,i=1,2,\ldots,N$ or the parabolic function 
$y_i=\beta_0+\beta_1x_i+\beta_2x^2_i+\epsilon_i,\,i=1,2,\ldots,N$ can be used in which the coefficient parameter $\beta$ are $\beta=[\beta_0,\beta_1]$ and $[\beta_0,\beta_1,\beta_2]$, respectively. By minimizing the sum of the errors $e_n,\,n=1,2,\ldots,N$ and solving the following minimization problem 
\[
\min_{\beta} \sum_{n=1}^{N}e_i=\min_{\beta}\sum_{n=1}^{N}(y_n-f(x_n,\beta)),
\]
the underlying coefficient vector $\beta$ is computed. The above procedure fits the whole data points by a single line or a single polynomial. However it is more desirable to fit each part of data by a line or a polynomial separately to get more accurate approximations. This is called {\it local regression} and is a special case of the regression problem. The LOESS and the WLOESS are two examples of the local regression techniques. We will exploit a variety of smoothing techniques in our experiments and compare their performance for reconstructing data tensors.

As discussed in Section \ref{Sec:TCUR}, the CMA methods can be generalized to tensors by sampling fibers, slices-tubes and slices. We proposed to first smooth out the selected fibers/slices and then compute the corresponding low CTA. We experimentally confirmed that this idea always provides better results than the algorithm without imposing smoothness. More importantly, in the situation that we have a high missing ratio of elements or structural missing components, e.g., sequential columns and rows, the basic Algorithm 1 may not work properly while its smooth variant works perfectly. 

\begin{figure}
\begin{center}
\includegraphics[width=0.9\columnwidth]{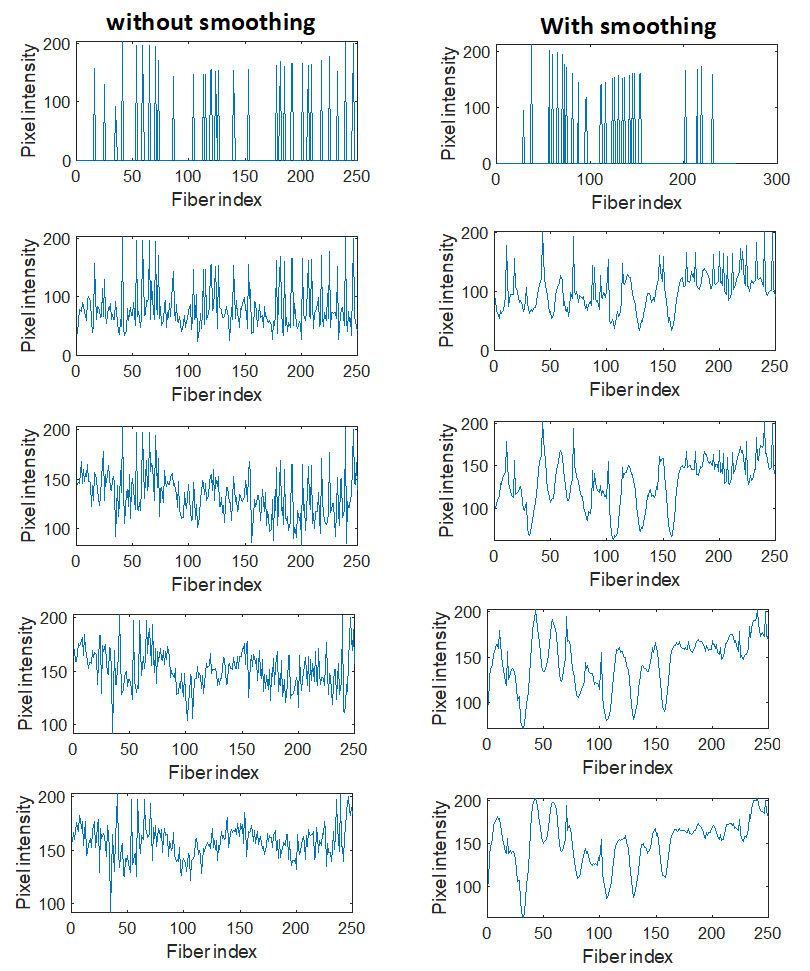}\\
\caption{\small{Comparison of the structure of smooth fibers and nonsmooth in CUR algorithms. The results are for the column fiber $\underline{\mathbf X}(:,50,1)$ using the Tucker CUR and the iterations $1,\,10,\,20,\,30,\,40$, (from the top to the bottom). Vertical axis represents the number of pixels, while the vertical axis represents the corresponding pixel intensity.}}\label{Smoothsignal}
\end{center}
\end{figure}

In our extensive simulations, we found that for difficult incomplete data tensors, e.g., removing sequential columns and rows or a high missing ratio ($95\%$), the tensor CUR algorithms may have a problem in reconstructing the incomplete data tensors. This is because the selected fibers of the incomplete data tensors have many zero components, and their components do not change continuously (smoothly) as expected in the case of images/videos. To overcome this difficulty, we proposed to first smooth out the selected fibers and then compute the CTA, see Figure \ref{Smoothsignal} for the difference between the structure of a smooth signal and nonsmooth ones. We have tried different types of smoothing techniques in our computations, including moving average, LOESS, LWOESS, LOESS, RLOESS, RLWOESS, and Savitzky-Golay. In all our experiment we used the MATLAB command "smooth" can be utilized for the mentioned smoothing techniques. In fact, in the simulation part, we show that Algorithm 1 with smoothing the fibers provides better results than those without smoothing. As we will show in the simulations, this idea significantly improves the Tucker CUR algorithm while it requires almost the same running time. 

\begin{figure}
    \centering
    \includegraphics{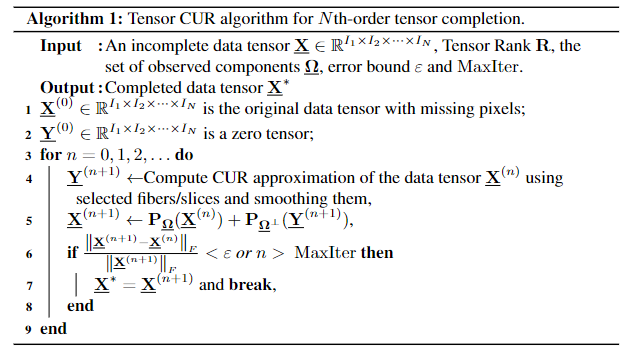}
\end{figure}

\section{Computational Complexity}\label{Sec:CCom} In this section, we make a brief comparison between the complexity of the completion algorithms. Let $\in\mathbb{R}^{I_1\times I_2\times \cdots\times I_N}$. For the simplicity, assume $I_1=I_2=\ldots=I_N=I$ and for the Tucker rank $(R_1,R_2,\ldots,R_N)$, let $R_1=R_2=\ldots=R_N=R$. The computational complexities of the TR-ALS, and TR-WOPT per iteration are $\mathcal{O}(PNR^4I^N+NR^6)$, and $\mathcal{O}(2NI^NR^2)$, respectively, where $P$ is the total number of observations. The computational complexity of Tucker CUR, Tubal CUR, FSTD and Slice-tube CUR are $\mathcal{O}(I^NR),\,\mathcal{O}(I^{N-1}\log(I)R),\,\mathcal{O}(IR^N)$ and $\mathcal{O}(I^NR)$, respectively, all of them are lower than the other completion algorithms. Our experimental results also confirm this. In our experiments, the TR-WOPT required more iterations for convergence and this is why in Figure \ref{Runtimecomhouse}, and \ref{RuntimeForeman}, the TR-WOPT has high running time.

\section{Simulations}\label{Sec:Sim}
In this section, we examine the proposed CUR algorithms on image/video completion. All numerical simulations were performed on a laptop computer with 2.60 GHz Intel(R) Core(TM) i7-5600U processor and 8GB memory. In all our experiments, the sampling of the fibers/slices are performed without replacement. We use Pick to Signal of Noise Ration (PSNR) and Structural SIMilarity (SSIM) to compare the performance of different completion algorithms. It is worth mentioning that it has been recently shown that the Hankelization completion algorithms can achieve high performance for recovering incomplete images with structured missing components, however we have not compared our algorithms with them as they are much slower and have higher computational complexity because the size and order of the data tensor increased in a pre-processing step (Hankelization step). 

Besides, since the running time of the tensor CUR algorithms are less than the other completion algorithms, we only compare the Tucker CUR with them as the baseline and the running time of tensor CUR algorithms are compared with themselves. We mainly use tubal CUR, Tucker CUR, FSTD, and slice-tube CUR algorithms and compare them with TR-WOPT \cite{yuan2018higher}, TR-ALS \cite{wang2017efficient}, SPC \cite{yokota2016smooth}, and Fast MDT \cite{yamamoto2022fast} which are among the state-of-the-art algorithms for the tensor completion task. 
The codes are accessible at 
\url{https://github.com/SalmanAhmadi-Asl/Cross_Tensor_Completion}.

\begin{figure}
\begin{center}
\includegraphics[width=0.9\columnwidth]{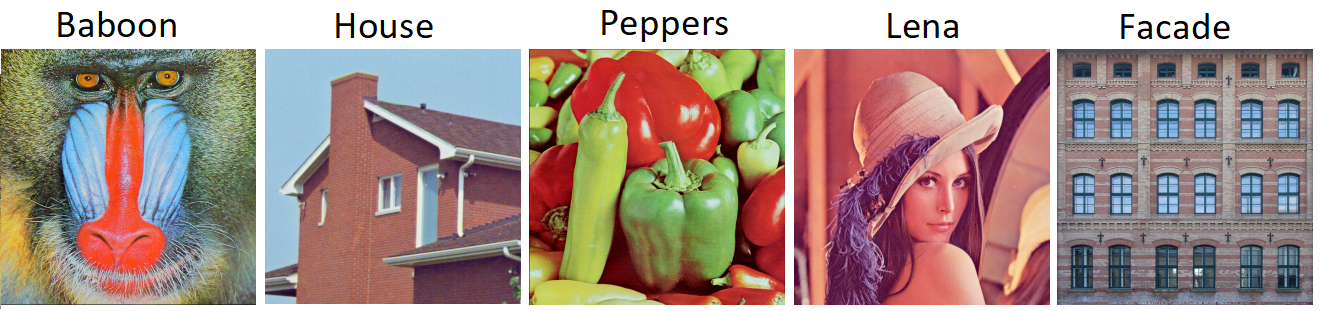}
\caption{\small{Benchmark images used in our experiments.}}\label{Pic:Samples}
\end{center}
\end{figure}

{\bf Image completion} We first consider the image completion case. The benchmark images used in our simulations are "Baboon", "House", "Peppers", "Lena", "Facade" shown in Figure \ref{Pic:Samples}. All of them have size $256\times 256 \times 3$. We considered different kinds of missing components (randomly and structured) shown in the second column of Figure  \ref{Compare1} and the reconstructed images using different completion algorithms along with their corresponding PSNR and SSIM and running times are reported in Figure \ref{Compare1}. In Table \ref{Tab1}, detailed information regarding each of the TR-WOPT, TR-ALS, and SPC are mentioned. {For the Fast MDT algorithm we considered the delay embedding parameter $(32,32,1),$ tolerance $10^{(-4)}$, and maximum number of iterations 10000}. Since each image has only three slices and all of them are necessary to be selected, it is not possible to apply Algorithm 6 directly. To this end, we reshaped the images to a 3rd order tensor of size $64\times 64\times 48$ and then selected some slices and tubes. In our algorithm for the Tucker CUR case, we considered Tucker rank $(70,70,3)$ in which we select 70 columns, 70 rows, and 3 channels. For the tubal CUR we considered 40 lateral and 40 horizontal slices, and for the slice-tube CUR, we considered 35 frontal slices and 2500 tubes. 

\begin{figure}
\begin{center}
\includegraphics[width=0.9\columnwidth]{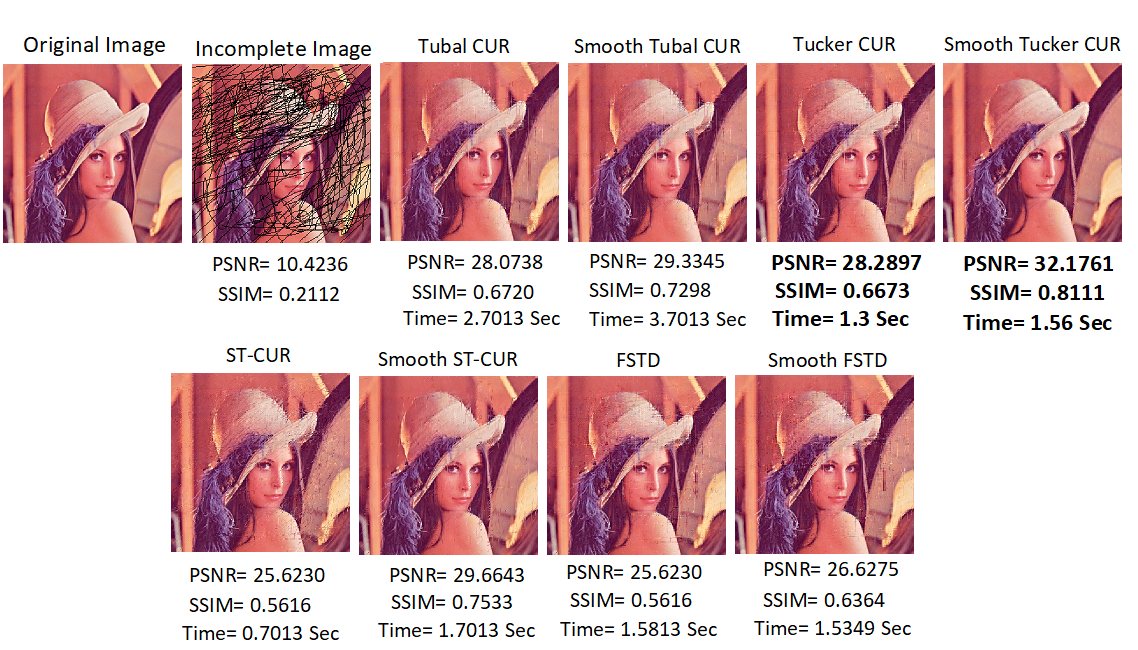}
\caption{\small{The reconstructed images using different tensor CUR completion algorithms.}}\label{CURCompletion}
\end{center}
\end{figure}

\begin{figure}
\begin{center}
\includegraphics[width=0.9\columnwidth]{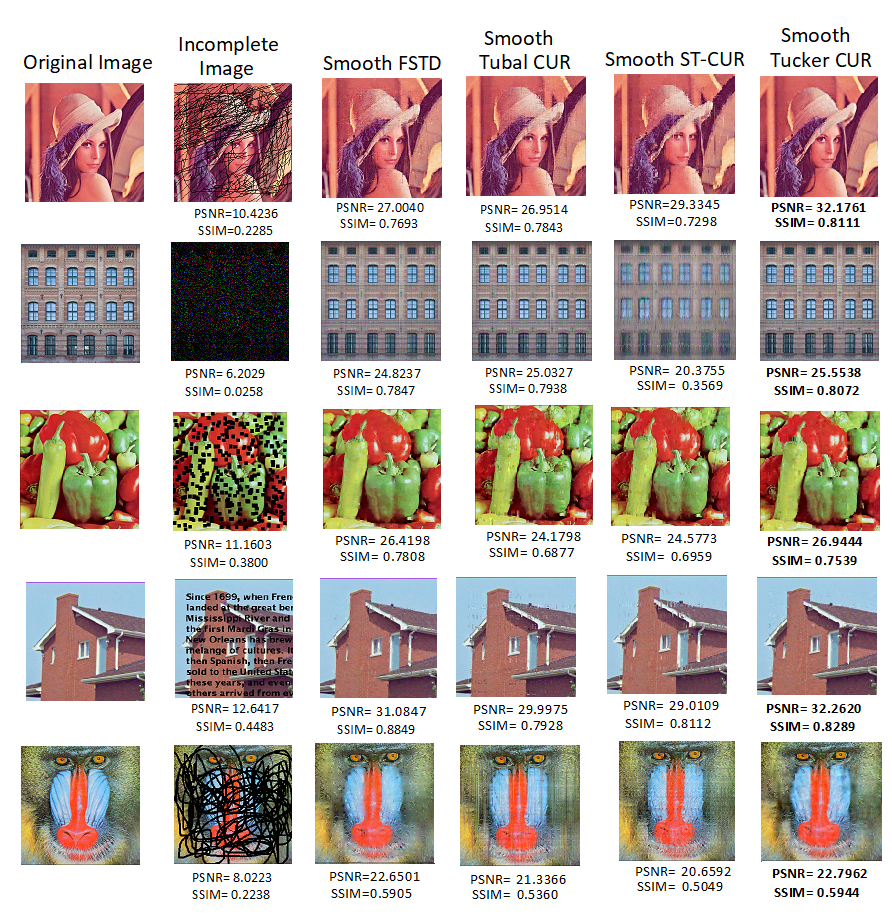}
\caption{\small{The reconstructed images using different tensor CUR completion algorithms.}}\label{Compare2}
\end{center}
\end{figure}

\begin{figure}
\begin{center}
\includegraphics[width=1\columnwidth]{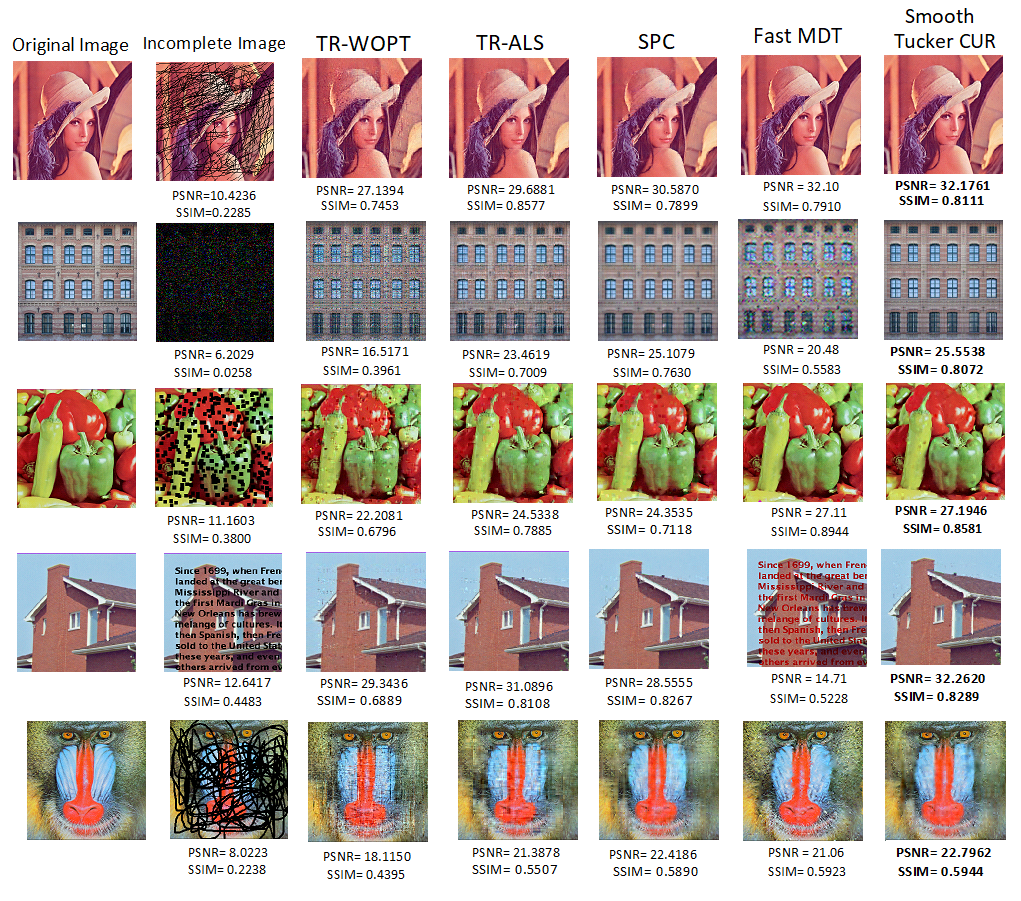}
\caption{\small{The reconstructed images using different tensor completion algorithms.}}\label{Compare1}
\end{center}
\end{figure}

In the first set of experiments we show the effect of smoothness for achieving better accuracy. To this end, we consider the "Lena" image and its degraded form as an image with incomplete pixels. We apply different CUR algorithms and their smooth variants on the degraded Lena image. The reconstructed images using different CUR completion algorithms are displayed in Figure \ref{CURCompletion}. The corresponding PSNR/SSIM and execution time required to perform the completion is also reported. The results show that the smooth CUR algorithms always work better than ordinary CUR completion algorithms. Besides, we experimentally found that the smooth Tucker CUR algorithm was the most promising algorithm for reconstructing the Lena image. We conducted new experiments with more incomplete images and again compared the performance of different CUR completion algorithms. The results are reported in Figure \ref{Compare2}. As can be seen we have used images with different missing patterns. All of them have structured missing expect the face for which we have removed $90\%$ pixels. The results again show  the superiority of the Tucker CUR algorithm over other CUR completion algorithms.  

\begin{figure}
     \centering
         \includegraphics[width=0.9\columnwidth]{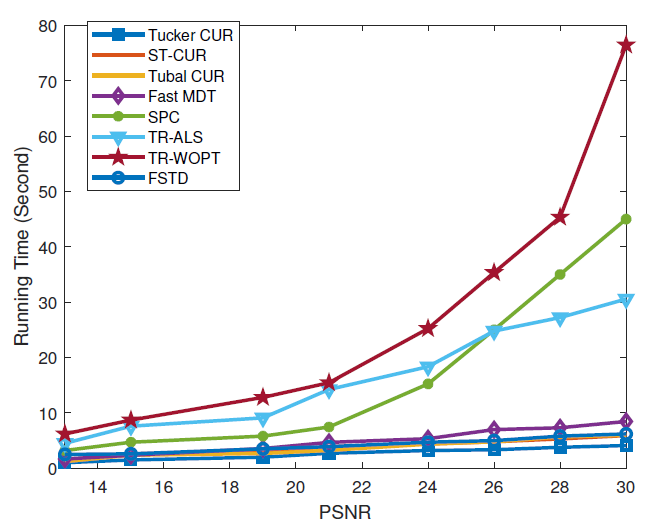}
     \caption{Comparing the performance and running times of different completion algorithms for the picture house.}\label{Runtimecomhouse}
\end{figure}

The smooth CUR Tucker completion algorithm which was the most efficient algorithm in our previous experiments is compared with the TR-WOPT, TR-ALS, and SPC algorithms. The results are reported in Figure \ref{Compare1}. From Figure \ref{Compare1}, it is seen the smooth Tucker CUR requires less time to complete the images with missing components while achieving roughly the same or even better performance compared to other completion algorithms. {The running time of different completion algorithms are summarized in Table \ref{tab:my_time}. We also report the running time against the PSNR of the house image for different completion algorithms.} From Figure \ref{Runtimecomhouse} and Table \ref{tab:my_time}, it is visible that the proposed CUR algorithms are more scalable for achieving higher performance that other completion algorithms. 
\begin{table*}
\begin{center}
\caption{The selected parameters for different completion algorithms used in our experiments for images/videos completion.}\label{Tab1}
\begin{tabular}{ |p{3.1cm}|p{4cm}|p{3cm}|}
 \hline
 \multicolumn{3}{|c|}{Images} \\
 \hline
  $\quad\quad$TR-WOPT & \hspace{1.3cm} TR-ALS & \hspace{1.2cm}SPC\\
 \hline
    Size=$(256,256,3)$   \hspace{.5cm} TR-Rank=$(5,5,5)$ \hspace{.5cm} Max $\#$ Iteration= 100&  \hspace*{.3cm} Size=$(4,4,16,4,4,16,3)$  TR-Rank=$(10,10,\ldots,10)$, \hspace*{.3cm} Max $\#$ Iteration= 10 &  Size=$(256,256,3)$ Type=Quadratic\,\,\, Max $\#$ Iteration=150 \\
 \hline
  \multicolumn{3}{|c|}{Videos} \\
  \hline
  $\quad\quad$TR-WOPT & \hspace{1.3cm} TR-ALS & \hspace{1.2cm}SPC\\
 \hline
  Size=$(176,144,30)$ \hspace{.5cm} TR-Rank=$(5,5,5)$ \hspace{.5cm} Max $\#$ Iteration= 100&  \hspace*{.3cm} Size=$(16,4,4,6,33,15)$  TR-Rank=$(10,10,\ldots,10)$, \hspace*{.3cm} Max $\#$ Iteration= 10 &  Size=$(176,144,30)$
 Type=Quadratic\quad\quad Max $\#$ Iteration=150 \\
 \hline
\end{tabular}
\end{center}
\end{table*}

\begin{figure}
\begin{center}
\includegraphics[width=0.9\columnwidth]{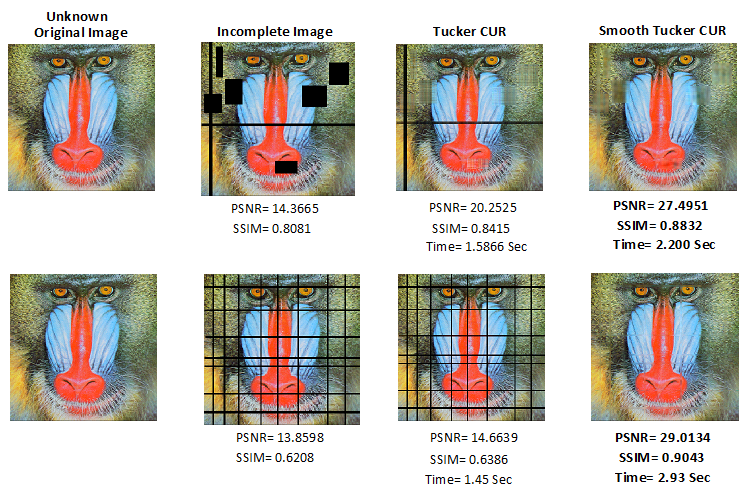}
\caption{\small{The reconstructed images using the Tucker CUR and the smooth Tucker CUR algorithm. The Tucker rank $(37,37,3)$ was used.}}\label{baboonsmooth}
\end{center}
\end{figure}

\begin{figure}
\begin{center}
\includegraphics[width=0.9\columnwidth]{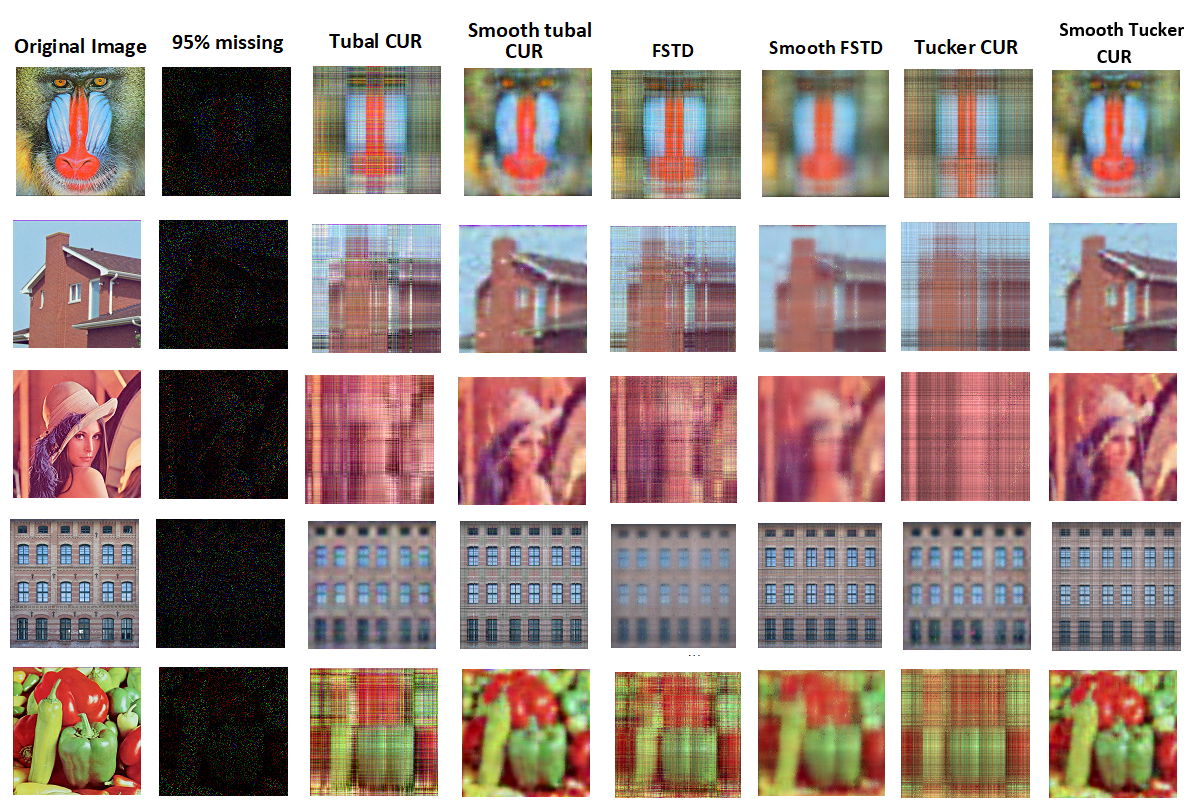}
\caption{\small{The reconstructed images using the Tucker CUR, the tubal CUR and their smooth variants for images with $95\%$ missing pixels. The Tucker rank $(37,37,3)$ and tubal rank $R_1=25,\,R_2=25$ were used.}}\label{95smooth}
\end{center}
\end{figure}

In another simulation, we considered images with more difficult structural missing components where several sequential columns/rows are removed and seen in Figure \ref{baboonsmooth} (Second column). The reconstructed images using the Tucker CUR and the smooth Tucker algorithms with the Tucker rank $(37,37,3)$ are shown in Figure \ref{baboonsmooth}. The superiority of the smooth Tucker CUR over the Tucker CUR algorithm is visible. Recovering images with a high missing ratio is also a difficult and more challenging problem. We evaluated the performance of the Tucker CUR, the FSTD, the tubal CUR, for $95\%$ missing components, we compared the reconstructed images of the CUR algorithms and their smooth variants using the Tucker rank $(37,37,3)$ and tubal rank $(25,25)$ are reported in Figure \ref{95smooth} and Table \ref{tab:my_label}. The results show that the proposed smooth Tucker CUR algorithm is able to recover images $95\%$ with missing components within a few seconds. To the best of our knowledge our algorithms is the fastest algorithm to recover such difficult images.

\begin{table*}\caption{The running time comparison (in second) of different completion algorithms.}\label{tab:my_time}
\begin{center}
\footnotesize
\begin{tabular}{||c |c| c| c| c| c||} 
 \hline
 Image & TR-WOPT & TR-ALS &  SPC & FAST MDT & Smooth Tucker CUR\\ [0.5ex] 
 \hline\hline
 Lena &  68.35 & 97.29 & 30.15 & 6.53 & {\bf 1.56}\\ 
 \hline
 House &  65.71 & 102.15  & 28.85 & 6.23 & {\bf 1.37} \\
 \hline
 Facade &  63.45 & 14.17  & 38.95 & 6.23 & {\bf 1.64}\\
 \hline
Peppers  & 67.14 & 260.68  & 28.18 & 7.25 & {\bf 1.58} \\
 \hline
 Baboon & 70.95 & 79.68  &  29.30 & 7.21 & {\bf 1.60}\\  
 \hline
\end{tabular}
\end{center}
\end{table*}

\begin{table*}\caption{The PSNR and SSIM of the initial incomplete images and those achieved by the smooth Tucker CUR, smooth tubal CUR, and the smooth FSTD with corresponding running time (in second) for images with $95\%$ missing components. The Tucker rank $(37,37,3)$ and tubal rank $R_1=25,\,R_2=25$ were used.}\label{tab:my_label}
\begin{center}
\footnotesize
\begin{tabular}{||c c c c c||} 
 \hline
 Image & Incomplete image & Smooth FSTD &  Smooth Tubal CUR & Smooth Tucker CUR\\ 
  & [PSNR,\,SSIM] & [PSNR,\,SSIM,\,Time]  & [PSNR,\,SSIM,\,Time]  & [PSNR,\,SSIM,\,Time]\\ [0.5ex] 
 \hline\hline
 Lena & [5.3439,\,\,0.0077] & [19.0626,\,0.8614,\,3.22]  & [21.4871,\,\,0.8919,\,\,5.97] &[{\bf 22.1557},\,\,{\bf 0.9011},\,\,4.11] \\ 
 \hline
 House & [4.8375,\,\,0.0042] & [19.6538,\,0.7140,\,2.08] & [21.9302\,\,0.7923,\,\,5.26] & [{\bf 22.6334},\,\,{\bf 0.8168},\,\,3.15]\\
 \hline
 Facade & [5.9693,\,\,0.0045] & [20.0580,\,\,0.4271,\,\,3.22] & [{ 22.3742},\,\,{ 0.6061},\,\,4.23] & [{\bf 22.5591},\,\,{\bf 0.6445},\,\,3.57]\\
 \hline
Peppers  & [6.1378,\,\,0.0088
] & [15.3699,\,0.6983,2.76] & [20.0312,\,\,0.8580,\,\,4.620] & [{\bf 20.6701},\,\,{\bf 0.8751},\,\,2.32]\\
 \hline
 Baboon & [5.5815,\,\,0.0064] & [18.1769,\,0.4988,\,2.16] & [19.5935,\,\,0.5624,\,\,5.82] & [{\bf 19.8792},\,\,{\bf 0.5763},\,\,3.38]\\ [1ex] 
 \hline
\end{tabular}
\end{center}
\end{table*}

\begin{figure}
\begin{center}
\includegraphics[width=0.9\columnwidth]{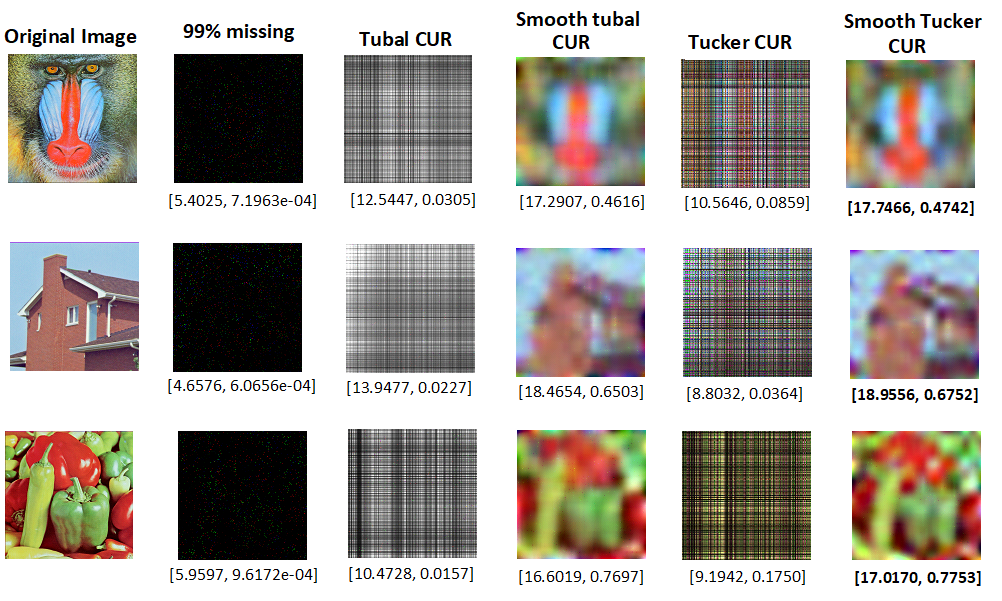}
\caption{\small{The reconstructed images using the Tucker CUR, the tubal CUR and their smooth variants for images with $99\%$ missing pixels. The Tucker rank $(25,25,3)$ and tubal rank $R_1=22,\,R_2=22$ were used.}}\label{99smooth}
\end{center}
\end{figure}

We also conducted simulations with $99\%$ missing pixels and applied Tucker, tubal CUR approximations and their smooth variants. In the case of Tucker CUR, we have used Tucker rank $(25,25,3)$ and for the tubal CUR, we have used tubal rank $(22,22)$. The reconstructed images and corresponding PSNR and SSIM are reported in Figure \ref{99smooth}. The results again show the effectiveness of the idea of smoothness for recovering incomplete images with a high  missing rate.

\begin{figure}
\begin{center}
\includegraphics[width=0.9\columnwidth]{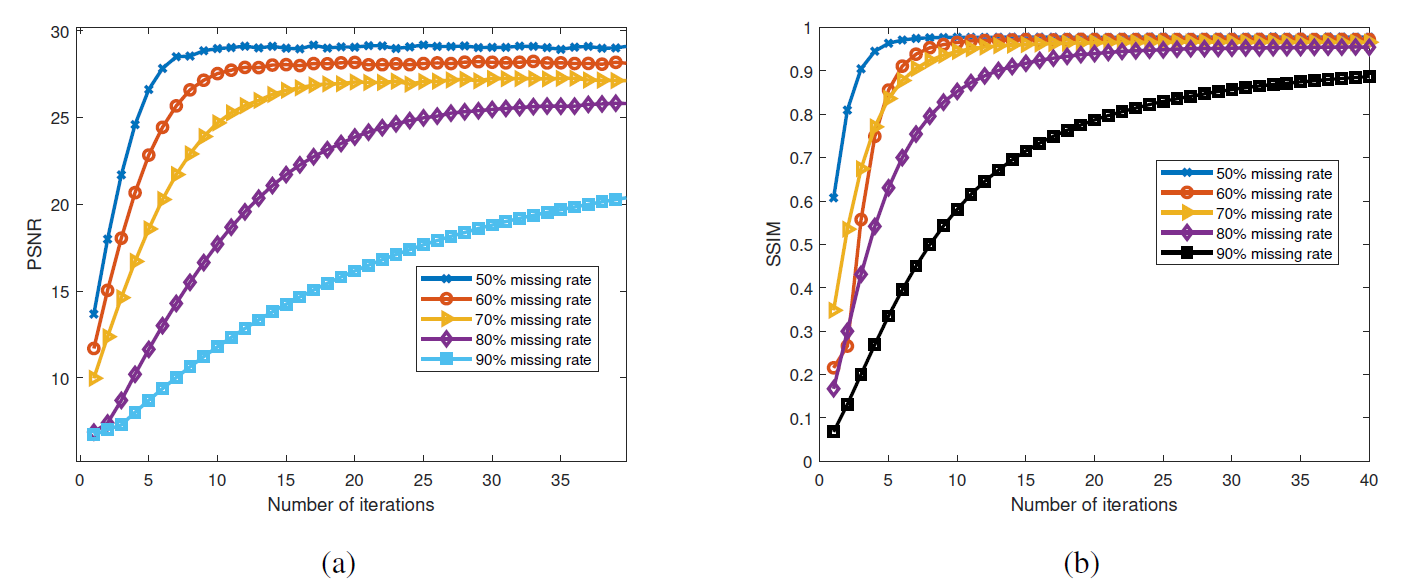}
\caption{\small{The results of the smooth Tucker CUR for reconstructing the pepper image with different missing rates a) The PSNR comparison b) The SSIM comparison. The Tucker rank $(70,70,3)$ was used in the simulations.}}\label{PSNRSSImTucker}
\end{center}
\end{figure}

\begin{figure}
\begin{center}
\includegraphics[width=1\columnwidth]{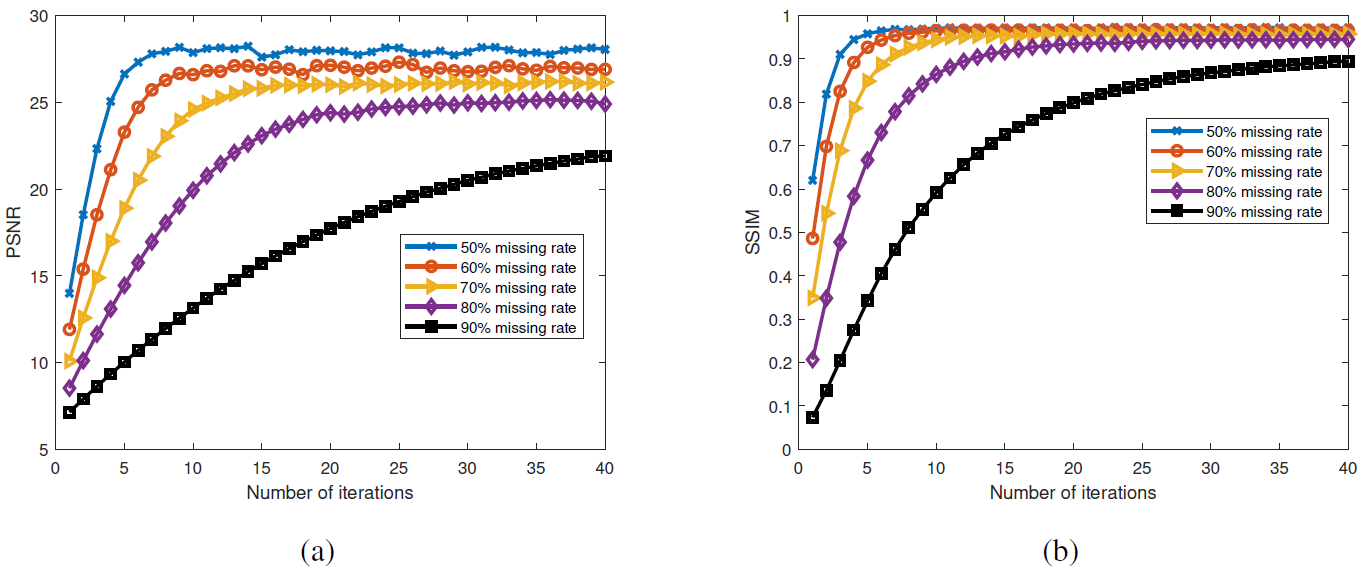}
\caption{\small{The results of the Smooth tubal CUR for the pepper image with different missing rates a) The PSNR comparison b) The SSIM comparison. The tubal rank $(50,50)$ was used.}}\label{PSNRSSImTubal}
\end{center}
\end{figure}

In Figures \ref{PSNRSSImTucker}-\ref{PSNRSSImTubal}, for the "pepper" image, we consider different missing rates and compare the number of iterations against the achieved PSNR and SSIM using the Tucker and tubal decompositions, respectively. 

\begin{figure}
\begin{center}
\includegraphics[width=1\columnwidth]{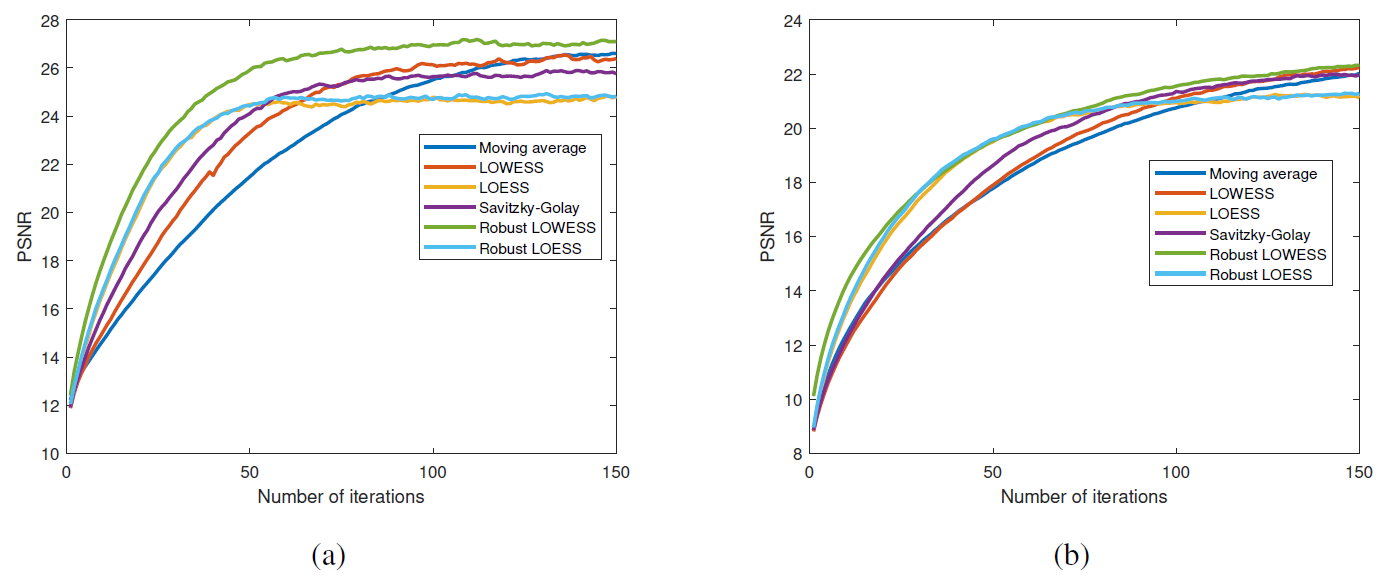}
\caption{\small{Comparing the performance of different smoothing strategies using the Tucker CUR approximation for reconstructing images, a) incompleted pepper in third row of Figure \ref{Compare1} b) incompleted baboon in the last row of Figure \ref{Compare1}. The Tucker rank $(70,70,3)$ was used.}}\label{smoothTucker}
\end{center}
\end{figure}

\begin{figure}
\begin{center}
\includegraphics[width=1\columnwidth]{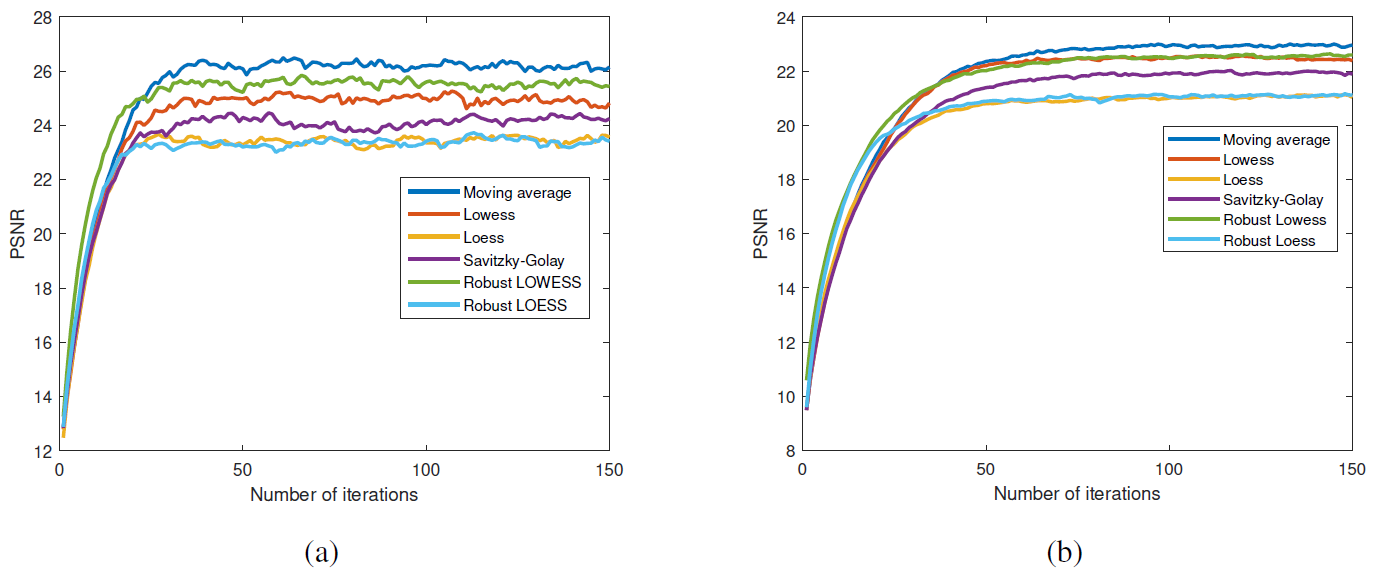}
\caption{\small{Comparing the performance of different smoothing strategies for reconstructing images using the tubal CUR approximation, a) incompleted pepper in Figure \ref{Compare1} b) incomplete baboon in Figure \ref{Compare1}. 40 lateral/horizontal slices were used.}}\label{smoothtubal}
\end{center}
\end{figure}

The quality of different smoothing strategies is compared in Figures \ref{smoothTucker} and \ref{smoothtubal}, for the Tucker CUR and the tubal CUR decompositions, respectively. The two images were used, the incompleted pepper in the third row of Figure 8 and the incompleted baboon in the last row of Figure 8. For the Tucker CUR, Tucker rank $(70,70,3)$ and for the tubal CUR, 40 lateral/horizontal slices were sampled.

For the Tucker decomposition and both images, the RLOWESS smoothing provided the best results while and the tubal decomposition, the moving average was the most promising smoother.

{\bf Video completion} In the second experiment, we consider "Akiyo" and "News" grayscale video datasets from\footnote{\url{http://trace.eas.asu.edu/yuv/}} each of which is a 3rd order tensor of size $176\times 144\times 300$. 
In the next step, we remove randomly $70\%$ pixels of the videos and apply the completion algorithms to complete them. Here similar to the images, the performance of algorithms are reported in two figures. For the Tucker CUR, we considered the Tucker rank $(90,90,120)$, for the tubal CUR we have selected 40 lateral and 40 horizontal slices and for the Slice-tube CUR algorithm, we have selected 80 frontal slices and 1500 tubes.

\begin{figure}
\begin{center}
\includegraphics[width=0.9\columnwidth]{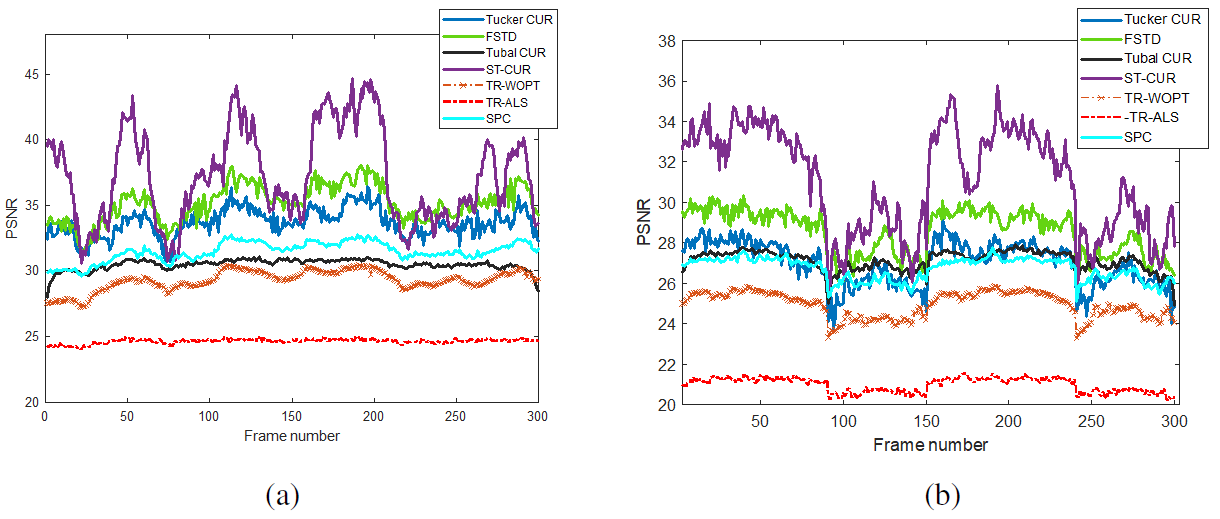}
\caption{\small{The PSNR comparison of different completion algorithms for a) Foreman video dataset and b) News video dataset.}}\label{PSNRcomparisionForeman}
\end{center}
\end{figure}

\begin{figure}
\begin{center}
\includegraphics[width=1\columnwidth]{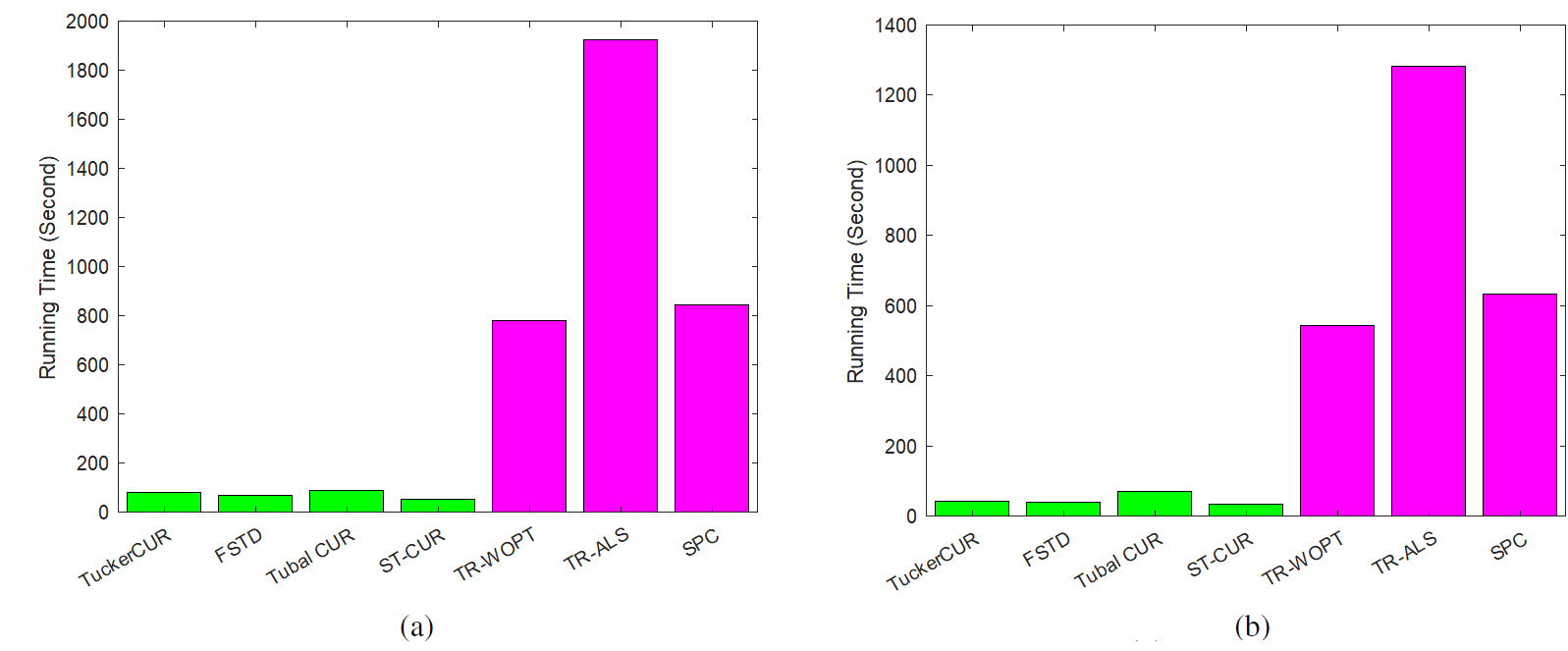}
\caption{\small{The running time comparison of different completion algorithms for a) Foreman video dataset and b) News video dataset.}}\label{RuntimeForeman}
\end{center}
\end{figure}

The PSNR achieved by different completion algorithms are reported in Figure \ref{PSNRcomparisionForeman} (a)-(b) for "Foreman" and "news" Videos, respectively.  
The ST-CUR algorithm outperformed other completion for both video datasets, while other CUR completion algorithms almost worked well. Besides, the running time of all completion algorithms for recovering the Akiyo video with $70\%$ missing components are reported in Figure \ref{RuntimeForeman} (a)-(b). The results show that the Tucker CUR algorithms need less running time for recovering the mentioned video with almost the same or even better performance compared to the other completion algorithms. The results confirmed that the proposed framework is very efficient for reconstructing incomplete video datasets.


{To show the efficiency of the proposed algorithm for higher order data tensors, we used Smoke 1 and Smoke 2 video datasets \cite{yokota2018missing} which are color videos. The Smoke-1 video is a of size $90\times 160\times 3\times 64$ while the Smoke-2 video is of size $90\times 160\times 3\times 100$. We removed randomly $70\%$ of the pixels of both videos and applied the Tucker CUR algorithm with the Tucker rank $(80,80,3,24)$. The reconstruction of some of the frames are displayed in Figures \ref{smoke1} (a) and \ref{smoke2} (a). Moreover the PSNR of all frames are reported in Figures \ref{smoke1} (b) and \ref{smoke2} (b). The time required for the Smoke-1 video was 13.48 seconds while the time for the Smoke-2 video was 17.34 seconds. This clearly shows that the Tucker CUR is a practical algorithm for real-time video reconstruction.}

\begin{figure}
\begin{center}
\includegraphics[width=1\columnwidth]{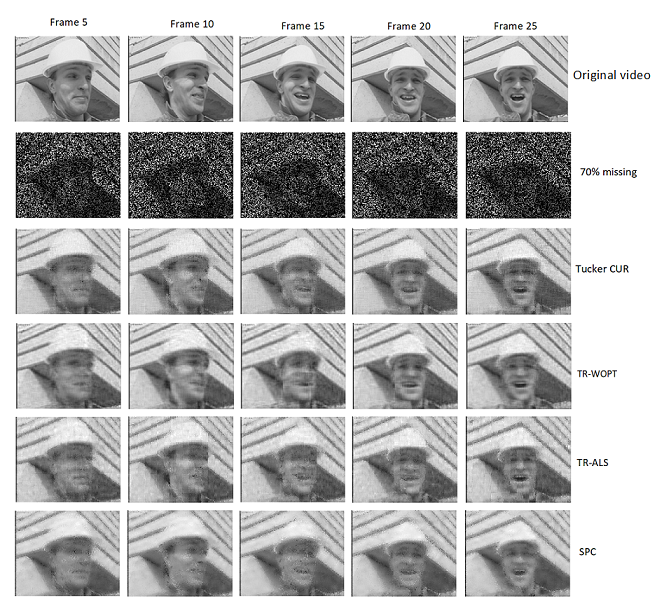}
\caption{\small{The reconstructed frames of the Foreman video dataset with $70\%$ missing pixels using different completion algorithms.}}\label{Comparevideoforeman}
\end{center}
\end{figure}

\begin{figure}
\begin{center}
\includegraphics[width=1\columnwidth]{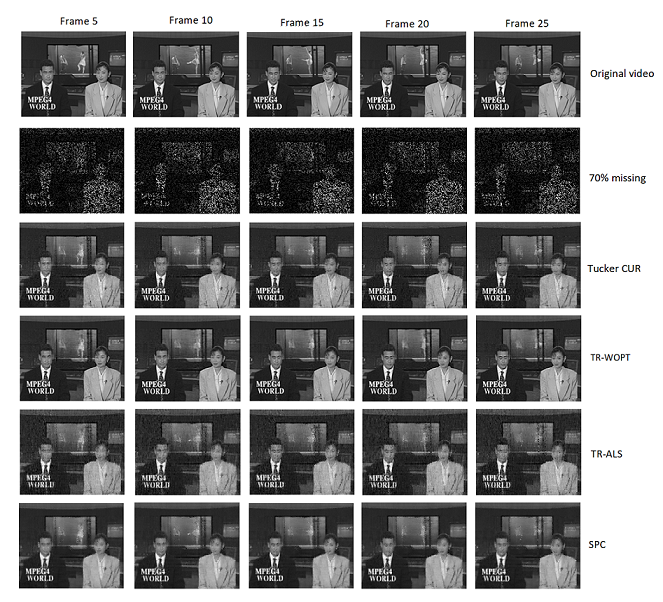}
\caption{\small{The reconstructed frames of the News video dataset using different completion algorithms.}}\label{NewsCompare}
\end{center}
\end{figure}

\begin{figure}
\begin{center}
\includegraphics[width=1\columnwidth]{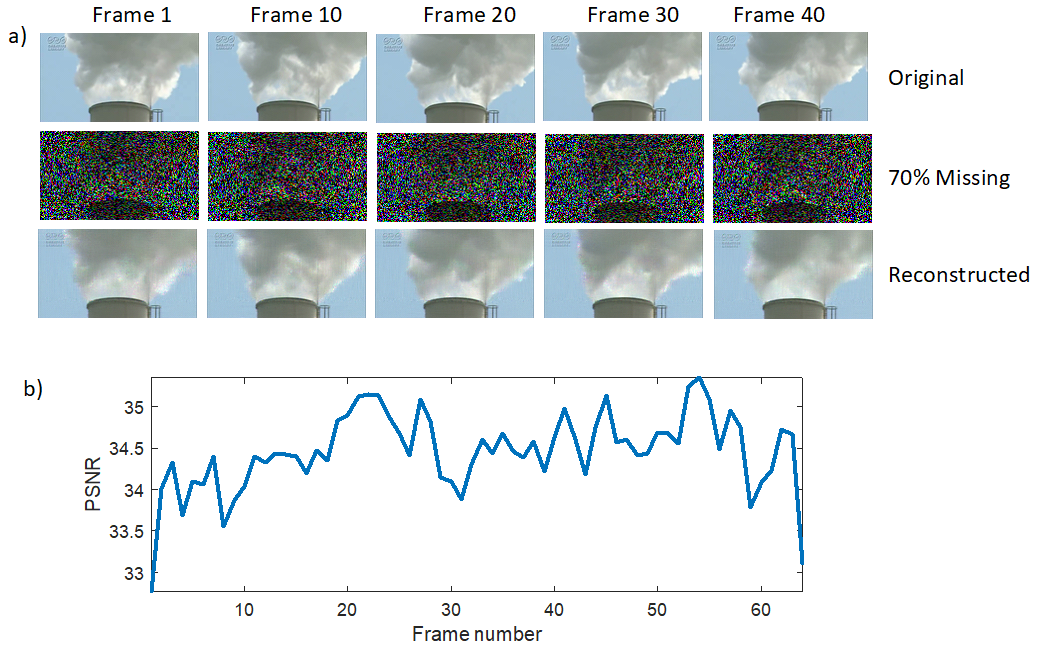}
\caption{\small{(a) The original, the reconstructed and the missing frames (with $70\%$ missing pixels) using the Tucker CUR algorithm for the Smoke-1 video dataset, (b) The PSNR of all reconstructed frames of the Smoke-1 video using the Tucker CUR algorithm.}}\label{smoke1}
\end{center}
\end{figure}

\begin{figure}
\begin{center}
\includegraphics[width=1\columnwidth]{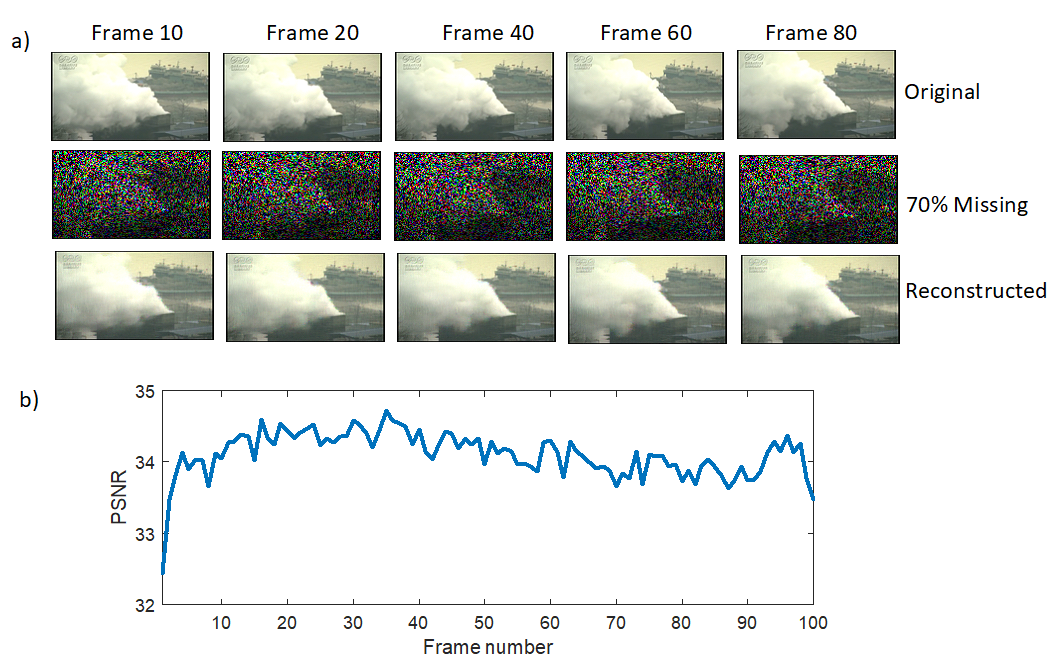}
\caption{\small{(a) The original, the reconstructed and the missing frames (with $70\%$ missing pixels) using the Tucker CUR algorithm for the Smoke-2 video dataset, (b) The PSNR of all reconstructed frames of the Smoke-2 video using the Tucker CUR algorithm.}}\label{smoke2}
\end{center}
\end{figure}

\begin{figure}
     \centering
     \begin{subfigure}[b]{0.95\textwidth}
         \centering
         \includegraphics[width=\textwidth]{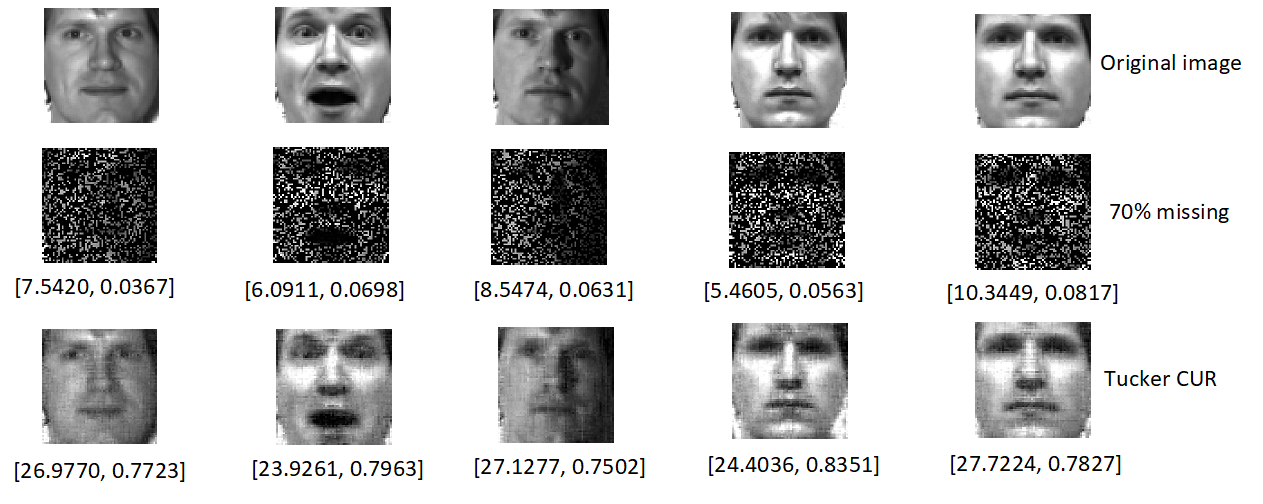}
         \caption{}
         
     \end{subfigure}
      ~
     \begin{subfigure}[b]{0.96\textwidth}
         \centering
         \includegraphics[width=\textwidth]{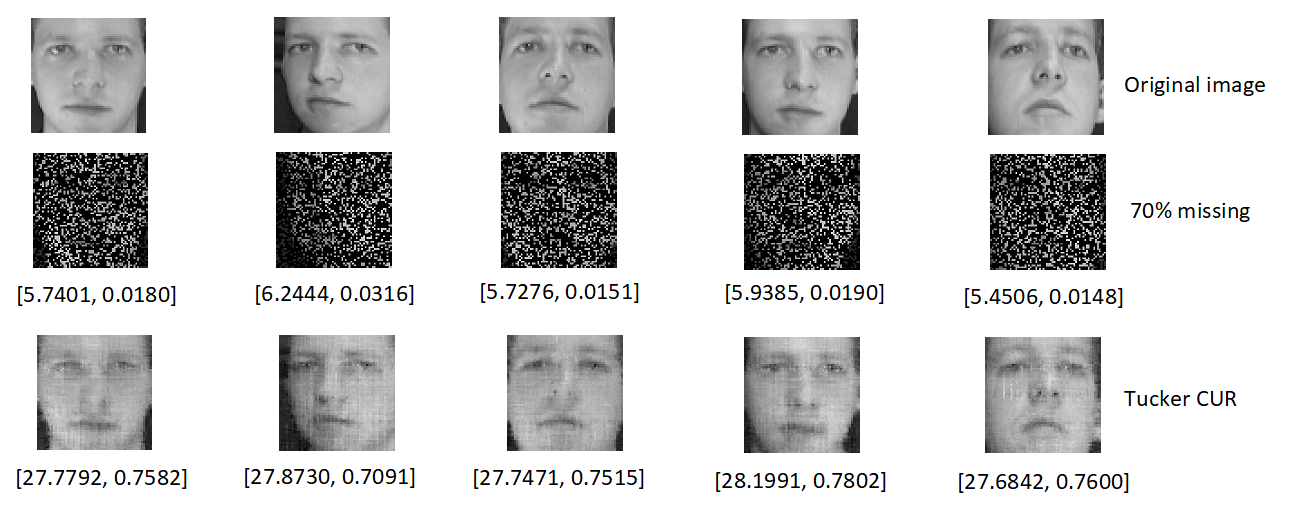}
         \caption{}
     \end{subfigure}
     \caption{The reconstructed images using Tucker CUR algorithm for the a) Yale dataset b) ORL dataset with $70\%$ missing components.} \label{YaleORLrecon}
\end{figure}

\begin{figure}
\begin{center}
\includegraphics[width=1\columnwidth]{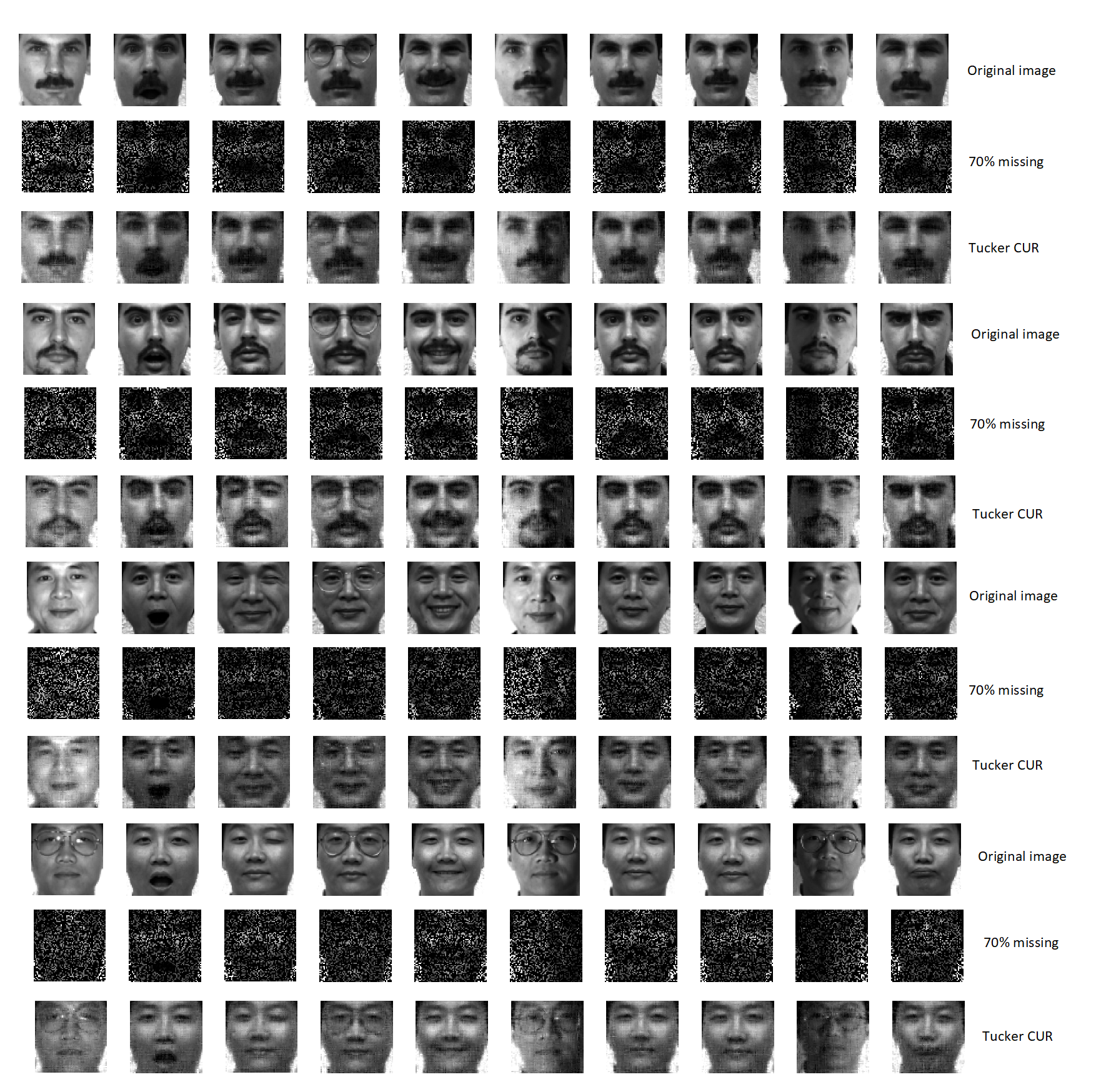}
\caption{\small{The reconstructed images of some samples of the Yale dataset under different expressions using the Tucker CUR algorithm.}}\label{CompareYale2}
\end{center}
\end{figure}

\begin{figure}
\begin{center}
\includegraphics[width=1\columnwidth]{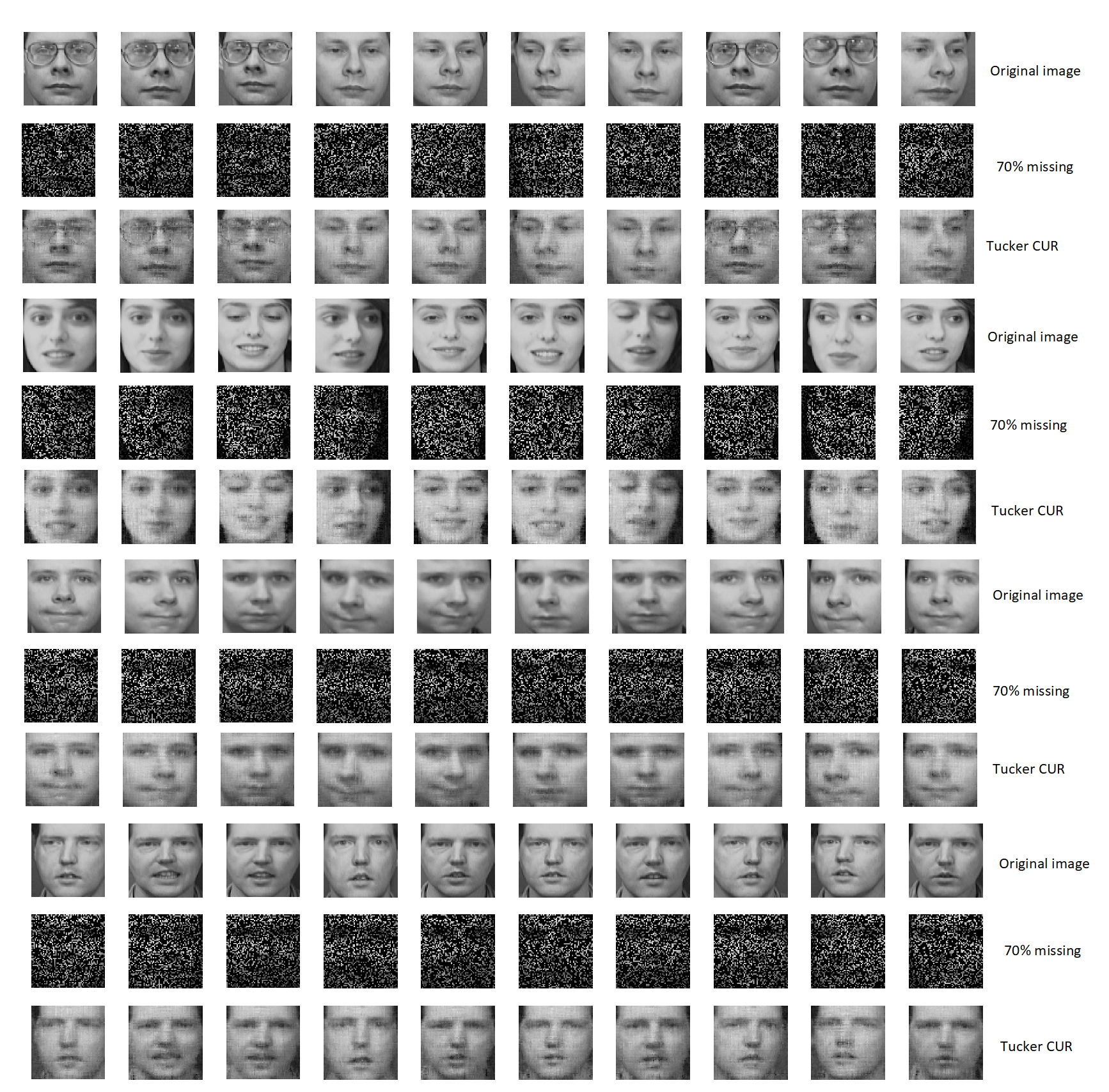}
\caption{\small{The reconstructed images of some samples of the ORAL dataset under different expressions using the Tucker CUR algorithm.}}\label{CompareOral2}
\end{center}
\end{figure}

{\bf Yale and ORAL datasets.} In the previous examples, the superiority of our CUR framework was visible. Now, in the rest of our experiments, we only consider the Tucker CUR algorithm, which almost was the cheapest and more efficient algorithm for the tensor completion task compared to the other tensor CUR algorithms. In the third experiment, we considered Yale and ORL datasets\footnote{\url{http://www.cad.zju.edu.cn/home/dengcai/Data/FaceData.html}} which have been extensively utilized for the face recognition task \cite{CHHH07,CHH07b,CHHZ06,HYHNZ05}. The Yale and ORL detests are a collection of images of size $64\times 64$. For the Yale dataset, there are 15 persons under 11 different expressions e.g., sad, sleepy, surprised, and as a result a fourth order tensor of size $11\times15\times64\times 64$ is produced while for the ORL dataset, 40 distinct subjects for each subject ten different images exist, and a fourth order dataset of size $10\times 40\times 64\times 64$ is generated. We uniformly remove $70\%$ pixels of both mentioned datasets and apply the Tucker CUR algorithms to reconstruct the  whole datasets. For the Yale dataset we considered $R_1=11,R_2=15,R_3=50,R_4=50$ and for the ORL we used $R_1=10,R_2=40,R_3=45,R_4=45$. We first report the accuracy of some of the reconstructed images in Figure \ref{YaleORLrecon}. Besides, in Figures \ref{CompareYale2} and \ref{CompareOral2}, more reconstructed images are displayed. The results show the efficiency and applicability of the proposed algorithm for the case we want to reconstruct a set of images instead of only one single image. 

\begin{figure}
\begin{center}
\includegraphics[width=1\columnwidth]{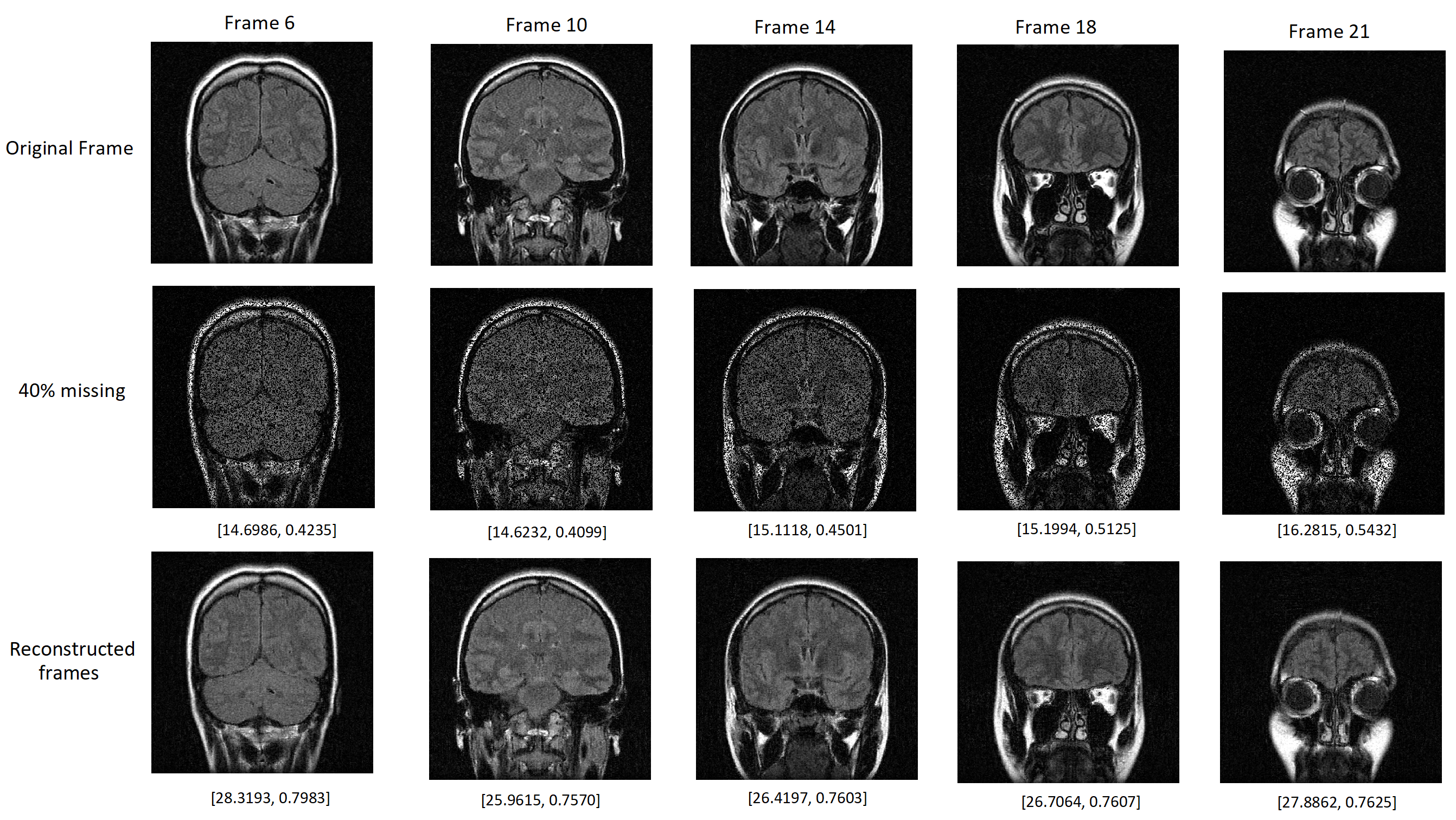}
\caption{\small{The reconstructed frames (number 6-10-14-18-21) using the Tucker CUR algorithm for the MRI brain dataset.}}\label{MRI}
\end{center}
\end{figure}

\begin{figure}
\begin{center}
\includegraphics[width=1\columnwidth]{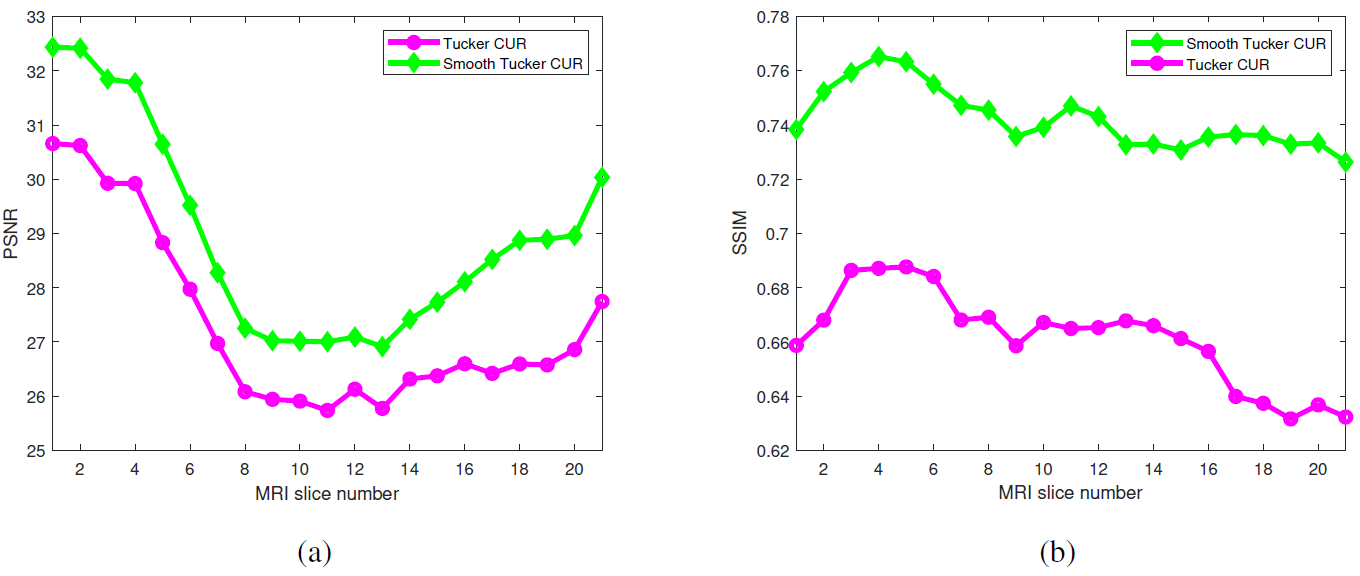}
\caption{\small{The comparison of Tucker CUR and Smooth Tucker CUR for the MRI dataset a) The PSNR comparison b) The SSIM comparison.}}\label{PSNRSSIM_MRI}
\end{center}
\end{figure}

{\bf MRI data set.} In the fourth experiment, we show that our proposed Tucker CUR algorithm is also applicable for medical images. To this end, we applied the Tucker CUR completion algorithm developed by us for the MRI dataset. The original dataset is a 3rd order tensor of size $256\times 256\times 21$, see Figure \ref{MRI} for some sample slices of the data tensor. We uniformly removed $40\%$ of the pixels of the original data tensor and We applied the Tucker CUR algorithm with $R_1=R_2=160,\,\,R_3=21$. The reconstructed images (also images with missing pixels) and their corresponding PSNR and SSIM are reported in Figure \ref{MRI}. We also apply  smoothing the fibers in the first and second modes and the PSNR and SSIM comparisons are reported in Figure \ref{PSNRSSIM_MRI}.

The results show that the proposed algorithm, especially the algorithm equipped with smoothing idea is a promising approach for completing the medical images.

\section{Conclusion}\label{Sec:Conclu}
In this paper we proposed a general framework of tensor ColUmn-Row (CUR) approximation algorithms for the tensor completion task. The smooth variants of the suggested CUR completion algorithms were also proposed to improve the results. We experimentally showed that the smooth algorithms tackle the problem of reconstructing data tensor with a high missing ratio. The proposed algorithms are simple and easy to be implemented only in a few lines of code. The algorithms have low computational complexity because of using CUR approximation within them instead of more expensive ones. The efficiency and performance of the proposed algorithms are substantiated by extensive simulations on images/videos. Our future work is using our proposed CUR framework combined with a Hankelization step as a pre-processing to further improve the performance of the results.   

\section*{Acknowledgement}
The work was partially supported by the Ministry of Education and Science of the Russian Federation (grant 075.10.2021.068).

\begin{appendices}
\section{Tubal operations}
The tensor CUR approximation can be defined based on t-product \cite{kilmer2011factorization,kilmer2013third}, and due to this, we introduce some concepts and definitions related to this concept.  
\begin{defn} \cite{kilmer2013third} ({t-product})
Let $\underline{\mathbf X}\in\mathbb{R}^{I_1\times I_2\times I_3}$ and $\underline{\mathbf Y}\in\mathbb{R}^{I_2\times I_4\times I_3}$, the t-product $\underline{\mathbf X}*\underline{\mathbf Y}\in\mathbb{R}^{I_1\times I_4\times I_3}$ is defined as follows
\begin{equation}\label{TPROD}
\underline{\mathbf C} = \underline{\mathbf X} * \underline{\mathbf Y} = {\rm fold}\left( {{\rm circ}\left( \underline{\mathbf X} \right){\rm unfold}\left( \underline{\mathbf Y} \right)} \right),
\end{equation}
where 
\[
{\rm circ} \left(\underline{\mathbf X}\right)
=
\begin{bmatrix}
\underline{\mathbf X}(:,:,1) & \underline{\mathbf X}(:,:,I_3) & \cdots & \underline{\mathbf X}(:,:,2)\\
\underline{\mathbf X}(:,:,2) & \underline{\mathbf X}(:,:,1) & \cdots & \underline{\mathbf X}(:,:,3)\\
 \vdots & \vdots & \ddots &  \vdots \\
 \underline{\mathbf X}(:,:,I_3) & \underline{\mathbf X}(:,:,I_3-1) & \cdots & \underline{\mathbf X}(:,:,1)
\end{bmatrix},
\]
and 
\[
{\rm unfold}(\underline{\mathbf Y})=
\begin{bmatrix}
\underline{\mathbf Y}(:,:,1)\\
\underline{\mathbf Y}(:,:,2)\\
\vdots\\
\underline{\mathbf Y}(:,:,I_3)
\end{bmatrix},\hspace*{.5cm}
\underline{\mathbf Y}={\rm fold} \left({\rm unfold}\left(\underline{\mathbf Y}\right)\right).
\]
\end{defn}
In view of \eqref{TPROD}, it is seen that the t-product operation is the circular convolution operator, and because of this, it can be easily computed through Fast Fourier Transform (FFT). More precisely, we first transform all tubes of two tensors $\underline{\bf X},\,\underline{\bf Y}$, into the frequency domain, then multiply frontal slices of the spectral tensors. Finally, we apply the Inverse FFT (IFFT) to all tubes of the last tensor. This procedure is summarized in Algorithm 2. Note that other types of tensor decompositions in the t-product format can be computed in a similar manner. For example, for computing the tubal QR computation of a tensor $\underline{\mathbf X}\in\mathbb{R}^{I_1\times I_2\times I_3},\,$ i.e., $\underline{\mathbf X} = \underline{\mathbf Q} * \underline{\mathbf R}$, we first compute the FFT of the tensor $\underline{\mathbf X}$ as \cite{kilmer2013third}
\begin{equation}\label{fftten}
\underline{\widehat{\mathbf X}}={\rm fft}({\underline{\mathbf X}},[],3),
\end{equation}
and then the QR decomposition of all frontal slices of the tensor $\underline{\widehat{\mathbf X}}$ is computed as follows
\[
\underline{\widehat{\mathbf X}}(:,:,i)=\underline{\widehat{\mathbf Q}}(:,:,i)\,\,\underline{\widehat{\mathbf R}}(:,:,i).
\]
Note that \eqref{fftten} is equivalent to computing the FFT of all tubes of the tensor $\underline{\bf X}$. Finally, the IFFT operator is applied to the tensors $\underline{\widehat{\mathbf Q}}$ and $\underline{\widehat{\mathbf R}}$, to compute the tensors $\underline{{\mathbf Q}}$ and $\underline{{\mathbf R}}$. 

\begin{figure}
    \centering
    \includegraphics[width=1\columnwidth]{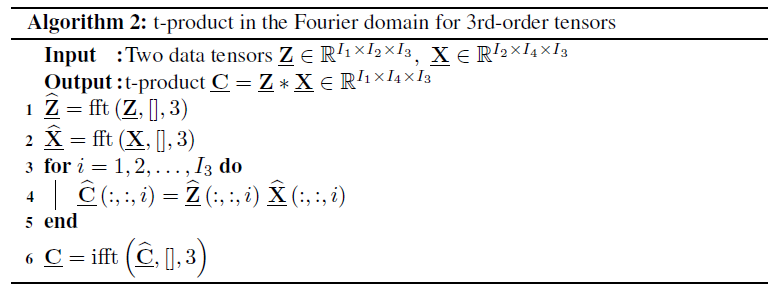}
\end{figure}

\begin{defn}\cite{kilmer2013third} ({Transpose})
Let $\underline{\mathbf X}\in\mathbb{R}^{I_1\times I_2\times I_3}$ be a given tensor. Then the transpose of tensor $\underline{\mathbf X}$ is denoted by $\underline{\mathbf X}^{T}\in\mathbb{R}^{I_2\times I_1\times I_3}$, which is constructed by applying transpose to all its frontal slices and reversing the order of second till last transposed frontal slices.
\end{defn}

\begin{defn} Assume $\underline{\mathbf{X}}\in\mathbb{R}^{I_1\times I_2\times I_3}$, then the MP of the tensor $\underline{\mathbf{X}}$ is denoted by $\underline{\mathbf{X}}^{+}\in\mathbb{R}^{I_2\times I_1\times I_3}$ and defined as a unique tensor satisfying the next four relations
\begin{eqnarray*}
\underline{\mathbf X}*{\underline{\mathbf X}^{+} }*\underline{\mathbf X} &=& \underline{\mathbf X},\quad\quad\quad {\underline{\mathbf X}^{+} }*\underline{\mathbf X}*{\underline{\mathbf X}^{+} } ={\underline{\mathbf X}^{+} },\\
{\left( {\underline{\mathbf X}*{\underline{\mathbf X}^{+} }} \right)^T} &=& \underline{\mathbf X}*{\underline{\mathbf X}^{+} },\,\,\,\,\,\,\,{\left( {{\underline{\mathbf X}^{+} }*\underline{\mathbf X}} \right)^T} = {\underline{\mathbf X}^{+} }*\underline{\mathbf X}.
\end{eqnarray*}

\end{defn}

Similar to other tensor operations, the MP pseudoinverse of tensors can be computed through FFT and this is described in Algorithm 3.
\begin{figure}
    \centering
    \includegraphics[width=1\columnwidth]{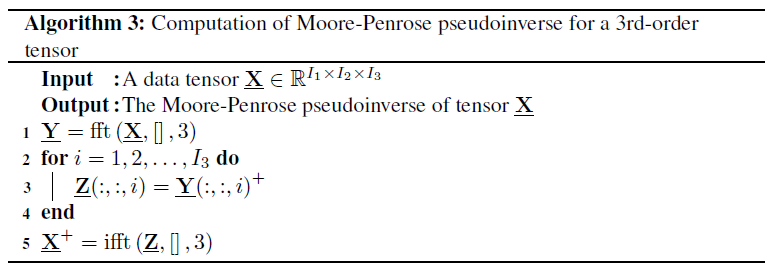}
\end{figure}
So, for computing the MP of tensor, it is first converted to the frequency domain, and then the MP of all frontal slices is computed. Finally, the resulting tensor is converted to the original space using inverse FFT.

In this part, we provide a comprehensive materials about different cross tensor approximation methods and algorithms. For more details we refer to our recent publication \cite{ahmadi2021cross}.
\section{CTA Based on Fiber Selection}
The CMA can be extended to tensors  based on the Tucker-2 model, see Figure \ref{Tucker2}. Here, some columns and rows $\mathbf{C},\,\mathbf{R}$ are selected. To be more precise, matrices are chosen and the middle core tensor $\underline{\mathbf{U}}$ should be computed to yield the smallest error. The optimal option for the middle core tensor is
 \[
 \underline{\mathbf U} = \underline{\mathbf X}{ \times _1}{\mathbf C}^{+}{ \times _2}{\mathbf R}^{+}.
 \]
The Tucker-2 model in suitable in situation that the 3rd mode is relatively small and reduction in that mode not required. For example, in the case of color images, the third mode is related to three colors (green-blue-red) and only columns/rows are selected. In the experimental results we show this procedure.
Natural and straightforward generalizations of the CMA to tensors are proposed in \cite{drinea2001randomized,oseledets2008tucker}, and \cite{caiafa2010generalizing}. Similar to the CMA in which a part of columns and rows of a given matrix is sampled, a set of fibers (along different modes) are selected, and the goal is to compute a Tucker approximation based on these sampled fibers. For instance, for 3rd-order tensors, we should sample columns, rows, and tubes. The approach proposed in \cite{drineas2007randomized} is randomized, while those proposed in \cite{oseledets2008tucker,caiafa2010generalizing} are deterministic. These works can be considered as the starting points of generalization of the CMA to tensors, and in the sequel, we explain the idea behind each of these Tucker approximations. 
For noiseless data tensor, the existence of an exact Tucker model whose factor matrices are taken from the fibers of the original data tensor is obvious. To be more precise, let $\underline{\mathbf X}$ be an $N$th-order tensor of size ${I_1\times I_2\times \cdots \times I_N}$ and of Tucker rank $(R_1,R_2,\ldots,R_N)$. Now, if we generate any full-rank factor matrices ${\mathbf A}_{n}\in\mathbb{R}^{I_n\times R_n},\,n=1,2,\ldots,N$, by sampling fibers in each mode and computing the core tensor as
\begin{eqnarray}\label{CorInter}
\underline{\mathbf S} = \underline{\mathbf X}{ \times _1}{\mathbf A}_1^{+} { \times _2}{\mathbf A}_2^{+}\cdots { \times _N} {\mathbf A}_N^{+}\in\mathbb{R}^{R_1\times R_2\times \cdots\times R_N},
\end{eqnarray}
then the obtained Tucker decomposition has the exact Tucker rank $(R_1,R_2,\ldots,R_N)$. So an exact Tucker decomposition of the tensor $\underline{\mathbf{X}}$ whose factor matrices are taken from the original data tensor is computed. For noisy tensors, similar to the CMA, the accuracy of approximation quite depends on the number of selected fibers and also the list of sampled fibers. The idea of sampling fibers and considering them as the factor matrices, first was proposed in \cite{drineas2007randomized}. In the first step, the factor matrices are generated, after which the core tensor is computed through \eqref{CorInter} as described above. This is summarized in Algorithm 4. In Algorithm 4. a shorthand notation is used for the Tucker decomposition as
\[
\underline{\mathbf X} = \left[\kern-0.15em\left[ {\underline{\mathbf S};{{\mathbf A}_{1}},{{\mathbf A}_{2}}, \ldots ,{{\mathbf A}_{N}}}
 \right]\kern-0.15em\right].
\]
\begin{figure}
    \centering
    \includegraphics[width=1\columnwidth]{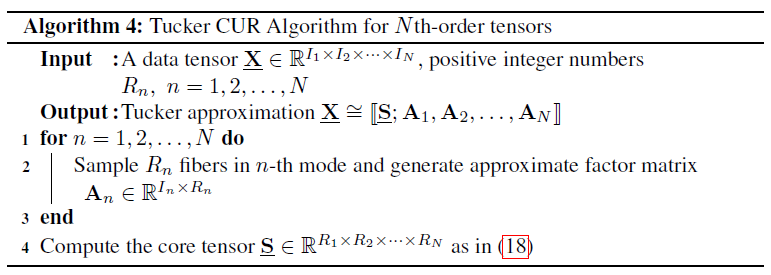}
\end{figure}
\begin{figure}
\begin{center}
\includegraphics[width=8.5 cm,height=4.3 cm]{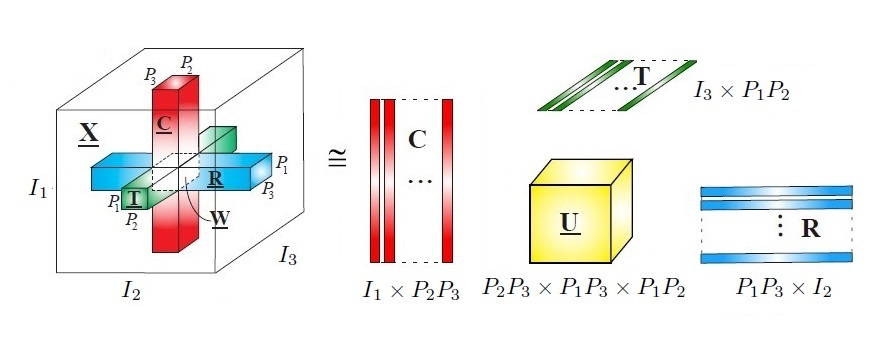}\\
\caption{\small{Illustration of the CTA (fiber selection version) for a 3rd-order low-rank tensor. For simplicity of presentation, we assume that all fibers build up to block sub-tensors, \cite{cichocki2016tensor}}.}\label{CTD}
\end{center}
\end{figure}
\begin{figure}
\begin{center}
\includegraphics[width=8 cm,height=2 cm]{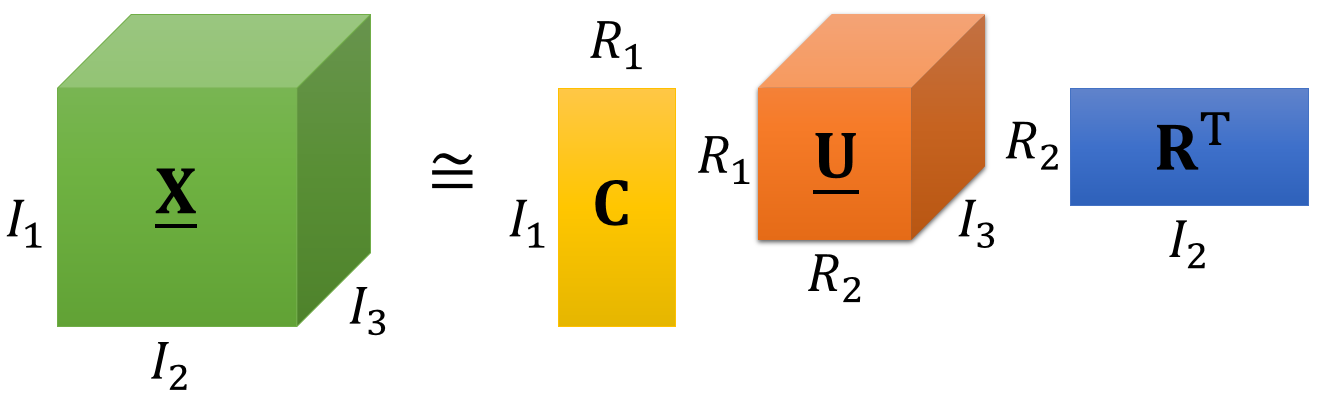}\\
\caption{\small{llustration of the Tucker-2 CUR approximation.}}\label{Tucker2}
\end{center}
\end{figure}

The question is, how to compute an approximate Tucker decomposition for the
tensor $\underline{\mathbf X}$ based on the intersection subtensor $\underline{\mathbf W}$? Motivated by the fact that 
\begin{eqnarray}\label{Motiv}
{\mathbf U}={\mathbf W}{ \times _1}{\mathbf W}_{\left( 1 \right)}^{+} { \times _2}{\mathbf W}_{\left( 2 \right)}^{+}  = {{\mathbf W}^{+} }{\mathbf W}{{\mathbf W}^{+} } = {\mathbf W}^{+} ,
\end{eqnarray}
which is used as the middle matrix in the CMA, it is suggested \cite{caiafa2010generalizing} to compute the approximate core tensor in the Tucker decomposition as 
\begin{eqnarray}\label{Intertensor}
\nonumber
\underline{\mathbf U} &=& \underline{\mathbf W}{ \times _1}{\mathbf W}_{(1)}^{+} { \times _2}{\mathbf W}_{(2)}^{+} { \times _3}{\mathbf W}_{(3)}^{+}\\
&\equiv & \left[\kern-0.15em\left[ {\underline{\mathbf W};{\mathbf W}_{\left( 1 \right)}^{+} ,{\mathbf W}_{\left( 2 \right)}^{+} ,{\mathbf W}_{\left( 3 \right)}^{+} } 
 \right]\kern-0.15em\right].
\end{eqnarray}
This is a direct generalization of \eqref{Motiv} to 3rd-order tensors. In view of \eqref{Intertensor}, the core tensor $\underline{\mathbf U}$ of the Tucker approximation is of size $P_2P_3\times P_1P_3 \times P_1P_2$ and as a result, we need to sample $P_2P_3$ columns, $P_1P_3$ rows, and $P_1P_2$ tubes, see Figure \ref{CTD}.
They should be selected in an appropriate way. It is shown in \cite{caiafa2010generalizing} that the corresponding factor matrices ${{\mathbf A}_1} \in {\mathbb{R}^{I_1 \times {P_2}{P_3}}},\,\,{{\mathbf A}_2} \in {\mathbb{R}^{I_2 \times {P_1}{P_3}}},{{\mathbf A}_3} \in {\mathbb{R}^{I_3 \times {P_1}{P_2}}}$ are the subsampled matrices from the unfolding matrices ${\mathbf X}_{(1)}(:,\mathcal{I}_2,\mathcal{I}_3),\,\,{\mathbf X}_{(2)}(\mathcal{I}_1,:,\mathcal{I}_3)$ and ${\mathbf X}_{(3)}(\mathcal{I}_1,\mathcal{I}_2,:)$ respectively and CTA approximation can be found as
\begin{eqnarray}\label{TuckCUR}
\nonumber
\underline{\mathbf X} &\cong & \left[\kern-0.15em\left[ {\underline{\mathbf U};{{\mathbf A}_1},{{\mathbf A}_2},{{\mathbf A}_3}} 
 \right]\kern-0.15em\right]\\ 
 &\equiv & 
\left[\kern-0.15em\left[ {\underline{\mathbf W};\underbrace {{{\mathbf A}_1}{\mathbf W}_{\left( 1 \right)}^{+} }_{{\widetilde{\mathbf C}_1}},\underbrace {{{\mathbf A}_2}{\mathbf W}_{\left( 2 \right)}^{+} }_{{\widetilde{\mathbf C}_2}},\underbrace {{{\mathbf A}_3}{\mathbf W}_{\left( 3 \right)}^{+} }_{{\widetilde{\mathbf C}_3}}} 
 \right]\kern-0.15em\right].
\end{eqnarray}
The procedure of this approach is summarized as follows:
\begin{itemize}
\item Consider indices $\mathcal{I}_n\in [I_n],\,\,n=1,2,3$ and produce the intersection subtensor $\underline{\mathbf W}$ and corresponding sampled columns, rows, and tubes ${\mathbf A}_1,\,\,{\mathbf A}_2$ and ${\mathbf A}_3$.

\item Compute the Tucker approximation \eqref{TuckCUR}.

\end{itemize}
We refer to this algorithm as Fast Sampling Tucker Decomposition (FSTD) and summarize it in Algorithm 5 \cite{caiafa2010generalizing}.  A similar approach is proposed in \cite{friedland2011fast} in the sense of the number of selected fibers in each mode. 
\begin{figure}
    \centering
    \includegraphics[width=1\columnwidth]{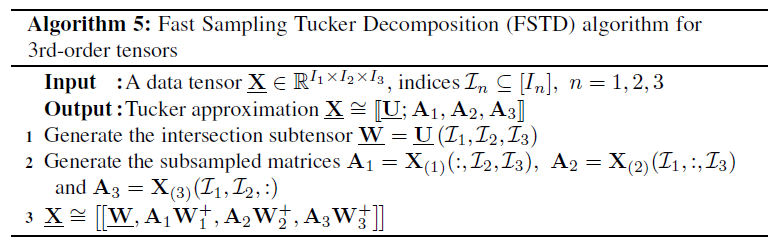}
\end{figure}

\section{CTA Based on Slice-tube Selection}
Motivated by some applications in hyperspectral medical image analysis and consumer recommender system analysis where one of the modes is qualitatively different from others, an alternative CTA is proposed in \cite{mahoney2008tensor}. We first briefly describe this idea for 3rd-order tensors. Let $\underline{\mathbf X}\in\mathbb{R}^{I_1\times I_2\times I_3}$ be a given data tensor, and without loss of generality, we assume that the last mode is qualitatively different from the others.

\begin{figure}
\begin{center}
\includegraphics[width=1\columnwidth]{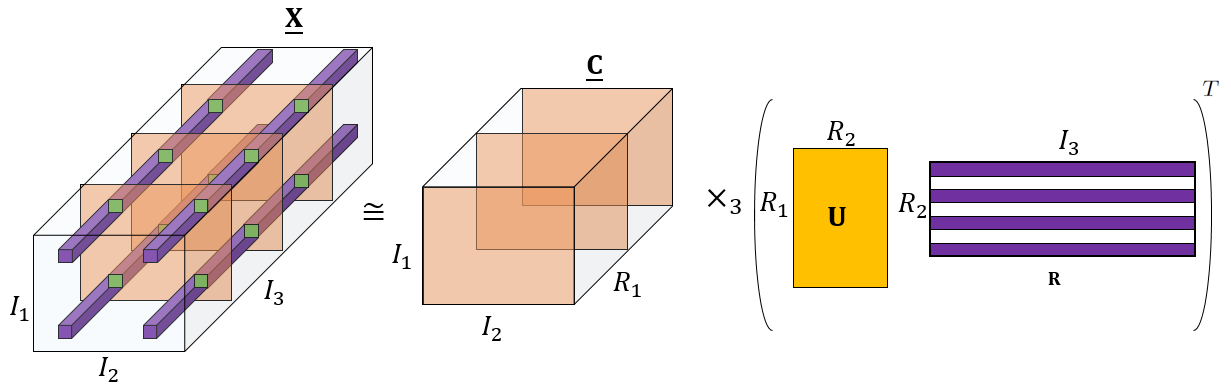}\\
\caption{\small{Illustration of the CTA based on frontal slice and tube selection.}}\label{Pic5}
\end{center}
\end{figure}

Given prior probability distributions for sampling frontal slices as $\{p_i\}_{i=1}^{I_3}$ and tubes as $\{q_{j}\}_{j=1}^{I_1I_2}$, in the first step, some frontal slices, say $L_1$, are sampled and they are stored in $\underline{\mathbf C}\in\mathbb{R}^{I_1\times I_2\times L_1}$. In the second step, we sample some tubes, say $L_2=R_1R_2$, and store them in $\underline{\mathbf R}\in\mathbb{R}^{R_1\times R_2\times I_3}$, or a matrix ${\mathbf R}\in\mathbb{R}^{L_2\times I_3}$, (see Figure \ref{Pic5} (a)). 
The CTA is then defined as (see Figure \ref{Pic5} (b))
\begin{equation}\label{FSTCD}
\underline{\mathbf X} \cong \underline{\mathbf C}{ \times _3}\left( {\mathbf U{\mathbf R}} \right)^T,
\end{equation} 
where the tensor $\underline{\mathbf C}\in\mathbb{R}^{I_1\times I_2\times L_1}$ and the matrix ${\mathbf R}\in\mathbb{R}^{L_2\times I_3}$ contain the sampled frontal slices and tubes, respectively. The matrix ${\mathbf U}\in \mathbb{R}^{L_1\times L_2}$ is defined as
\[
{\mathbf U} = {{\mathbf D}_1}{\left( {{{\mathbf D}_2}{\mathbf W}{{\mathbf D}_1}} \right)^{+}}{{\mathbf D}_2}\in\mathbb{R}^{L_1\times L_2},
\]
\[
{\mathbf W}={\rm reshape}(\underline{\mathbf{W}},[L_2,L_1]),
\]
where ${\mathbf D}_1\in\mathbb{R}^{L_1\times L_1}$ and ${\mathbf D}_2\in\mathbb{R}^{L_2\times L_2}$ are scaling diagonal matrices corresponding to the slice and fiber sampling, respectively and defined as follows
\begin{eqnarray*}
{\left( {{{\mathbf D}_1}} \right)_{tt}} = \frac{1}{{\sqrt {{L_1}{p_{{i_t}}}} }},\,\,\,t = 1,2, \ldots ,{L_1},
\\
{\left( {{{\mathbf D}_2}} \right)_{tt}} = \frac{1}{{\sqrt {{L_2}{q_{{i_t}}}} }},\,\,\,t = 1,2, \ldots ,{L_2},
\end{eqnarray*}
where $\{p_i\}_{i=1}^{I_3}$ are $\{q_j\}_{j=1}^{I_1I_2}$ are probability distributions under which the frontal slices and fibers are sampled. This procedure is summarized in Algorithm 6. The length-squared probability distributions are defined as follows
\begin{eqnarray}\label{LenSquaProb}
\nonumber
p_i&=&\frac{\left\|\underline{\mathbf X}(:,:,i_3)\right\|^2_F}{\left\|\underline{\mathbf{X}}\right\|^2_F},\,\,i_3=1,2,\ldots,I_3,\\
q_j&=&\frac{\underline{\mathbf{X}}(j_1,j_2,:)}{\left\|\underline{\mathbf{X}}\right\|^2_F},\,\,\ j_1,j_2\in J_1,J_2
\end{eqnarray}
where $J_1$ and $J_2$ are subsets of the indices $I_1$ and $I_2$, are used in \cite{mahoney2008tensor} for selecting the slices/tubes.
\begin{figure}
    \centering
    \includegraphics[width=1\columnwidth]{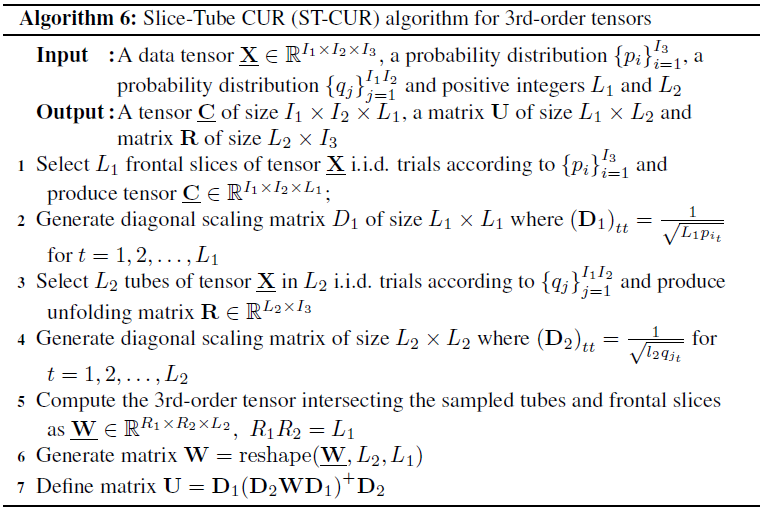}
\end{figure}
\section{CTA Based on Tubal Product (t-Product)}\label{TCDTubal}
\begin{figure}
\begin{center}
\includegraphics[width=8.75 cm,height=3.75 cm]{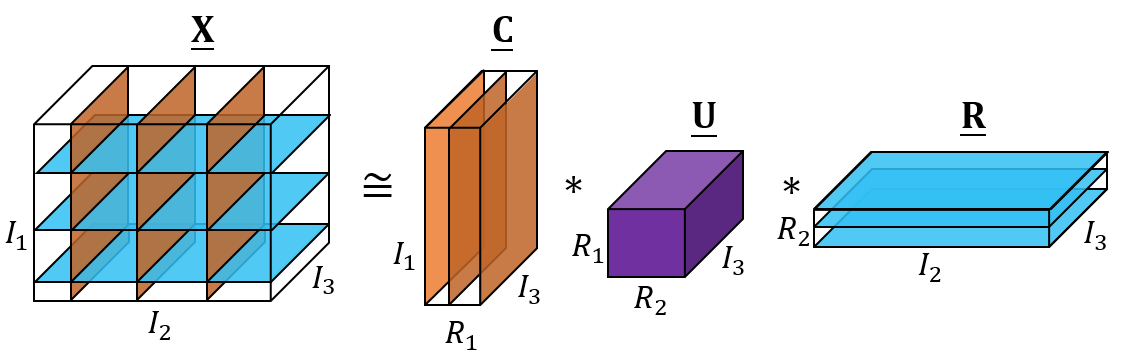}\\
\caption{\small{Illustration of the tubal CTA for a 3rd-order tensor.}}\label{Pic2}
\end{center}
\end{figure}

Here, the tensor variants of the CMA and matrix column selection are called the tubal CTA and lateral slice selection respectively. To be precise, let $\underline{\mathbf X}$ be a given 3rd-order tensor. The tubal cross approximation based on t-product is formulated as follows
\begin{equation}\label{CURT}
\underline{\mathbf X} \cong \underline{\mathbf C} * \underline{\mathbf U} * \underline{\mathbf R},
\end{equation}
where $\underline{\mathbf C}\in\mathbb{R}^{I_1\times L_1\times I_3}$ and $\underline{\mathbf R}\in\mathbb{R}^{L_2\times I_2\times I_3}$ are some sampled lateral and horizontal slices of the original tensor $\underline{\mathbf X}$ respectively and the middle tensor $\underline{\mathbf U}\in\mathbb{R}^{L_1\times L_2\times I_3}$ is computed in such a way that the approximation \eqref{CURT} should be as small as possible, see Figure \ref{Pic2}, for graphical illustration concerning the tubal CTA. 
Note that similar to other CTA algorithms discussed so far, the lateral and horizontal slices sampling procedure can be performed based on prior probability distributions\footnote{Here different various probability distributions (with/without replacement) can be used but we will not go through the theoretical details \cite{tarzanagh2018fast}.}. The probability distributions used in \cite{tarzanagh2018fast} are uniform and nonuniform (length-squared and leverage scores) distributions. 
Let us first consider the tubal CTA. Similar to the CMA, the best solution for the middle tensor $\underline{\mathbf U}$ in the least-squares sense is
\begin{eqnarray}\label{TCUR}
\underline{\mathbf U} = {\underline{\mathbf C}^{+} } * \underline{\mathbf X} * {\underline{\mathbf R}^{+} },
\end{eqnarray}
because
\[
{{\underline{\mathbf C}^{+} } * \underline{\mathbf X} * {\underline{\mathbf R}^{+} }} = \mathop {\rm argmin }\limits_{\underline{\mathbf U} \in {\mathbb{R}^{L \times {I_2} \times {I_3}}}} \,\,{\left\| {\underline{\mathbf X} - \underline{\mathbf C}*\underline{\mathbf U}*\underline{\mathbf U}} \right\|_F},
\]
This is a straightforward generalization of the CMA \cite{goreinov1997theory} to tensors. Formula \eqref{TCUR} can be computed in the Fourier domain and these computations are summarized in Algorithm 7. However, it is clear that \eqref{TCUR} needs to pass the data tensor $\underline{\mathbf X}$ once again and this is of less practical interest for very large-scale data tensors, especially when the data tensors do not fit into the memory and communication between memory and disk is expensive \cite{halko2011finding}. To solve this problem, the MP pseudoinverse of the intersection subtensor $\underline{\mathbf W}\in\mathbb{R}^{L_2\times L_1\times I_3}$, which is obtained based on intersecting the sampled horizontal and lateral slices, should be approximated as
\[
\underline{\mathbf U} = \underline{\mathbf C} * {\underline{\mathbf W}^{+} } * \underline{\mathbf R}.
\]
It is not difficult to see that the tensor $\underline{\mathbf W}$ consists of some tubes of the original data tensor $\underline{\mathbf X}$. 
\begin{figure}
    \centering
    \includegraphics[width=1\columnwidth]{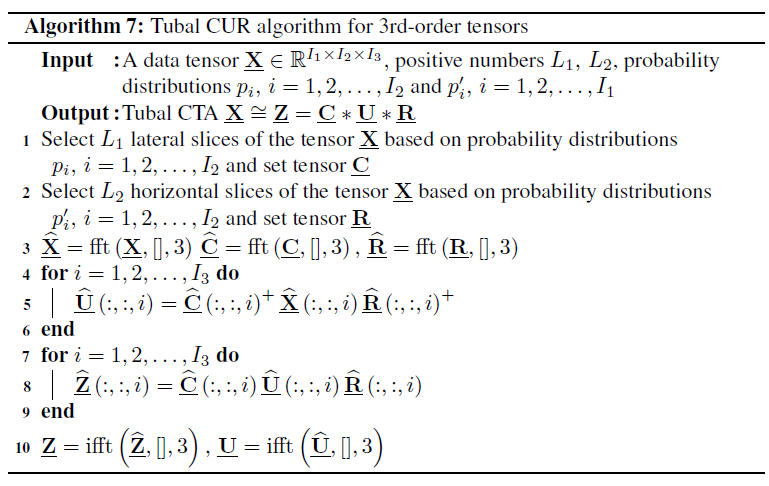}
\end{figure}
Also the tensor lateral slice selection based on the t-product is formulated as follows
\begin{equation}\label{TCX}
\underline{\mathbf X} \cong \underline{\mathbf C} * \underline{\mathbf Y}.
\end{equation}
where $\underline{\mathbf C}\in\mathbb{R}^{I_1\times L\times I_3}$ is a part of lateral slices of the tensor $\underline{\mathbf X}$ and the tensor $\underline{\mathbf Y}\in\mathbb{R}^{L\times I_2\times I_3}$ is computed in such a way that the reconstruction error \eqref{TCX} is as small as possible \cite{tarzanagh2018fast}. The best solution for the tensor $\underline{\mathbf Y}$ is 
\[
\underline{\mathbf Y}={\underline{\mathbf C}^{+}} * \underline{\mathbf X},
\]
which provides the best approximation in a least-squares sense, that is 
\[
{\left\| {\underline{\mathbf X} - \underline{\mathbf C}*\left( {{\underline{\mathbf C}^+ }*\underline{\mathbf X}} \right)} \right\|_{_F}} = \mathop {\min }\limits_{\underline{\mathbf Y} \in {\mathbb{R}^{L \times {I_2} \times {I_3}}}} \,\,{\left\| {\underline{\mathbf X} - \underline{\mathbf C}*\underline{\mathbf Y}} \right\|_F},
\]
through the projection approximation  $\underline{\mathbf X} \cong \underline{\mathbf C} * {\underline{\mathbf C}^{+} } * \underline{\mathbf X}$.

\end{appendices}

\section*{Funding and/or conflicts of interests/competing interests}
The authors have no relevant financial or non-financial interests to disclose.

\section*{Data Availability}
Data sharing is not applicable to this article as no new data were created or analysed in this study.

\section*{Acknowledgements}
The authors would like to thank four reviewers for their constructive and
insightful comments which have greatly improved the quality of the paper. The authors were partially supported by the Ministry of Education and Science of the
Russian Federation (grant 075.10.2021.068).

\bibliography{sn-bibliography}


\begin{thebibliography}{89}
\ifx \bisbn   \undefined \def \bisbn  #1{ISBN #1}\fi
\ifx \binits  \undefined \def \binits#1{#1}\fi
\ifx \bauthor  \undefined \def \bauthor#1{#1}\fi
\ifx \batitle  \undefined \def \batitle#1{#1}\fi
\ifx \bjtitle  \undefined \def \bjtitle#1{#1}\fi
\ifx \bvolume  \undefined \def \bvolume#1{\textbf{#1}}\fi
\ifx \byear  \undefined \def \byear#1{#1}\fi
\ifx \bissue  \undefined \def \bissue#1{#1}\fi
\ifx \bfpage  \undefined \def \bfpage#1{#1}\fi
\ifx \blpage  \undefined \def \blpage #1{#1}\fi
\ifx \burl  \undefined \def \burl#1{\textsf{#1}}\fi
\ifx \doiurl  \undefined \def \doiurl#1{\url{https://doi.org/#1}}\fi
\ifx \betal  \undefined \def \betal{\textit{et al.}}\fi
\ifx \binstitute  \undefined \def \binstitute#1{#1}\fi
\ifx \binstitutionaled  \undefined \def \binstitutionaled#1{#1}\fi
\ifx \bctitle  \undefined \def \bctitle#1{#1}\fi
\ifx \beditor  \undefined \def \beditor#1{#1}\fi
\ifx \bpublisher  \undefined \def \bpublisher#1{#1}\fi
\ifx \bbtitle  \undefined \def \bbtitle#1{#1}\fi
\ifx \bedition  \undefined \def \bedition#1{#1}\fi
\ifx \bseriesno  \undefined \def \bseriesno#1{#1}\fi
\ifx \blocation  \undefined \def \blocation#1{#1}\fi
\ifx \bsertitle  \undefined \def \bsertitle#1{#1}\fi
\ifx \bsnm \undefined \def \bsnm#1{#1}\fi
\ifx \bsuffix \undefined \def \bsuffix#1{#1}\fi
\ifx \bparticle \undefined \def \bparticle#1{#1}\fi
\ifx \barticle \undefined \def \barticle#1{#1}\fi
\bibcommenthead
\ifx \bconfdate \undefined \def \bconfdate #1{#1}\fi
\ifx \botherref \undefined \def \botherref #1{#1}\fi
\ifx \url \undefined \def \url#1{\textsf{#1}}\fi
\ifx \bchapter \undefined \def \bchapter#1{#1}\fi
\ifx \bbook \undefined \def \bbook#1{#1}\fi
\ifx \bcomment \undefined \def \bcomment#1{#1}\fi
\ifx \oauthor \undefined \def \oauthor#1{#1}\fi
\ifx \citeauthoryear \undefined \def \citeauthoryear#1{#1}\fi
\ifx \endbibitem  \undefined \def \endbibitem {}\fi
\ifx \bconflocation  \undefined \def \bconflocation#1{#1}\fi
\ifx \arxivurl  \undefined \def \arxivurl#1{\textsf{#1}}\fi
\csname PreBibitemsHook\endcsname

\bibitem{cichocki2016tensor}
\begin{barticle}
\bauthor{\bsnm{Cichocki}, \binits{A.}},
\bauthor{\bsnm{Lee}, \binits{N.}},
\bauthor{\bsnm{Oseledets}, \binits{I.}},
\bauthor{\bsnm{Phan}, \binits{A.-H.}},
\bauthor{\bsnm{Zhao}, \binits{Q.}},
\bauthor{\bsnm{Mandic}, \binits{D.P.}}:
\batitle{Tensor networks for dimensionality reduction and large-scale
  optimization: {P}art 1 low-rank tensor decompositions}.
\bjtitle{Foundations and Trends{\textregistered} in Machine Learning}
\bvolume{9}(\bissue{4-5}),
\bfpage{249}--\blpage{429}
(\byear{2016})
\end{barticle}
\endbibitem

\bibitem{cichocki2017tensor}
\begin{barticle}
\bauthor{\bsnm{Cichocki}, \binits{A.}},
\bauthor{\bsnm{Phan}, \binits{A.-H.}},
\bauthor{\bsnm{Zhao}, \binits{Q.}},
\bauthor{\bsnm{Lee}, \binits{N.}},
\bauthor{\bsnm{Oseledets}, \binits{I.}},
\bauthor{\bsnm{Sugiyama}, \binits{M.}},
\bauthor{\bsnm{Mandic}, \binits{D.P.}}:
\batitle{Tensor networks for dimensionality reduction and large-scale
  optimization: {P}art 2 applications and future perspectives}.
\bjtitle{Foundations and Trends{\textregistered} in Machine Learning}
\bvolume{9}(\bissue{6}),
\bfpage{431}--\blpage{673}
(\byear{2017})
\end{barticle}
\endbibitem

\bibitem{cichocki2015tensor}
\begin{barticle}
\bauthor{\bsnm{Cichocki}, \binits{A.}},
\bauthor{\bsnm{Mandic}, \binits{D.}},
\bauthor{\bsnm{De~Lathauwer}, \binits{L.}},
\bauthor{\bsnm{Zhou}, \binits{G.}},
\bauthor{\bsnm{Zhao}, \binits{Q.}},
\bauthor{\bsnm{Caiafa}, \binits{C.}},
\bauthor{\bsnm{Phan}, \binits{H.A.}}:
\batitle{Tensor decompositions for signal processing applications: From two-way
  to multiway component analysis}.
\bjtitle{IEEE signal processing magazine}
\bvolume{32}(\bissue{2}),
\bfpage{145}--\blpage{163}
(\byear{2015})
\end{barticle}
\endbibitem

\bibitem{cichocki2009nonnegative}
\begin{botherref}
\oauthor{\bsnm{Cichocki}, \binits{A.}},
\oauthor{\bsnm{Zdunek}, \binits{R.}},
\oauthor{\bsnm{Phan}, \binits{A.H.}},
\oauthor{\bsnm{Amari}, \binits{S.-i.}}:
Nonnegative matrix and tensor factorizations: applications to exploratory
  multi-way data analysis and blind source separation
(2009)
\end{botherref}
\endbibitem

\bibitem{long2019low}
\begin{barticle}
\bauthor{\bsnm{Long}, \binits{Z.}},
\bauthor{\bsnm{Liu}, \binits{Y.}},
\bauthor{\bsnm{Chen}, \binits{L.}},
\bauthor{\bsnm{Zhu}, \binits{C.}}:
\batitle{Low rank tensor completion for multiway visual data}.
\bjtitle{Signal Processing}
\bvolume{155},
\bfpage{301}--\blpage{316}
(\byear{2019})
\end{barticle}
\endbibitem

\bibitem{song2019tensor}
\begin{barticle}
\bauthor{\bsnm{Song}, \binits{Q.}},
\bauthor{\bsnm{Ge}, \binits{H.}},
\bauthor{\bsnm{Caverlee}, \binits{J.}},
\bauthor{\bsnm{Hu}, \binits{X.}}:
\batitle{Tensor completion algorithms in big data analytics}.
\bjtitle{ACM Transactions on Knowledge Discovery from Data (TKDD)}
\bvolume{13}(\bissue{1}),
\bfpage{1}--\blpage{48}
(\byear{2019})
\end{barticle}
\endbibitem

\bibitem{asante2021matrix}
\begin{barticle}
\bauthor{\bsnm{Asante-Mensah}, \binits{M.G.}},
\bauthor{\bsnm{Ahmadi-Asl}, \binits{S.}},
\bauthor{\bsnm{Cichocki}, \binits{A.}}:
\batitle{Matrix and tensor completion using tensor ring decomposition with
  sparse representation}.
\bjtitle{Machine Learning: Science and Technology}
\bvolume{2}(\bissue{3}),
\bfpage{035008}
(\byear{2021})
\end{barticle}
\endbibitem

\bibitem{frolov2017tensor}
\begin{barticle}
\bauthor{\bsnm{Frolov}, \binits{E.}},
\bauthor{\bsnm{Oseledets}, \binits{I.}}:
\batitle{Tensor methods and recommender systems}.
\bjtitle{Wiley Interdisciplinary Reviews: Data Mining and Knowledge Discovery}
\bvolume{7}(\bissue{3}),
\bfpage{1201}
(\byear{2017})
\end{barticle}
\endbibitem

\bibitem{liu2012tensor}
\begin{barticle}
\bauthor{\bsnm{Liu}, \binits{J.}},
\bauthor{\bsnm{Musialski}, \binits{P.}},
\bauthor{\bsnm{Wonka}, \binits{P.}},
\bauthor{\bsnm{Ye}, \binits{J.}}:
\batitle{Tensor completion for estimating missing values in visual data}.
\bjtitle{IEEE transactions on pattern analysis and machine intelligence}
\bvolume{35}(\bissue{1}),
\bfpage{208}--\blpage{220}
(\byear{2012})
\end{barticle}
\endbibitem

\bibitem{tomasi2005parafac}
\begin{barticle}
\bauthor{\bsnm{Tomasi}, \binits{G.}},
\bauthor{\bsnm{Bro}, \binits{R.}}:
\batitle{Parafac and missing values}.
\bjtitle{Chemometrics and Intelligent Laboratory Systems}
\bvolume{75}(\bissue{2}),
\bfpage{163}--\blpage{180}
(\byear{2005})
\end{barticle}
\endbibitem

\bibitem{acar2009link}
\begin{bchapter}
\bauthor{\bsnm{Acar}, \binits{E.}},
\bauthor{\bsnm{Dunlavy}, \binits{D.M.}},
\bauthor{\bsnm{Kolda}, \binits{T.G.}}:
\bctitle{Link prediction on evolving data using matrix and tensor
  factorizations}.
In: \bbtitle{2009 IEEE International Conference on Data Mining Workshops},
pp. \bfpage{262}--\blpage{269}
(\byear{2009}).
\bcomment{IEEE}
\end{bchapter}
\endbibitem

\bibitem{hitchcock1928multiple}
\begin{barticle}
\bauthor{\bsnm{Hitchcock}, \binits{F.L.}}:
\batitle{Multiple invariants and generalized rank of a p-way matrix or tensor}.
\bjtitle{Journal of Mathematics and Physics}
\bvolume{7}(\bissue{1-4}),
\bfpage{39}--\blpage{79}
(\byear{1928})
\end{barticle}
\endbibitem

\bibitem{hitchcock1927expression}
\begin{barticle}
\bauthor{\bsnm{Hitchcock}, \binits{F.L.}}:
\batitle{The expression of a tensor or a polyadic as a sum of products}.
\bjtitle{Journal of Mathematics and Physics}
\bvolume{6}(\bissue{1-4}),
\bfpage{164}--\blpage{189}
(\byear{1927})
\end{barticle}
\endbibitem

\bibitem{tucker1963implications}
\begin{barticle}
\bauthor{\bsnm{Tucker}, \binits{L.R.}}:
\batitle{Implications of factor analysis of three-way matrices for measurement
  of change}.
\bjtitle{Problems in measuring change}
\bvolume{15},
\bfpage{122}--\blpage{137}
(\byear{1963})
\end{barticle}
\endbibitem

\bibitem{tucker1964extension}
\begin{botherref}
\oauthor{\bsnm{Tucker}, \binits{L.R.}}, et al.:
The extension of factor analysis to three-dimensional matrices.
Contributions to mathematical psychology
\textbf{110119}
(1964)
\end{botherref}
\endbibitem

\bibitem{tucker1966some}
\begin{barticle}
\bauthor{\bsnm{Tucker}, \binits{L.R.}}:
\batitle{Some mathematical notes on three-mode factor analysis}.
\bjtitle{Psychometrika}
\bvolume{31}(\bissue{3}),
\bfpage{279}--\blpage{311}
(\byear{1966})
\end{barticle}
\endbibitem

\bibitem{de2000multilinear}
\begin{barticle}
\bauthor{\bsnm{De~Lathauwer}, \binits{L.}},
\bauthor{\bsnm{De~Moor}, \binits{B.}},
\bauthor{\bsnm{Vandewalle}, \binits{J.}}:
\batitle{A multilinear singular value decomposition}.
\bjtitle{SIAM journal on Matrix Analysis and Applications}
\bvolume{21}(\bissue{4}),
\bfpage{1253}--\blpage{1278}
(\byear{2000})
\end{barticle}
\endbibitem

\bibitem{de2008decompositionsI}
\begin{barticle}
\bauthor{\bsnm{De~Lathauwer}, \binits{L.}}:
\batitle{Decompositions of a higher-order tensor in block terms—{P}art {II}:
  Definitions and uniqueness}.
\bjtitle{SIAM Journal on Matrix Analysis and Applications}
\bvolume{30}(\bissue{3}),
\bfpage{1033}--\blpage{1066}
(\byear{2008})
\end{barticle}
\endbibitem

\bibitem{de2008decompositionsII}
\begin{barticle}
\bauthor{\bsnm{De~Lathauwer}, \binits{L.}},
\bauthor{\bsnm{Nion}, \binits{D.}}:
\batitle{Decompositions of a higher-order tensor in block terms—{P}art {III}:
  Alternating least squares algorithms}.
\bjtitle{SIAM journal on Matrix Analysis and Applications}
\bvolume{30}(\bissue{3}),
\bfpage{1067}--\blpage{1083}
(\byear{2008})
\end{barticle}
\endbibitem

\bibitem{de2008decompositionsIII}
\begin{barticle}
\bauthor{\bsnm{De~Lathauwer}, \binits{L.}}:
\batitle{Decompositions of a higher-order tensor in block terms—{P}art {I}:
  Lemmas for partitioned matrices}.
\bjtitle{SIAM Journal on Matrix Analysis and Applications}
\bvolume{30}(\bissue{3}),
\bfpage{1022}--\blpage{1032}
(\byear{2008})
\end{barticle}
\endbibitem

\bibitem{oseledets2011tensor}
\begin{barticle}
\bauthor{\bsnm{Oseledets}, \binits{I.V.}}:
\batitle{Tensor-train decomposition}.
\bjtitle{SIAM Journal on Scientific Computing}
\bvolume{33}(\bissue{5}),
\bfpage{2295}--\blpage{2317}
(\byear{2011})
\end{barticle}
\endbibitem

\bibitem{zhao2016tensor}
\begin{botherref}
\oauthor{\bsnm{Zhao}, \binits{Q.}},
\oauthor{\bsnm{Zhou}, \binits{G.}},
\oauthor{\bsnm{Xie}, \binits{S.}},
\oauthor{\bsnm{Zhang}, \binits{L.}},
\oauthor{\bsnm{Cichocki}, \binits{A.}}:
Tensor ring decomposition.
arXiv preprint arXiv:1606.05535
(2016)
\end{botherref}
\endbibitem

\bibitem{espig2012note}
\begin{barticle}
\bauthor{\bsnm{Espig}, \binits{M.}},
\bauthor{\bsnm{Naraparaju}, \binits{K.K.}},
\bauthor{\bsnm{Schneider}, \binits{J.}}:
\batitle{A note on tensor chain approximation}.
\bjtitle{Computing and Visualization in Science}
\bvolume{15}(\bissue{6}),
\bfpage{331}--\blpage{344}
(\byear{2012})
\end{barticle}
\endbibitem

\bibitem{kilmer2011factorization}
\begin{barticle}
\bauthor{\bsnm{Kilmer}, \binits{M.E.}},
\bauthor{\bsnm{Martin}, \binits{C.D.}}:
\batitle{Factorization strategies for third-order tensors}.
\bjtitle{Linear Algebra and its Applications}
\bvolume{435}(\bissue{3}),
\bfpage{641}--\blpage{658}
(\byear{2011})
\end{barticle}
\endbibitem

\bibitem{kilmer2013third}
\begin{barticle}
\bauthor{\bsnm{Kilmer}, \binits{M.E.}},
\bauthor{\bsnm{Braman}, \binits{K.}},
\bauthor{\bsnm{Hao}, \binits{N.}},
\bauthor{\bsnm{Hoover}, \binits{R.C.}}:
\batitle{Third-order tensors as operators on matrices: A theoretical and
  computational framework with applications in imaging}.
\bjtitle{SIAM Journal on Matrix Analysis and Applications}
\bvolume{34}(\bissue{1}),
\bfpage{148}--\blpage{172}
(\byear{2013})
\end{barticle}
\endbibitem

\bibitem{braman2010third}
\begin{barticle}
\bauthor{\bsnm{Braman}, \binits{K.}}:
\batitle{Third-order tensors as linear operators on a space of matrices}.
\bjtitle{Linear Algebra and its Applications}
\bvolume{433}(\bissue{7}),
\bfpage{1241}--\blpage{1253}
(\byear{2010})
\end{barticle}
\endbibitem

\bibitem{ahmadi2021randomized}
\begin{barticle}
\bauthor{\bsnm{Ahmadi-Asl}, \binits{S.}},
\bauthor{\bsnm{Abukhovich}, \binits{S.}},
\bauthor{\bsnm{Asante-Mensah}, \binits{M.G.}},
\bauthor{\bsnm{Cichocki}, \binits{A.}},
\bauthor{\bsnm{Phan}, \binits{A.H.}},
\bauthor{\bsnm{Tanaka}, \binits{T.}},
\bauthor{\bsnm{Oseledets}, \binits{I.}}:
\batitle{Randomized algorithms for computation of {T}ucker decomposition and
  higher order svd ({HOSVD})}.
\bjtitle{IEEE Access}
\bvolume{9},
\bfpage{28684}--\blpage{28706}
(\byear{2021})
\end{barticle}
\endbibitem

\bibitem{ahmadiTRrandomized2020}
\begin{barticle}
\bauthor{\bsnm{Ahmadi-Asl}, \binits{S.}},
\bauthor{\bsnm{Cichocki}, \binits{A.}},
\bauthor{\bsnm{Phan}, \binits{A.H.}},
\bauthor{\bsnm{Asante-Mensah}, \binits{M.G.}},
\bauthor{\bsnm{Ghazani}, \binits{M.M.}},
\bauthor{\bsnm{Tanaka}, \binits{T.}},
\bauthor{\bsnm{Oseledets}, \binits{I.}}:
\batitle{Randomized algorithms for fast computation of low rank tensor ring
  model}.
\bjtitle{Machine Learning: Science and Technology}
\bvolume{2}(\bissue{1}),
\bfpage{011001}
(\byear{2020})
\end{barticle}
\endbibitem

\bibitem{xu2015cur}
\begin{bchapter}
\bauthor{\bsnm{Xu}, \binits{M.}},
\bauthor{\bsnm{Jin}, \binits{R.}},
\bauthor{\bsnm{Zhou}, \binits{Z.-H.}}:
\bctitle{{CUR} algorithm for partially observed matrices}.
In: \bbtitle{International Conference on Machine Learning},
pp. \bfpage{1412}--\blpage{1421}
(\byear{2015})
\end{bchapter}
\endbibitem

\bibitem{wang2017missing}
\begin{barticle}
\bauthor{\bsnm{Wang}, \binits{L.}},
\bauthor{\bsnm{Xie}, \binits{K.}},
\bauthor{\bsnm{Semong}, \binits{T.}},
\bauthor{\bsnm{Zhou}, \binits{H.}}:
\batitle{Missing data recovery based on tensor-{CUR} decomposition}.
\bjtitle{IEEE Access}
\bvolume{6},
\bfpage{532}--\blpage{544}
(\byear{2017})
\end{barticle}
\endbibitem

\bibitem{goreinov1997theory}
\begin{barticle}
\bauthor{\bsnm{Goreinov}, \binits{S.A.}},
\bauthor{\bsnm{Tyrtyshnikov}, \binits{E.E.}},
\bauthor{\bsnm{Zamarashkin}, \binits{N.L.}}:
\batitle{A theory of pseudoskeleton approximations}.
\bjtitle{Linear algebra and its applications}
\bvolume{261}(\bissue{1-3}),
\bfpage{1}--\blpage{21}
(\byear{1997})
\end{barticle}
\endbibitem

\bibitem{goreinov2010find}
\begin{botherref}
\oauthor{\bsnm{Goreinov}, \binits{S.A.}},
\oauthor{\bsnm{Oseledets}, \binits{I.V.}},
\oauthor{\bsnm{Savostyanov}, \binits{D.V.}},
\oauthor{\bsnm{Tyrtyshnikov}, \binits{E.E.}},
\oauthor{\bsnm{Zamarashkin}, \binits{N.L.}}:
How to find a good submatrix,
247--256
(2010)
\end{botherref}
\endbibitem

\bibitem{goreinov2001maximal}
\begin{barticle}
\bauthor{\bsnm{Goreinov}, \binits{S.A.}},
\bauthor{\bsnm{Tyrtyshnikov}, \binits{E.E.}}:
\batitle{The maximal-volume concept in approximation by low-rank matrices}.
\bjtitle{Contemporary Mathematics}
\bvolume{280},
\bfpage{47}--\blpage{52}
(\byear{2001})
\end{barticle}
\endbibitem

\bibitem{savostyanov2006polilinear}
\begin{botherref}
\oauthor{\bsnm{Savostyanov}, \binits{D.}}:
Polilinear approximation of matrices and integral equations.
PhD thesis,
PhD thesis, INM RAS, Moscow, 2006.(in Russian)
(2006)
\end{botherref}
\endbibitem

\bibitem{tyrtyshnikov2000incomplete}
\begin{barticle}
\bauthor{\bsnm{Tyrtyshnikov}, \binits{E.}}:
\batitle{Incomplete cross approximation in the mosaic-skeleton method}.
\bjtitle{Computing}
\bvolume{64}(\bissue{4}),
\bfpage{367}--\blpage{380}
(\byear{2000})
\end{barticle}
\endbibitem

\bibitem{chaturantabut2009discrete}
\begin{bchapter}
\bauthor{\bsnm{Chaturantabut}, \binits{S.}},
\bauthor{\bsnm{Sorensen}, \binits{D.C.}}:
\bctitle{Discrete empirical interpolation for nonlinear model reduction}.
In: \bbtitle{Proceedings of the 48h IEEE Conference on Decision and Control
  (CDC) Held Jointly with 2009 28th Chinese Control Conference},
pp. \bfpage{4316}--\blpage{4321}
(\byear{2009}).
\bcomment{IEEE}
\end{bchapter}
\endbibitem

\bibitem{sorensen2016deim}
\begin{barticle}
\bauthor{\bsnm{Sorensen}, \binits{D.C.}},
\bauthor{\bsnm{Embree}, \binits{M.}}:
\batitle{A {DEIM} induced {CUR} factorization}.
\bjtitle{SIAM Journal on Scientific Computing}
\bvolume{38}(\bissue{3}),
\bfpage{1454}--\blpage{1482}
(\byear{2016})
\end{barticle}
\endbibitem

\bibitem{frieze2004fast}
\begin{barticle}
\bauthor{\bsnm{Frieze}, \binits{A.}},
\bauthor{\bsnm{Kannan}, \binits{R.}},
\bauthor{\bsnm{Vempala}, \binits{S.}}:
\batitle{Fast {M}onte-{C}arlo algorithms for finding low-rank approximations}.
\bjtitle{Journal of the ACM (JACM)}
\bvolume{51}(\bissue{6}),
\bfpage{1025}--\blpage{1041}
(\byear{2004})
\end{barticle}
\endbibitem

\bibitem{boutsidis2009improved}
\begin{bchapter}
\bauthor{\bsnm{Boutsidis}, \binits{C.}},
\bauthor{\bsnm{Mahoney}, \binits{M.W.}},
\bauthor{\bsnm{Drineas}, \binits{P.}}:
\bctitle{An improved approximation algorithm for the column subset selection
  problem}.
In: \bbtitle{Proceedings of the Twentieth Annual ACM-SIAM Symposium on Discrete
  Algorithms},
pp. \bfpage{968}--\blpage{977}
(\byear{2009}).
\bcomment{SIAM}
\end{bchapter}
\endbibitem

\bibitem{boutsidis2014near}
\begin{barticle}
\bauthor{\bsnm{Boutsidis}, \binits{C.}},
\bauthor{\bsnm{Drineas}, \binits{P.}},
\bauthor{\bsnm{Magdon-Ismail}, \binits{M.}}:
\batitle{Near-optimal column-based matrix reconstruction}.
\bjtitle{SIAM Journal on Computing}
\bvolume{43}(\bissue{2}),
\bfpage{687}--\blpage{717}
(\byear{2014})
\end{barticle}
\endbibitem

\bibitem{deshpande2006adaptive}
\begin{botherref}
\oauthor{\bsnm{Deshpande}, \binits{A.}},
\oauthor{\bsnm{Vempala}, \binits{S.}}:
Adaptive sampling and fast low-rank matrix approximation,
292--303
(2006)
\end{botherref}
\endbibitem

\bibitem{deshpande2006matrix}
\begin{barticle}
\bauthor{\bsnm{Deshpande}, \binits{A.}},
\bauthor{\bsnm{Rademacher}, \binits{L.}},
\bauthor{\bsnm{Vempala}, \binits{S.}},
\bauthor{\bsnm{Wang}, \binits{G.}}:
\batitle{Matrix approximation and projective clustering via volume sampling}.
\bjtitle{Theory of Computing}
\bvolume{2}(\bissue{1}),
\bfpage{225}--\blpage{247}
(\byear{2006})
\end{barticle}
\endbibitem

\bibitem{deshpande2010efficient}
\begin{bchapter}
\bauthor{\bsnm{Deshpande}, \binits{A.}},
\bauthor{\bsnm{Rademacher}, \binits{L.}}:
\bctitle{Efficient volume sampling for row/column subset selection}.
In: \bbtitle{2010 Ieee 51st Annual Symposium on Foundations of Computer
  Science},
pp. \bfpage{329}--\blpage{338}
(\byear{2010}).
\bcomment{IEEE}
\end{bchapter}
\endbibitem

\bibitem{guruswami2012optimal}
\begin{bchapter}
\bauthor{\bsnm{Guruswami}, \binits{V.}},
\bauthor{\bsnm{Sinop}, \binits{A.K.}}:
\bctitle{Optimal column-based low-rank matrix reconstruction}.
In: \bbtitle{Proceedings of the Twenty-third Annual ACM-SIAM Symposium on
  Discrete Algorithms},
pp. \bfpage{1207}--\blpage{1214}
(\byear{2012}).
\bcomment{SIAM}
\end{bchapter}
\endbibitem

\bibitem{mai2020vgg}
\begin{bchapter}
\bauthor{\bsnm{Mai}, \binits{A.}},
\bauthor{\bsnm{Tran}, \binits{L.}},
\bauthor{\bsnm{Tran}, \binits{L.}},
\bauthor{\bsnm{Trinh}, \binits{N.}}:
\bctitle{{VGG} deep neural network compression via {SVD} and {CUR}
  decomposition techniques}.
In: \bbtitle{2020 7th NAFOSTED Conference on Information and Computer Science
  (NICS)},
pp. \bfpage{118}--\blpage{123}
(\byear{2020}).
\bcomment{IEEE}
\end{bchapter}
\endbibitem

\bibitem{hendryx2018finding}
\begin{barticle}
\bauthor{\bsnm{Hendryx}, \binits{E.P.}},
\bauthor{\bsnm{Rivi{\`e}re}, \binits{B.M.}},
\bauthor{\bsnm{Sorensen}, \binits{D.C.}},
\bauthor{\bsnm{Rusin}, \binits{C.G.}}:
\batitle{Finding representative electrocardiogram beat morphologies with
  {CUR}}.
\bjtitle{Journal of biomedical informatics}
\bvolume{77},
\bfpage{97}--\blpage{110}
(\byear{2018})
\end{barticle}
\endbibitem

\bibitem{cai2020rapid}
\begin{botherref}
\oauthor{\bsnm{Cai}, \binits{H.}},
\oauthor{\bsnm{Hamm}, \binits{K.}},
\oauthor{\bsnm{Huang}, \binits{L.}},
\oauthor{\bsnm{Li}, \binits{J.}},
\oauthor{\bsnm{Wang}, \binits{T.}}:
Rapid robust principal component analysis: {CUR} accelerated inexact low rank
  estimation.
IEEE Signal Processing Letters
(2020)
\end{botherref}
\endbibitem

\bibitem{drineas2006fastI}
\begin{barticle}
\bauthor{\bsnm{Drineas}, \binits{P.}},
\bauthor{\bsnm{Kannan}, \binits{R.}},
\bauthor{\bsnm{Mahoney}, \binits{M.W.}}:
\batitle{Fast {M}onte {C}arlo algorithms for matrices {I}: Approximating matrix
  multiplication}.
\bjtitle{SIAM Journal on Computing}
\bvolume{36}(\bissue{1}),
\bfpage{132}--\blpage{157}
(\byear{2006})
\end{barticle}
\endbibitem

\bibitem{drineas2006fastII}
\begin{barticle}
\bauthor{\bsnm{Drineas}, \binits{P.}},
\bauthor{\bsnm{Kannan}, \binits{R.}},
\bauthor{\bsnm{Mahoney}, \binits{M.W.}}:
\batitle{Fast {M}onte {C}arlo algorithms for matrices ii: Computing a low-rank
  approximation to a matrix}.
\bjtitle{SIAM Journal on computing}
\bvolume{36}(\bissue{1}),
\bfpage{158}--\blpage{183}
(\byear{2006})
\end{barticle}
\endbibitem

\bibitem{drineas2006fastIII}
\begin{barticle}
\bauthor{\bsnm{Drineas}, \binits{P.}},
\bauthor{\bsnm{Kannan}, \binits{R.}},
\bauthor{\bsnm{Mahoney}, \binits{M.W.}}:
\batitle{Fast {M}onte {C}arlo algorithms for matrices {III}: Computing a
  compressed approximate matrix decomposition}.
\bjtitle{SIAM Journal on Computing}
\bvolume{36}(\bissue{1}),
\bfpage{184}--\blpage{206}
(\byear{2006})
\end{barticle}
\endbibitem

\bibitem{mahoney2009cur}
\begin{barticle}
\bauthor{\bsnm{Mahoney}, \binits{M.W.}},
\bauthor{\bsnm{Drineas}, \binits{P.}}:
\batitle{{CUR} matrix decompositions for improved data analysis}.
\bjtitle{Proceedings of the National Academy of Sciences}
\bvolume{106}(\bissue{3}),
\bfpage{697}--\blpage{702}
(\byear{2009})
\end{barticle}
\endbibitem

\bibitem{li2018joint}
\begin{barticle}
\bauthor{\bsnm{Li}, \binits{C.}},
\bauthor{\bsnm{Wang}, \binits{X.}},
\bauthor{\bsnm{Dong}, \binits{W.}},
\bauthor{\bsnm{Yan}, \binits{J.}},
\bauthor{\bsnm{Liu}, \binits{Q.}},
\bauthor{\bsnm{Zha}, \binits{H.}}:
\batitle{Joint active learning with feature selection via {CUR} matrix
  decomposition}.
\bjtitle{IEEE transactions on pattern analysis and machine intelligence}
\bvolume{41}(\bissue{6}),
\bfpage{1382}--\blpage{1396}
(\byear{2018})
\end{barticle}
\endbibitem

\bibitem{aldroubi2019cur}
\begin{barticle}
\bauthor{\bsnm{Aldroubi}, \binits{A.}},
\bauthor{\bsnm{Hamm}, \binits{K.}},
\bauthor{\bsnm{Koku}, \binits{A.B.}},
\bauthor{\bsnm{Sekmen}, \binits{A.}}:
\batitle{{CUR} decompositions, similarity matrices, and subspace clustering}.
\bjtitle{Frontiers in Applied Mathematics and Statistics}
\bvolume{4},
\bfpage{65}
(\byear{2019})
\end{barticle}
\endbibitem

\bibitem{ahmadi2021cross}
\begin{barticle}
\bauthor{\bsnm{Ahmadi-Asl}, \binits{S.}},
\bauthor{\bsnm{Caiafa}, \binits{C.F.}},
\bauthor{\bsnm{Cichocki}, \binits{A.}},
\bauthor{\bsnm{Phan}, \binits{A.H.}},
\bauthor{\bsnm{Tanaka}, \binits{T.}},
\bauthor{\bsnm{Oseledets}, \binits{I.}},
\bauthor{\bsnm{Wang}, \binits{J.}}:
\batitle{Cross tensor approximation methods for compression and dimensionality
  reduction}.
\bjtitle{IEEE Access}
\bvolume{9},
\bfpage{150809}--\blpage{150838}
(\byear{2021})
\end{barticle}
\endbibitem

\bibitem{drinea2001randomized}
\begin{bchapter}
\bauthor{\bsnm{Drinea}, \binits{E.}},
\bauthor{\bsnm{Drineas}, \binits{P.}},
\bauthor{\bsnm{Huggins}, \binits{P.}}:
\bctitle{A randomized singular value decomposition algorithm for image
  processing applications}.
In: \bbtitle{Proceedings of the 8th Panhellenic Conference on Informatics},
pp. \bfpage{278}--\blpage{288}
(\byear{2001}).
\bcomment{Citeseer}
\end{bchapter}
\endbibitem

\bibitem{oseledets2008tucker}
\begin{barticle}
\bauthor{\bsnm{Oseledets}, \binits{I.V.}},
\bauthor{\bsnm{Savostianov}, \binits{D.}},
\bauthor{\bsnm{Tyrtyshnikov}, \binits{E.E.}}:
\batitle{{T}ucker dimensionality reduction of three-dimensional arrays in
  linear time}.
\bjtitle{SIAM Journal on Matrix Analysis and Applications}
\bvolume{30}(\bissue{3}),
\bfpage{939}--\blpage{956}
(\byear{2008})
\end{barticle}
\endbibitem

\bibitem{caiafa2010generalizing}
\begin{barticle}
\bauthor{\bsnm{Caiafa}, \binits{C.F.}},
\bauthor{\bsnm{Cichocki}, \binits{A.}}:
\batitle{Generalizing the column--row matrix decomposition to multi-way
  arrays}.
\bjtitle{Linear Algebra and its Applications}
\bvolume{433}(\bissue{3}),
\bfpage{557}--\blpage{573}
(\byear{2010})
\end{barticle}
\endbibitem

\bibitem{mahoney2008tensor}
\begin{barticle}
\bauthor{\bsnm{Mahoney}, \binits{M.W.}},
\bauthor{\bsnm{Maggioni}, \binits{M.}},
\bauthor{\bsnm{Drineas}, \binits{P.}}:
\batitle{Tensor-{CUR} decompositions for tensor-based data}.
\bjtitle{SIAM Journal on Matrix Analysis and Applications}
\bvolume{30}(\bissue{3}),
\bfpage{957}--\blpage{987}
(\byear{2008})
\end{barticle}
\endbibitem

\bibitem{tarzanagh2018fast}
\begin{barticle}
\bauthor{\bsnm{Tarzanagh}, \binits{D.A.}},
\bauthor{\bsnm{Michailidis}, \binits{G.}}:
\batitle{Fast randomized algorithms for t-product based tensor operations and
  decompositions with applications to imaging data}.
\bjtitle{SIAM Journal on Imaging Sciences}
\bvolume{11}(\bissue{4}),
\bfpage{2629}--\blpage{2664}
(\byear{2018})
\end{barticle}
\endbibitem

\bibitem{fazel2002matrix}
\begin{botherref}
\oauthor{\bsnm{Fazel}, \binits{M.}}:
Matrix rank minimization with applications.
PhD thesis,
PhD thesis, Stanford University
(2002)
\end{botherref}
\endbibitem

\bibitem{fang2018sequentially}
\begin{barticle}
\bauthor{\bsnm{Fang}, \binits{Z.}},
\bauthor{\bsnm{Yang}, \binits{X.}},
\bauthor{\bsnm{Han}, \binits{L.}},
\bauthor{\bsnm{Liu}, \binits{X.}}:
\batitle{A sequentially truncated higher order singular value
  decomposition-based algorithm for tensor completion}.
\bjtitle{IEEE transactions on cybernetics}
\bvolume{49}(\bissue{5}),
\bfpage{1956}--\blpage{1967}
(\byear{2018})
\end{barticle}
\endbibitem

\bibitem{acar2011scalable}
\begin{barticle}
\bauthor{\bsnm{Acar}, \binits{E.}},
\bauthor{\bsnm{Dunlavy}, \binits{D.M.}},
\bauthor{\bsnm{Kolda}, \binits{T.G.}}:
\batitle{A scalable optimization approach for fitting canonical tensor
  decompositions}.
\bjtitle{Journal of Chemometrics}
\bvolume{25}(\bissue{2}),
\bfpage{67}--\blpage{86}
(\byear{2011})
\end{barticle}
\endbibitem

\bibitem{wang2017efficient}
\begin{bchapter}
\bauthor{\bsnm{Wang}, \binits{W.}},
\bauthor{\bsnm{Aggarwal}, \binits{V.}},
\bauthor{\bsnm{Aeron}, \binits{S.}}:
\bctitle{Efficient low rank tensor ring completion}.
In: \bbtitle{Proceedings of the IEEE International Conference on Computer
  Vision},
pp. \bfpage{5697}--\blpage{5705}
(\byear{2017})
\end{bchapter}
\endbibitem

\bibitem{candes2009exact}
\begin{barticle}
\bauthor{\bsnm{Cand{\`e}s}, \binits{E.J.}},
\bauthor{\bsnm{Recht}, \binits{B.}}:
\batitle{Exact matrix completion via convex optimization}.
\bjtitle{Foundations of Computational mathematics}
\bvolume{9}(\bissue{6}),
\bfpage{717}--\blpage{772}
(\byear{2009})
\end{barticle}
\endbibitem

\bibitem{liu2012tensor2}
\begin{botherref}
\oauthor{\bsnm{Liu}, \binits{J.}},
\oauthor{\bsnm{Musialski}, \binits{P.}},
\oauthor{\bsnm{Wonka}, \binits{P.}},
\oauthor{\bsnm{Ye}, \binits{J.}}:
Tensor completion for estimating missing values in visual data.
ICCV
(2009)
\end{botherref}
\endbibitem

\bibitem{zhao2015bayesian}
\begin{barticle}
\bauthor{\bsnm{Zhao}, \binits{Q.}},
\bauthor{\bsnm{Zhang}, \binits{L.}},
\bauthor{\bsnm{Cichocki}, \binits{A.}}:
\batitle{Bayesian cp factorization of incomplete tensors with automatic rank
  determination}.
\bjtitle{IEEE transactions on pattern analysis and machine intelligence}
\bvolume{37}(\bissue{9}),
\bfpage{1751}--\blpage{1763}
(\byear{2015})
\end{barticle}
\endbibitem

\bibitem{yokota2016smooth}
\begin{barticle}
\bauthor{\bsnm{Yokota}, \binits{T.}},
\bauthor{\bsnm{Zhao}, \binits{Q.}},
\bauthor{\bsnm{Cichocki}, \binits{A.}}:
\batitle{Smooth parafac decomposition for tensor completion}.
\bjtitle{IEEE Transactions on Signal Processing}
\bvolume{64}(\bissue{20}),
\bfpage{5423}--\blpage{5436}
(\byear{2016})
\end{barticle}
\endbibitem

\bibitem{zhang2016exact}
\begin{barticle}
\bauthor{\bsnm{Zhang}, \binits{Z.}},
\bauthor{\bsnm{Aeron}, \binits{S.}}:
\batitle{Exact tensor completion using t-svd}.
\bjtitle{IEEE Transactions on Signal Processing}
\bvolume{65}(\bissue{6}),
\bfpage{1511}--\blpage{1526}
(\byear{2016})
\end{barticle}
\endbibitem

\bibitem{bengua2017efficient}
\begin{barticle}
\bauthor{\bsnm{Bengua}, \binits{J.A.}},
\bauthor{\bsnm{Phien}, \binits{H.N.}},
\bauthor{\bsnm{Tuan}, \binits{H.D.}},
\bauthor{\bsnm{Do}, \binits{M.N.}}:
\batitle{Efficient tensor completion for color image and video recovery:
  Low-rank tensor train}.
\bjtitle{IEEE Transactions on Image Processing}
\bvolume{26}(\bissue{5}),
\bfpage{2466}--\blpage{2479}
(\byear{2017})
\end{barticle}
\endbibitem

\bibitem{yuan2018higher}
\begin{bchapter}
\bauthor{\bsnm{Yuan}, \binits{L.}},
\bauthor{\bsnm{Cao}, \binits{J.}},
\bauthor{\bsnm{Zhao}, \binits{X.}},
\bauthor{\bsnm{Wu}, \binits{Q.}},
\bauthor{\bsnm{Zhao}, \binits{Q.}}:
\bctitle{Higher-dimension tensor completion via low-rank tensor ring
  decomposition}.
In: \bbtitle{2018 Asia-Pacific Signal and Information Processing Association
  Annual Summit and Conference (APSIPA ASC)},
pp. \bfpage{1071}--\blpage{1076}
(\byear{2018}).
\bcomment{IEEE}
\end{bchapter}
\endbibitem

\bibitem{yokota2018missing}
\begin{bchapter}
\bauthor{\bsnm{Yokota}, \binits{T.}},
\bauthor{\bsnm{Erem}, \binits{B.}},
\bauthor{\bsnm{Guler}, \binits{S.}},
\bauthor{\bsnm{Warfield}, \binits{S.K.}},
\bauthor{\bsnm{Hontani}, \binits{H.}}:
\bctitle{Missing slice recovery for tensors using a low-rank model in embedded
  space}.
In: \bbtitle{Proceedings of the IEEE Conference on Computer Vision and Pattern
  Recognition},
pp. \bfpage{8251}--\blpage{8259}
(\byear{2018})
\end{bchapter}
\endbibitem

\bibitem{yamamoto2022fast}
\begin{bchapter}
\bauthor{\bsnm{Yamamoto}, \binits{R.}},
\bauthor{\bsnm{Hontani}, \binits{H.}},
\bauthor{\bsnm{Imakura}, \binits{A.}},
\bauthor{\bsnm{Yokota}, \binits{T.}}:
\bctitle{Fast algorithm for low-rank tensor completion in delay-embedded
  space}.
In: \bbtitle{Proceedings of the IEEE/CVF Conference on Computer Vision and
  Pattern Recognition},
pp. \bfpage{2058}--\blpage{2066}
(\byear{2022})
\end{bchapter}
\endbibitem

\bibitem{boyd2011distributed}
\begin{barticle}
\bauthor{\bsnm{Boyd}, \binits{S.}},
\bauthor{\bsnm{Parikh}, \binits{N.}},
\bauthor{\bsnm{Chu}, \binits{E.}},
\bauthor{\bsnm{Peleato}, \binits{B.}},
\bauthor{\bsnm{Eckstein}, \binits{J.}}, \betal:
\batitle{Distributed optimization and statistical learning via the alternating
  direction method of multipliers}.
\bjtitle{Foundations and Trends{\textregistered} in Machine learning}
\bvolume{3}(\bissue{1}),
\bfpage{1}--\blpage{122}
(\byear{2011})
\end{barticle}
\endbibitem

\bibitem{halko2011finding}
\begin{barticle}
\bauthor{\bsnm{Halko}, \binits{N.}},
\bauthor{\bsnm{Martinsson}, \binits{P.-G.}},
\bauthor{\bsnm{Tropp}, \binits{J.A.}}:
\batitle{Finding structure with randomness: Probabilistic algorithms for
  constructing approximate matrix decompositions}.
\bjtitle{SIAM review}
\bvolume{53}(\bissue{2}),
\bfpage{217}--\blpage{288}
(\byear{2011})
\end{barticle}
\endbibitem

\bibitem{ulyanov2018deep}
\begin{bchapter}
\bauthor{\bsnm{Ulyanov}, \binits{D.}},
\bauthor{\bsnm{Vedaldi}, \binits{A.}},
\bauthor{\bsnm{Lempitsky}, \binits{V.}}:
\bctitle{Deep image prior}.
In: \bbtitle{Proceedings of the IEEE Conference on Computer Vision and Pattern
  Recognition},
pp. \bfpage{9446}--\blpage{9454}
(\byear{2018})
\end{bchapter}
\endbibitem

\bibitem{rombach2022high}
\begin{bchapter}
\bauthor{\bsnm{Rombach}, \binits{R.}},
\bauthor{\bsnm{Blattmann}, \binits{A.}},
\bauthor{\bsnm{Lorenz}, \binits{D.}},
\bauthor{\bsnm{Esser}, \binits{P.}},
\bauthor{\bsnm{Ommer}, \binits{B.}}:
\bctitle{High-resolution image synthesis with latent diffusion models}.
In: \bbtitle{Proceedings of the IEEE/CVF Conference on Computer Vision and
  Pattern Recognition},
pp. \bfpage{10684}--\blpage{10695}
(\byear{2022})
\end{bchapter}
\endbibitem

\bibitem{qin2021image}
\begin{barticle}
\bauthor{\bsnm{Qin}, \binits{Z.}},
\bauthor{\bsnm{Zeng}, \binits{Q.}},
\bauthor{\bsnm{Zong}, \binits{Y.}},
\bauthor{\bsnm{Xu}, \binits{F.}}:
\batitle{Image inpainting based on deep learning: A review}.
\bjtitle{Displays}
\bvolume{69},
\bfpage{102028}
(\byear{2021})
\end{barticle}
\endbibitem

\bibitem{oseledets2010tt}
\begin{barticle}
\bauthor{\bsnm{Oseledets}, \binits{I.}},
\bauthor{\bsnm{Tyrtyshnikov}, \binits{E.}}:
\batitle{{TT}-cross approximation for multidimensional arrays}.
\bjtitle{Linear Algebra and its Applications}
\bvolume{432}(\bissue{1}),
\bfpage{70}--\blpage{88}
(\byear{2010})
\end{barticle}
\endbibitem

\bibitem{liu2022deep}
\begin{botherref}
\oauthor{\bsnm{Liu}, \binits{R.W.}},
\oauthor{\bsnm{Guo}, \binits{Y.}},
\oauthor{\bsnm{Lu}, \binits{Y.}},
\oauthor{\bsnm{Chui}, \binits{K.T.}},
\oauthor{\bsnm{Gupta}, \binits{B.B.}}:
Deep network-enabled haze visibility enhancement for visual iot-driven
  intelligent transportation systems.
IEEE Transactions on Industrial Informatics
(2022)
\end{botherref}
\endbibitem

\bibitem{savitzky1964smoothing}
\begin{barticle}
\bauthor{\bsnm{Savitzky}, \binits{A.}},
\bauthor{\bsnm{Golay}, \binits{M.J.}}:
\batitle{Smoothing and differentiation of data by simplified least squares
  procedures.}
\bjtitle{Analytical chemistry}
\bvolume{36}(\bissue{8}),
\bfpage{1627}--\blpage{1639}
(\byear{1964})
\end{barticle}
\endbibitem

\bibitem{cleveland1981lowess}
\begin{barticle}
\bauthor{\bsnm{Cleveland}, \binits{W.S.}}:
\batitle{Lowess: A program for smoothing scatterplots by robust locally
  weighted regression}.
\bjtitle{American Statistician}
\bvolume{35}(\bissue{1}),
\bfpage{54}
(\byear{1981})
\end{barticle}
\endbibitem

\bibitem{cleveland1979robust}
\begin{barticle}
\bauthor{\bsnm{Cleveland}, \binits{W.S.}}:
\batitle{Robust locally weighted regression and smoothing scatterplots}.
\bjtitle{Journal of the American statistical association}
\bvolume{74}(\bissue{368}),
\bfpage{829}--\blpage{836}
(\byear{1979})
\end{barticle}
\endbibitem

\bibitem{garimella2017simple}
\begin{botherref}
\oauthor{\bsnm{Garimella}, \binits{R.V.}}:
A simple introduction to moving least squares and local regression estimation.
Technical report,
Los Alamos National Lab.(LANL), Los Alamos, NM (United States)
(2017)
\end{botherref}
\endbibitem

\bibitem{CHHH07}
\begin{bchapter}
\bauthor{\bsnm{Cai}, \binits{D.}},
\bauthor{\bsnm{He}, \binits{X.}},
\bauthor{\bsnm{Hu}, \binits{Y.}},
\bauthor{\bsnm{Han}, \binits{J.}},
\bauthor{\bsnm{Huang}, \binits{T.}}:
\bctitle{Learning a spatially smooth subspace for face recognition}.
In: \bbtitle{Proc. IEEE Conf. Computer Vision and Pattern Recognition Machine
  Learning (CVPR'07)}
(\byear{2007})
\end{bchapter}
\endbibitem

\bibitem{CHH07b}
\begin{bchapter}
\bauthor{\bsnm{Cai}, \binits{D.}},
\bauthor{\bsnm{He}, \binits{X.}},
\bauthor{\bsnm{Han}, \binits{J.}}:
\bctitle{Spectral regression for efficient regularized subspace learning}.
In: \bbtitle{Proc. Int. Conf. Computer Vision (ICCV'07)}
(\byear{2007})
\end{bchapter}
\endbibitem

\bibitem{CHHZ06}
\begin{barticle}
\bauthor{\bsnm{Cai}, \binits{D.}},
\bauthor{\bsnm{He}, \binits{X.}},
\bauthor{\bsnm{Han}, \binits{J.}},
\bauthor{\bsnm{Zhang}, \binits{H.-J.}}:
\batitle{Orthogonal laplacianfaces for face recognition}.
\bjtitle{IEEE Transactions on Image Processing}
\bvolume{15}(\bissue{11}),
\bfpage{3608}--\blpage{3614}
(\byear{2006})
\end{barticle}
\endbibitem

\bibitem{HYHNZ05}
\begin{barticle}
\bauthor{\bsnm{He}, \binits{X.}},
\bauthor{\bsnm{Yan}, \binits{S.}},
\bauthor{\bsnm{Hu}, \binits{Y.}},
\bauthor{\bsnm{Niyogi}, \binits{P.}},
\bauthor{\bsnm{Zhang}, \binits{H.-J.}}:
\batitle{Face recognition using laplacianfaces}.
\bjtitle{IEEE Trans. Pattern Anal. Mach. Intelligence}
\bvolume{27}(\bissue{3}),
\bfpage{328}--\blpage{340}
(\byear{2005})
\end{barticle}
\endbibitem

\bibitem{drineas2007randomized}
\begin{barticle}
\bauthor{\bsnm{Drineas}, \binits{P.}},
\bauthor{\bsnm{Mahoney}, \binits{M.W.}}:
\batitle{A randomized algorithm for a tensor-based generalization of the
  singular value decomposition}.
\bjtitle{Linear algebra and its applications}
\bvolume{420}(\bissue{2-3}),
\bfpage{553}--\blpage{571}
(\byear{2007})
\end{barticle}
\endbibitem

\bibitem{friedland2011fast}
\begin{botherref}
\oauthor{\bsnm{Friedland}, \binits{S.}},
\oauthor{\bsnm{Mehrmann}, \binits{V.}},
\oauthor{\bsnm{Miedlar}, \binits{A.}},
\oauthor{\bsnm{Nkengla}, \binits{M.}}:
Fast low rank approximations of matrices and tensors.
The Electronic Journal of Linear Algebra
\textbf{22}
(2011)
\end{botherref}
\endbibitem

\end{thebibliography}
\end{document}